\documentclass[11pt]{article}
\usepackage{amssymb,epsf,latexsym, amsmath}
\newcommand{\bfi}{\bfseries\itshape}

\makeatletter
\@addtoreset{figure}{section}
\def\thefigure{\thesection.\@arabic\c@figure}
\def\fps@figure{h, t}
\@addtoreset{table}{bsection}
\def\thetable{\thesection.\@arabic\c@table}
\def\fps@table{h, t}
\@addtoreset{equation}{section}

\makeatother

\newtheorem{thm}{Theorem}[section]
\newtheorem{lemma}{Lemma}[section]
\newtheorem{cor}{Corollary}[section]
\newtheorem{prop}{Proposition}[section]
\newtheorem{defn}{Definition}[section]
\allowdisplaybreaks
\def\intprod{\mathbin{\hbox to 6pt{%
                 \vrule height0.4pt width5pt depth0pt
                 \kern-.4pt
                 \vrule height6pt width0.4pt depth0pt\hss}}}

\begin{document}

\title{Multisymplectic Geometry, Variational Integrators,
and Nonlinear PDEs}

\author{
Jerrold E. Marsden
\\Control and Dynamical Systems
\\California Institute of Technology 107-81
\\Pasadena, CA 91125 USA
\\{\footnotesize email: marsden@cds.caltech.edu}
\and George W. Patrick
\\Department of Mathematics and Statistics
\\University of Saskatchewan
\\Saskatoon, SK S7N5E6 CANADA
\\{\footnotesize email: patrick@math.usask.ca}
\and Steve Shkoller
\\Center for Nonlinear Studies, MS-B258
\\Los Alamos National Laboratory
\\Los Alamos, NM 87545 USA
\\{\footnotesize email: shkoller@cds.caltech.edu}
}
\date{February 1997; this version, May 15, 1998\\
{\small To appear in {\it Communications in Mathematical Physics}}}

\maketitle

\begin{abstract}

This paper presents a geometric-variational approach to
continuous and discrete mechanics and field theories.  Using
multisymplectic geometry, we show that the existence of the
fundamental geometric structures as well as their preservation
along solutions can be obtained directly from the variational
principle.  In particular, we prove that a unique
multisymplectic structure is obtained by taking the derivative
of an action function, and use this structure to prove
covariant generalizations of conservation of symplecticity and
Noether's theorem.  Natural discretization schemes for PDEs,
which have these important preservation properties, then follow
by choosing a discrete action functional.  In the case of
mechanics, we recover the variational symplectic integrators of
Veselov type, while for PDEs we obtain covariant spacetime
integrators which conserve the corresponding discrete
multisymplectic form as well as the discrete momentum mappings
corresponding to symmetries.  We show that the usual notion of
symplecticity along an infinite-dimensional space of fields can
be naturally obtained by making a spacetime split.  All of the
aspects of our method are demonstrated with a nonlinear
sine-Gordon equation, including computational results and a
comparison with other discretization schemes.
\end{abstract}

\tableofcontents

\section{Introduction} 
The purpose of this paper is to develop the geometric
foundations for multi\-sym\-plectic\--momen\-tum integrators for
variational partial differential equations (PDEs).  These
integrators are the PDE generalizations of symplectic
integrators that are popular for Hamiltonian ODEs (see, for
example, the articles in Marsden, Patrick and Shadwick [1996],
and especially the review article of McLachlan and Scovel
[1996]) in that they are covariant spacetime integrators which
preserve the geometric structures of the system.

Because of the covariance of our method which we shall describe
below, the resulting integrators are spacetime localizable in
the context of hyperbolic PDEs, and generalize the notion of
symplecticity and symmetry preservation in the context of
elliptic problems.  Herein, we shall primarily focus on
spacetime integrators; however, we shall remark on the
connection of our method with the finite element method for
elliptic problems, as well as the Gregory and Lin [1991] method
in optimal control.

Historically, in the setting of ODEs, there have been many
approaches devised for constructing symplectic integrators,
beginning with the original derivations based on generating
functions (see de Vogelaere [1956]) and proceeding to
symplectic Runge-Kutta algorithms, the shake algorithm, and
many others. In fact, in many areas of molecular dynamics,
symplectic integrators such as the Verlet algorithm and
variants thereof are quite popular, as are symplectic
integrators for the integration of the solar system.  In these
domains, integrators that are either symplectic or which are
adaptations of symplectic integrators, are amongst the most
widely used.

A fundamentally new approach to symplectic integration is that
of Veselov [1988], [1991] who developed a discrete mechanics
based on a discretization of Hamilton's principle.  Their
method leads in a natural way to symplectic-momentum
integrators which include the shake and Verlet integrators as
special cases (see Wendlandt and Marsden [1997]). In addition,
Veselov integrators often have amazing properties with regard
to preservation of integrable structures, as has been shown by
Moser and Veselov [1991]. This aspect has yet to be exploited
numerically, but it seems to be quite important.

The approach we take in this paper is to develop a Veselov-type
discretization for PDE's in variational form. The relevant
geometry for this situation is multisymplectic geometry (see
Gotay, Isenberg, and Marsden [1997] and Marsden and Shkoller
[1997]) and we develop it in a variational  framework. As we
have mentioned, this naturally leads to multisymplectic-momentum
integrators.  It is well-known that such integrators
cannot in general preserve the Hamiltonian {\it exactly} (Ge and
Marsden [1988]). However, these integrators have, under
appropriate circumstances, very good energy performance in the sense of the
conservation of a nearby Hamiltonian up to exponentially small errors,
assuming
small time steps, due to a result of Neishtadt [1984]. See
also Dragt and Finn [1979], and Simo and Gonzales [1993]. This is
related to backward error analysis; see Sanz-Serna and Calvo
[1994], Calvo and Hairer [1995], and the recent work of Hyman, Newman and
coworkers and references therein. It would be quite interesting to
develop the links with Neishtadt's analysis more thoroughly.

An important part of our approach is to understand how the
symplectic nature of the integrators is implied by the
variational structure.  In this way we are able to identify the
symplectic and momentum conserving properties after discretizing
the variational principle itself.  Inspired by a paper of Wald
[1993],  we obtain a formal method for locating the symplectic
or multisymplectic  structures directly from the action function
and its derivatives. We present the method  in  the context of
ordinary Lagrangian  mechanics, and apply it to discrete Lagrangian mechanics,
and both continuous and discrete multisymplectic field theory.
While in these contexts our variational method merely uncovers 
the  well-known differential-geometric structures, our method
forms an excellent pedagogical approach to those theories. 

\paragraph{Outline of paper.}
\begin{itemize}
\item[\S2] In this section we sketch the three main aspects
of our variational approach in the familiar context of particle
mechanics. We show that the usual symplectic $2$-form on the
tangent bundle of the configuration manifold arises naturally
as the boundary term in the first variational principle.  We
then show that application of $d^2=0$ to the variational
principle restricted to the space of solutions of the
Euler-Lagrange equations  produces the familiar concept of
conservation of the symplectic form; this statement is obtained
variationally in a non-dynamic context; that is, we do not
require an evolutionary flow.  We then show that if the action
function is left invariant by a symmetry group, then Noether's
theorem follows directly and simply from the variational
principle as well. 

\item[\S3] Here we use our variational approach to construct
discretization schemes for mechanics which preserve the
discrete symplectic form and the associated discrete momentum
mappings.  

\item[\S4] This section defines the three aspects of our
variational approach in the multisymplectic field-theoretic
setting. Unlike the traditional approach of defining the
canonical multisymplectic form on the dual of the first jet
bundle and then pulling back to the Lagrangian side using the
covariant Legendre transform, we obtain the geometric structure
by staying entirely on the Lagrangian side.  We prove the
covariant analogue of the fact that the flow of conservative
systems consists of symplectic maps; we call this result the
{\it multisymplectic form formula}. After variationally proving
a covariant version of Noether's theorem, we show that one can
use the multisymplectic form formula to recover the usual
notion of symplecticity of the flow in an infinite-dimensional
space of fields by making a spacetime split.  We demonstrate
this machinery using a nonlinear wave equation as an example.

\item[\S5] In this section we develop discrete field
theories from which the covariant integrators follow.  We
define discrete analogues of the first jet bundle of the
configuration bundle whose sections are the fields of interest,
and proceed to define the discrete action sum. We then apply
our variational algorithm to this discrete action function to
produce the discrete Euler-Lagrange equations and the discrete
multisymplectic forms.  As a consequence of our methodology, we
show that the solutions of the discrete Euler-Lagrange
equations satisfy the discrete version of the multisymplectic
form formula as well as the discrete version of  our
generalized Noether's theorem. Using our  nonlinear wave
equation example, we develop various multisymplectic-momentum
integrators for the sine-Gordon equations, and compare our
resulting numerical scheme with the energy-conserving methods of
Li and Vu-Quoc [1995] 
and Guo, Pascual, Rodriguez, and Vazquez [1986].  
Results are presented for  long-time
simulations of kink-antikink solutions for over 5000 soliton
collisions.

\item[\S6] This section contains some important remarks
concerning the variational integrator methodology. For example,
we discuss integrators for reduced systems, the role of grid
uniformity, and the interesting connections with the finite-element
methods for elliptic problems.  We also make some comments on future work.
\end{itemize}

\section{Lagrangian Mechanics}
\label{mechanics}
\paragraph{Hamilton's Principle.}
We begin by recalling a problem going back to Euler, Lagrange
and Hamilton in the period 1740--1830. Consider an
$n$-dimensional configuration manifold $Q$ with its tangent
bundle $TQ$. We denote coordinates on $Q$ by $q^i$ and those on
$TQ$ by $(q^i, \dot{q}^i)$. Consider a  Lagrangian $L : TQ
\rightarrow \mathbb{R}$. Construct the  corresponding action
functional $S$ on $C^2$ curves $ q (t) $ in $Q$ by integration
of $L$ along the tangent to the curve. In coordinate notation,
this reads
\begin{equation}
S\bigl(q(t)\bigr)
\equiv\int_a^bL\left(q^i(t),\frac{dq^i}{dt}(t)\right)\,dt.
\label{1}
\end{equation}
The action functional depends on $a$ and $b$, but this is not
explicit in the notation. {\bfi Hamilton's principle} seeks the
curves $q(t)$ for which the functional $S$ is stationary under
variations of $q(t)$ with fixed endpoints; namely, we seek
curves $q(t)$ which satisfy
\begin{equation}
dS\bigl(q(t)\bigr) \cdot \delta q(t) \equiv
{\displaystyle \left. \frac{d}{d\epsilon }
\right|_{\epsilon =0}} S
\bigl(q_\epsilon(t)\bigr)=0
\label{1b}
\end{equation}
for all
$\delta q(t) \text{ with } \delta q(a) = \delta q(b) = 0$,  where
$q_\epsilon$ is a smooth family of curves with $q_0=q$ and
$(d/d\epsilon)|_{\epsilon=0}q_\epsilon = \delta q$.
Using integration by parts, the calculation for this is
simply
\begin{eqnarray}
dS\bigl(q(t)\bigr)\cdot \delta q(t)&=&
{\displaystyle \left.
\frac{d}{d\epsilon }\right|_{\epsilon =0}}
  \int_a^bL\left(q_\epsilon^i(t),
  \frac{dq_\epsilon^i}{dt}(t)\right)\,dt
  \nonumber\\
&=&\int_a^b\delta q^i\left(
\frac{\partial L}{\partial q^i}-\frac d{dt}
  \frac{\partial L}{\partial\dot q^i}\right)dt+
  \left.\frac{\partial L}{\partial\dot q^i}
  \delta q^i\right|_a^b.
\label{2}
\end{eqnarray}
The last term in (\ref{2}) vanishes since $\delta q(a)=\delta
q(b)=0$, so that the requirement (\ref{1b}) for $S$ to be
stationary yields the {\bfi Euler-Lagrange equations}
\begin{equation}
\frac{\partial L}{\partial q^i}-\frac d{dt}
  \frac{\partial L}{\partial\dot q^i}=0.\label{3}
\end{equation}
Recall that $L$ is called {\bfi regular} when the symmetric
matrix
$[\partial^2L/\partial\dot  q^i\partial\dot q^j]$ is everywhere 
nonsingular. If $L$ is regular, the Euler-Lagrange equations
are second order ordinary differential  equations
for the required curves.

\paragraph{The standard geometric setting.} The
action~(\ref{1}) is independent of the choice of coordinates,
and thus the Euler-Lagrange equations are coordinate
independent as well. Consequently, it is natural that the
Euler-Lagrange equations may be intrinsically  expressed using
the language of differential geometry. This intrinsic
development of mechanics is now  standard, and can be seen, for
example, in Arnold [1978], Abraham and Marsden [1978], and
Marsden and Ratiu [1994].

The {\bfi canonical $1$-form} $\theta_0$ on the $2n$-dimensional
cotangent bundle of $Q$, $T^*Q$ is defined by
\begin{equation}
\theta_0(\alpha_q)w_{\alpha_q}\equiv
\alpha_q\cdot T\pi_Qw_{\alpha_q},
  \qquad\mbox{$\alpha_q\in T^*_qQ$,
  $w_{\alpha_q}\in T_{\alpha_q}T^*Q$},
\nonumber
\end{equation}
where $\pi_Q:T^*Q\rightarrow Q$ is the canonical projection. The
Lagrangian $L$ intrinsically defines a fiber preserving bundle
map $\mathbb{F}L:TQ\rightarrow T^*Q$, the {\bfi Legendre 
transformation}, by vertical differentiation:
\begin{equation}
\mathbb{F}L(v_q)w_q\equiv
{\displaystyle \left. \frac{d}{d\epsilon }
\right|_{\epsilon =0}}L(v_q+\epsilon
w_q).
\nonumber
\end{equation}
We define the {\bfi Lagrange $1$-form} on $TQ$, the Lagrangian
side, by pull-back $\theta_L\equiv \mathbb{F} L^*\theta_0$,
and the {\bfi Lagrange $2$-form} by $\omega_L=-d\theta_L$. We
then seek a vector field $X_E$ (called the {\bfi Lagrange vector
field}) on $TQ$ such that $X_E\intprod\omega_L=dE$, where the
{\bfi energy} $E$ is defined by $E(v_q)\equiv \mathbb{F}
L(v_q)v_q-L(v_q)$.

If $\mathbb{F}L$ is a local diffeomorphism then $X_E$ exists
and is unique, and its integral curves solve the Euler-Lagrange
equations~(\ref{3}). In addition, the flow $F_t$ of $X_E$
preserves $\omega_L$; that is, $F_t^*\omega_L= \omega_L$. Such
maps are {\bfi symplectic}, and the form $\omega_L$ is called a
{\bfi symplectic $2$-form}.  This is an example of a {\bfi
symplectic manifold}:  a pair $(M, \omega)$ where $M$ is a
manifold and $\omega$ is closed nondegenerate $2$-form.

Despite the compactness and precision of this
differential-geometric approach, it is difficult to motivate
and, furthermore, is not entirely contained on the Lagrangian
side. The canonical $1$-form~$\theta_0$ seems to appear from
nowhere,  as does the Legendre transform~$\mathbb{F}L$.
Historically, after the Lagrangian picture on $TQ$ was
constructed, the canonical picture  on $T^*Q$ emerged through
the work of Hamilton, but the modern approach described above
treats the relation between the Hamiltonian and Lagrangian
pictures of mechanics as a mathematical tautology, rather than
what it is---a discovery of the highest order. 

\paragraph{The variational approach.} More and more, one is
finding that there are advantages to staying on the
``Lagrangian side''. Many examples can be given, but the
theory of Lagrangian reduction (the Euler-Poincar\'e equations
being an instance) is one example (see, for example,
 Marsden and Ratiu [1994] and Holm, Marsden and Ratiu
[1998a,b]); another, of many, is the direct variational approach
to questions in black hole dynamics given by Wald [1993]. In
such studies, it is the variational principle that is the
center of attention.

We next show that one can derive in a natural way the
fundamental differential geometric structures, including
momentum mappings, directly from the variational approach. This
development begins by removing the boundary condition $\delta
q(a)=\delta q(b)=0$ from~(\ref{2}). Equation~(\ref{2}) becomes
\begin{equation}
dS\bigl(q(t)\bigr)\cdot \delta q(t)
=\int_a^b\delta q^i\left(
\frac{\partial L}{\partial q^i}-\frac d{dt}
  \frac{\partial L}{\partial\dot q^i}\right)dt+
  \left.\frac{\partial L}{\partial\dot q^i}
  \delta q^i\right|_a^b,
\label{5}
\end{equation}
where the left side now operates on more general $\delta q$ (this
generalization will be described in detail in Section (4)),
while the last term on the right side does not vanish.  That
last term of~(\ref{5}) is a linear pairing of the
function~$\partial L/\partial\dot q^i$, a function of $q^i$
and $\dot q^i$, with the tangent vector $\delta q^i$. Thus,
one may consider it to be a $1$-form on $TQ$; namely the $1$-form
$(\partial L/\partial\dot q^i)dq^i$. This is exactly the
Lagrange $1$-form, and we can turn this into a formal
theorem/definition:
\begin{thm}
\label{thm2.1}
Given a $C^k$ Lagrangian $L$, $k\ge2$, there exists a
unique $C^{k-2}$ mapping $D_{\mbox{\scriptsize EL}}L:\ddot
Q\rightarrow T^*Q$, defined on the second order submanifold
\begin{displaymath}
\ddot Q\equiv\left\{ \left. \frac{d^2q}{dt^2}(0) \; \right| \;
\mbox{$q$ a $C^2$ curve in $Q$} \right\}
\end{displaymath}
of $TTQ$, and a unique $C^{k-1}$ $1$-form $\theta_L$ on $TQ$,
such that, for all~$C^2$ variations $q_\epsilon(t)$,
\begin{equation}\label{7}
dS\bigl(q(t)\bigr)\cdot \delta q(t)
=\int_a^b D_{\mbox{\scriptsize EL}}L\left(
\frac{d^2q}{dt^2}\right)
  \cdot\delta q\,dt+
  \left.\theta_L\left(\frac{dq}{dt}\right)
  \cdot\hat{\delta q}\right|_a^b,
\end{equation}
where
\begin{displaymath}\delta q(t)\equiv
{\displaystyle \left. \frac{d}{d\epsilon }
\right|_{\epsilon = 0}}
q_\epsilon(t),\qquad
\hat{\delta q}(t)\equiv
{\displaystyle \left. \frac{d}{d\epsilon }
\right|_{\epsilon = 0}}
{\displaystyle \left. \frac{d}{dt }\right|_{t = 0}}
q_\epsilon(t).
\end{displaymath}
The $1$-form so defined is called the {\bfi Lagrange
$1$-form.}
\end{thm}

Indeed, uniqueness and local existence follow from the
calculation~(\ref{2}), and the coordinate independence of the
action, and then global existence is immediate. Here then, is
the first aspect of our method:
\begin{quote}\it Using the variational principle, the
Lagrange~$1$-form $\theta_L$ is the ``boundary part''  of the
the  functional derivative of the action when the boundary is
varied.  The analogue  of the symplectic form is the (negative
of)  the exterior derivative of $\theta_L$.
\end{quote} For the mechanics example being discussed, we
imagine a development wherein $\theta_L$ is so defined and we
define $\omega_L\equiv-d\theta_L$.

\paragraph{Lagrangian flows are symplectic.} One of Lagrange's
basic discoveries was that the solutions of the Euler-Lagrange
equations give rise to a symplectic map. It is a curious twist
of history that he did this without the machinery of either
differential forms, of the Hamiltonian formalism or of
Hamilton's principle itself. (See Marsden and Ratiu [1994] for
an account of some of this history.)

Assuming that $L$ is regular, the variational principle then
gives coordinate independent second order ordinary
differential equations, as we have noted. We temporarily
denote the vector field on $TQ$ so obtained by $X$, and its
flow by $F_t$. Our further development relies on a change of
viewpoint: we focus on the restriction of $S$ to the subspace
$\mathcal{C}_L$ of  solutions of the variational principle. The
space $\mathcal{C}_L$ may be identified with the initial
conditions, elements of $TQ$, for the flow: to $v_q \in T
Q$, we associate the integral curve $s\mapsto F_s(v_q)$,
$s\in[0,t]$. The value of $S$ on that curve is denoted by
$S_t$, and again called the {\bfi action}. Thus, we define the
map $S_t : TQ \rightarrow \mathbb{R}$ by
\begin{equation} \label{actiont.equation} 
       S_t (v_q) = \int _0 ^t L(q(s), \dot{q}(s)) \,ds, 
\end{equation}  
where $(q(s), \dot{q}(s)) = F_s (v _q) $. The fundamental
equation~(\ref{7}) becomes
\begin{displaymath}
dS_t(v_q)w_{v_q}=\theta_L\bigl(F_t(v_q)\bigr) \cdot
{\displaystyle \left. 
\frac{d}{d\epsilon }\right|_{\epsilon = 0}}
F_t(v_q^\epsilon ) -\theta_L(v_q) \cdot w_{v_q},
\end{displaymath}
where $\epsilon\mapsto v_q^\epsilon$ is an arbitrary curve in
$TQ$ such that $v_q^0=v_q$ and $(d/d\epsilon)|_0 v^\epsilon_q =
w_{v_q}$.  We have thus derived the equation
\begin{equation}\label{8} 
     dS_t=F_t^*\theta_L-\theta_L.
\end{equation}
Taking the exterior derivative of~(\ref{8}) 
yields the fundamental
fact that the flow of $X$ is symplectic:
\begin{displaymath}
0=ddS_t=d(F_t^*\theta_L-\theta_L)=-F_t^*\omega_L+\omega_L
\end{displaymath}
which is equivalent to
\[
F_t^*\omega_L=\omega_L.
\]
This leads to the following:
\begin{quote}
\it Using the variational principle, the fact
that the evolution is  symplectic is a consequence of the
equation $d^2 =0$, applied to the action restricted to the space
of solutions of the variational principle.
\end{quote}
In passing, we note that~(\ref{8}) also provides the
differential-geometric equations for~$X$. Indeed, one
time derivative of~(\ref{8})  and using
(\ref{actiont.equation}) gives 
$dL = {\mathfrak L}_X\theta_\mathcal{L}$,
so that
\begin{displaymath}
X\intprod\omega_L=-X\intprod d\theta_L
=-{\mathfrak L}_X\theta_L+d(X\intprod\theta_L)
=d(X\intprod\theta_L-L)=dE,
\end{displaymath}
if we define $E\equiv X\intprod\theta_L-L$.
Thus, we quite naturally find that $X=X_E$.

Of course, this set up also leads directly to Hamilton-Jacobi
theory, which was one of the ways in which symplectic
integrators were developed (see McLachlan and Scovel [1996]
and references therein.) However, we shall not pursue the
Hamilton-Jacobi aspect of the theory here.

\paragraph{Momentum maps.} Suppose that a Lie group~$G$, with
Lie algebra~$\mathfrak{g}$, acts on $Q$, and hence on curves
in~$Q$, in such a way that the action~$S$ is invariant.
Clearly, $G$ leaves the set of solutions of the variational
principle invariant,  so the action of~$G$
restricts to~$\mathcal{C}_L$, and the group action commutes
with~$F_t$. Denoting the infinitesimal generator of
$\xi\in\mathfrak{g}$ on $TQ$ by $\xi_{TQ}$, we have
by~(\ref{8}),
\begin{equation}\label{9}
0 = \xi_{TQ}\intprod dS_t
  = \xi_{TQ}\intprod (F_t^*\theta_L-\theta_L)
  = F_t^*(\xi_{TQ}\intprod\theta_L) - \xi_{TQ}\intprod\theta_L.
\end{equation}
For $\xi\in\mathfrak{g}$, define 
$J_\xi:TQ\rightarrow\mathbb R$ by
$J_\xi\equiv\xi_{TQ}\intprod\theta_L$. Then~(\ref{9}) says
that $J_\xi$ is an integral of the flow of $X_E$. 
We have arrived at a version of
Noether's theorem 
(rather close to the original derivation of Noether):
\begin{quote}\it Using the variational principle, Noether's
theorem results from the infinitesimal invariance of the action
restricted to space of solutions of the variational principle.
The conserved momentum associated to a Lie algebra element
$\xi$ is $J_\xi =  \xi_{TQ}\intprod\theta_L$, where
$\theta_L$ is the Lagrange one-form.
\end{quote}

\paragraph{Reformulation in terms of first variations.}
We have just seen that symplecticity of the flow and
Noether's theorem result from restricting the action to the
space of solutions. One tacit assumption is that the space of
solutions is a manifold in some appropriate sense. This is a
potential problem, since solution spaces for field theories
are  known to have singularities (see, e.g., Arms, Marsden and
Moncrief [1982]). More seriously there is the problem of finding
a  multisymplectic  analogue of the statement  that the
Lagrangian  flow map is symplectic, since for multisymplectic
field theory one obtains an evolution picture only after
splitting  spacetime into space and time and adopting the
``function space'' point of view. Having   the general
formalism  depend either on a spacetime split or  an analysis
of the associated   Cauchy problem would be contrary to the
general thrust of this article. We now give a formal argument,
in the context of Lagrangian mechanics, which shows how both
these problems can be  simultaneously avoided.

Given a solution $q(t) \in \mathcal{C}_L$, a {\bfi first
variation} at $q(t)$ is  a vector field $V$ on $Q$ such that
$t\mapsto F^V_\epsilon\circ q(t)$ is also a solution curve
(i.e. a curve in $\mathcal{C}_L$). We think of the solution
space $\mathcal{C}_L$ as being a (possibly) singular subset of
the smooth space of all putative curves $\mathcal{C}$ in $TQ$,
and the restriction of $V$ to $q(t)$ as being the derivative of
some curve in $\mathcal{C}_L$ at $q(t)$. When $\mathcal{C}_L$
is a manifold, a first variation is a vector at $q(t)$
tangent to $\mathcal{C}_L$. Temporarily define $\alpha\equiv
dS-\theta_L$ where by abuse of notation $\theta_L$ is the one
form on $\mathcal{C}$ defined by
\begin{displaymath}
\theta_L\bigl(q(t)\bigr)\delta q(t)\equiv 
\theta_L(b)\delta q(b)-
\theta_L(a)\delta q(a).\end{displaymath}
Then $\mathcal{C}_L$
is defined by $\alpha=0$ and we have the equation
\begin{displaymath}dS=\alpha+\theta_L,\end{displaymath}
so if $V$ and $W$ are first variations at $q(t)$, we obtain
\begin{equation}\label{54}
0=V\intprod W\intprod d^2S=V\intprod W\intprod d\alpha
+V\intprod W\intprod
d\theta_L.
\end{equation}
We have the identity
\begin{equation}\label{50}
d\alpha(V,W)\bigl(q(t)\bigr)
  =V\bigl(\alpha(W)\bigr)-
  W\bigl(\alpha(V)\bigr)-\alpha([V,W]),
\end{equation}
which we will use to evaluate~(\ref{54}) at the curve 
$q(t)$.
Let $F^V_\epsilon$ denote the flow of $V$,  define
$q_\epsilon^V(t)\equiv F_\epsilon^V\bigl(q(t)\bigr)$,
and make similar definitions with $W$ replacing $V$.
For the first term of~(\ref{50}), we have
\begin{displaymath}
V\bigl(\alpha(W)\bigr)\bigl(q(t)\bigr)=
{\displaystyle \left. \frac{d}{d\epsilon }
\right|_{\epsilon = 0}}
\alpha(W)(
  q^V_\epsilon),\end{displaymath}
which vanishes, since $\alpha$ is zero along
$q^V_\epsilon$ for every $\epsilon$. 
Similarly the second term of~(\ref{50})
at $q(t)$ also
vanishes,  while  the third term vanishes since
$\alpha\bigl(q(t)\bigr)=0$.
Consequently, symplecticity of the Lagrangian flow $F_t$
may be written
\begin{displaymath}
V\intprod W\intprod d\theta_L=0,
\end{displaymath}
for all first variations $V$ and $W$. This formulation is valid
whether or not the solution space is a manifold, and it does not
explicitly refer to any temporal notion.
Similarly, Noether's theorem may  be written in this way.
Summarizing,

\begin{quote}\it Using the variational principle, the analogue
of the evolution is  symplectic is the equation $d^2S =0$
restricted to first variations of the space of solutions of
the variational principle. The analogue of Noether's theorem
is infinitesimal invariance of $dS$ restricted to first
variations  of the space of solutions of the variational
principle.
\end{quote}

The variational route to the differential-geometric formalism
has obvious pedagogical advantages. More than that, however,
it systematizes searching for the corresponding formalism in
other contexts. We shall in the next sections show how this
works in the context of discrete mechanics, classical field
theory and multisymplectic geometry.

\section{Veselov Discretizations of Mechanics}
The discrete Lagrangian formalism in Veselov~[1988], [1991]
fits nicely into our variational framework. Veselov
uses $Q\times Q$ for the discrete version of the tangent bundle
of a configuration space~$Q$; heuristically, given some a
priori choice of time interval
$\Delta t$, a point
$(q_1,q_0)\in Q\times Q$ corresponds to the tangent vector
$(q_1-q_0)/\Delta t$. Define a {\bfi discrete
Lagrangian\/} to be a smooth map
$L:Q\times Q=\{q_1,q_0\}\rightarrow\mathbb R$,
and the corresponding action to be
\begin{equation}\label{850}
S\equiv\sum_{k=1}^{n}L(q_k,q_{k-1}).
\end{equation}
The variational principle is to extremize~$S$ for variations
holding the endpoints $q_0$ and $q_n$
fixed. This variational principle determines a ``discrete flow''
$F:Q\times Q\rightarrow Q\times Q$ by 
$F(q_1,q_0)=(q_2,q_1)$, where
$q_2$ is found from the 
{\bfi discrete Euler-Lagrange equations}
(DEL equations):
\begin{equation}\label{40}
\frac{\partial L}{\partial q_1}(q_1,q_0)
+\frac{\partial L}{\partial q_0}(q_2,q_1)=0.
\end{equation}
In this section we work out the basic differential-geometric
objects of this discrete mechanics directly from the
variational point of view, consistent with our philosophy in
the last section.

A mathematically significant aspect of this theory is how it
relates to integrable systems, a point taken up by Moser and
Veselov~[1991]. We will not explore this aspect in any detail
in this paper, although later, we will briefly discuss the
reduction process and we shall test an integrator for an
integrable pde, the sine-Gordon equation. 

\paragraph{The Lagrange $1$-form.}
We begin by calculating $dS$ for variations that do not fix the
endpoints:
\begin{eqnarray}
\lefteqn{dS(q_0,\cdots,q_n)\cdot(\delta q_0,\cdots,\delta q_n)}
\nonumber\\
&=&
  \sum_{k=0}^{n-1}\left(
  \frac{\partial L}{\partial q_1}(q_{k+1},q_k)\delta q_{k+1}
  +\frac{\partial L}{\partial q_0}(q_{k+1},q_k)\delta q_k
  \right)
\nonumber\\
&=&
  \sum_{k=1}^{n}
  \frac{\partial L}{\partial q_1}(q_{k},q_{k-1})\delta q_k
  +\sum_{k=0}^{n-1}
  \frac{\partial L}{\partial q_0}(q_{k+1},q_k)\delta q_k
\nonumber\\
&=&
  \sum_{k=1}^{n-1}\left(
    \frac{\partial L}{\partial q_1}(q_{k},q_{k-1})
    +\frac{\partial L}{\partial q_0}(q_{k+1},q_k)
  \right)\delta q_k \nonumber\\
&&\qquad\mbox{}
  +\frac{\partial  L}{\partial q_0}(q_1,q_0)\delta q_0
  +\frac{\partial  L}{\partial q_1}(q_n,q_{n-1})
  \delta q_n\label{12}.
\end{eqnarray}
It is the last two terms that arise from the boundary
variations~(i.e.~these are the ones 
that are zero if the boundary is fixed),
and so these are the terms amongst which we expect to find
the discrete analogue of the Lagrange $1$-form.
Actually, interpretation of the boundary terms gives
the  {\em two} $1$-forms on $Q\times Q$
\begin{equation}
 \theta_L^-(q_1,q_0)\cdot(\delta q_1,\delta q_0)\equiv
 \frac{\partial L}{\partial q_0}(q_1,q_0)\delta q_0,
\label{13}\end{equation}
and
\begin{equation}
  \theta_L^+(q_1,q_0)\cdot(\delta q_1,\delta q_0)\equiv
  \frac{\partial L}{\partial q_1}(q_1,q_0)\delta q_1,
\label{14}\end{equation}
and we regard {\em the pair} $(\theta^-,\theta^+)$ 
as being the analogue
of the one form in this situation.

\paragraph{Symplecticity of the flow.}
We
parameterize the solutions of the variational principle by the
initial conditions~$(q_1,q_0)$, and restrict  $S$ to that solution space.
Then Equation~(\ref{12}) becomes
\begin{equation}
dS=\theta_L^-+ F^*\theta_L^+.\label{15}
\end{equation}
We should be able to obtain the symplecticity 
of~$F$ by determining
what the equation $ddS=0$ means 
for the right-hand-side of~(\ref{15}).
At first, this does not appear to work, since~$ddS=0$ gives
\begin{equation}
 F^*(d\theta_L^+)=-d\theta_L^-,\label{16}
\end{equation}
which apparently says that $F$ pulls a certain  
$2$-form back to a different
$2$-form.
The situation is aided by the observation that, from~(\ref{13})
and~(\ref{14}),
\begin{equation}\label{20}
\theta_L^-+\theta_L^+=dL,
\end{equation}
and consequently,
\begin{equation}\label{2010}
d\theta_L^-+d\theta_L^+=0.
\end{equation}
Thus, there are {\em two\/} generally distinct $1$-forms, 
but (up to sign) only
{\em one\/} $2$-form. If we make the definition
\begin{displaymath}
\omega_L\equiv d\theta_L^-=-d\theta_L^+,
\end{displaymath}
then~(\ref{16}) becomes $F^*\omega_L=\omega_L$.
Equation~(\ref{13}), in coordinates, gives
\begin{displaymath}
\omega_L=\frac{\partial^2L}{\partial q_0^i\partial q_1^j}
dq_0^i\wedge dq_1^j,
\end{displaymath}
which agrees with the discrete symplectic form found in 
Veselov~[1988],~[1991].

\paragraph{Noether's Theorem.} Suppose a Lie group $G$ with
Lie algebra $\mathfrak{g}$ acts on $Q$,  and hence diagonally on
$Q\times Q$, and that $L$ is $G$-invariant. Clearly, $S$ is
also $G$-invariant and $G$ sends critical points of $S$ to
themselves. Thus, the action of $G$ restricts to the space of
solutions, the map $F$ is $G$-equivariant,
and  from~(\ref{15}),
\begin{displaymath}
0=\xi_{Q\times Q}\intprod dS=\xi_{Q\times Q}\intprod\theta_L^-+
  \xi_{Q\times Q}\intprod(F^*\theta_L^+),
\end{displaymath}
for $\xi \in \mathfrak{g}$, or equivalently, using the
equivariance of~$F$,
\begin{equation}\label{21}
\xi_{Q\times Q}\intprod\theta_L^-=-F^*(\xi_{Q\times Q}
\intprod\theta^+).
\end{equation}
Since $L$ is $G$-invariant,~(\ref{20}) gives
$\xi_{Q\times Q}\intprod\theta_L^-=-\xi_{Q\times Q}
\intprod\theta_L^+$,
which in
turn converts~(\ref{21}) to the conservation equation
\begin{equation}\label{24}
\xi_{Q\times Q}\intprod\theta_L^+=F^*(\xi_{Q\times Q}
\intprod\theta^+).
\end{equation}
Defining the discrete momentum to be
\begin{displaymath}J_\xi\equiv\xi_{Q\times Q}\intprod\theta_L^+,
\end{displaymath}
we see that~(\ref{24}) becomes conservation of momentum.
A virtually identical derivation of this discrete 
Noether theorem is found in Marsden and~Wendlant~[1997].

\paragraph{Reduction.} As we mentioned above, this formalism
lends itself to a discrete version of the theory of Lagrangian
reduction (see Marsden and Scheurle [1993a,b], Holm, Marsden and
Ratiu [1998a] and Cendra, Marsden and Ratiu [1998]). This
theory is not the focus of this article, so we shall defer a
brief discussion of it until the conclusions.

\section{Variational Principles for Classical Field Theory}
\label{sec4}
\paragraph{Multisymplectic geometry.}

We now review some aspects of multisymplectic geometry, 
following
Gotay, Isenberg and Marsden [1997] and 
Marsden and Shkoller [1997].

We let $\pi_{XY}:Y \rightarrow X$ be a fiber bundle over an
oriented manifold $X$. Denote the first jet bundle over $Y$ by
$J^1(Y)$ or $J^1Y$ and identify it with the {\it affine}
bundle over $Y$ whose fiber over $y\in Y_x:=\pi_{XY}^{-1}(x)$
consists of Aff$(T_xX, T_yY)$, those linear mappings
$\gamma:T_xX \rightarrow T_y Y$ satisfying
$$ 
   T \pi_{XY} \circ \gamma = \text{Identity on } T_xX. 
$$

We let ${\rm dim}X = n+1$ and the fiber dimension of $Y$ be
$N$. Coordinates on $X$ are denoted $x^\mu, \mu = 1, 2, \dots
, n, 0$, and fiber coordinates on $Y$ are denoted by $y^A, A =
1, \dots , N$. These induce coordinates $v^A{}_\mu$ on the
fibers of $J^1(Y)$.  If
$\phi : X \rightarrow Y$ is a section of $\pi_{XY}$, its
tangent map at $x \in X$, denoted $T_x\phi$, is an element of
$J^1(Y)_{\phi(x)}$. Thus, the map $x \mapsto T_x\phi$ defines a
section of $J^1(Y)$ regarded as a bundle over $X$.  This
section is denoted $j^1(\phi)$ or $j^1\phi$ and is called the first jet of
$\phi$.  In coordinates, $j^1(\phi)$ is given by
\begin{equation}
x^\mu \mapsto (x^\mu, \phi^A(x^\mu), 
\partial_\nu \phi^A (x^\mu)),
\label{A4}
\end{equation}
where $\partial_\nu = \partial / \partial x^\nu$.

Higher order jet bundles of $Y$, $J^m(Y)$, then follow as $J^1(
\cdot \cdot \cdot (J^1(Y))$.   Analogous to the tangent map of
the projection
$\pi_{Y,J^1(Y)}$, $T\pi_{Y,J^1(Y)}:TJ^1(Y)\rightarrow TY$, we
may define the jet map of this projection which takes $J^2(Y)$
onto $J^1(Y)$
\begin{defn}
Let $\gamma \in J^1(Y)$ so that $\pi_{X,J^1(Y)}(\gamma)=x$. Then
$$ J\pi_{Y,J^1(Y)}:\text{\rm Aff}(T_xX, T_\gamma J^1(Y))
\rightarrow \text{\rm Aff}(T_xX, T\pi_{Y,J^1(Y)} \cdot T_\gamma
J^1(Y)).$$ We define the subbundle $Y''$ of $J^2(Y)$ over $X$
which consists of second-order jets so that on each fiber
$$Y''_x = \{ s\in J^2(Y)_\gamma \ | \
J\pi_{Y,J^1(Y)}(s)=\gamma\}.
\qquad$$
\end{defn}
In coordinates, if $\gamma \in J^1(Y)$ is given by
$(x^\mu, y^A, {v^A}_\mu)$, and $s\in J^2(Y)_\gamma$ is given by
$(x^\mu, y^A, {v^A}_\mu, {w^A}_\mu, {\kappa^A}_{\mu\nu})$, then
$s$ is a second-order jet if
${v^A}_\mu = {w^A}_\mu$.  Thus, the second jet of $\phi \in
\Gamma (\pi_{XY})$, $j^2(\phi)$, given in coordinates by the map
$x^\mu \mapsto (x^\mu, \phi^A, \partial _\nu \phi^A,
\partial _\mu \partial _\nu \phi^A)$, is an example of a
second-order jet.

\begin{defn}
The {\bfi dual jet bundle} $J^1(Y)^\star$ is the {\it vector}
bundle over $Y$ whose fiber at $y\in Y_x$ is the set of {\it
affine} maps from
$J^1(Y)_y$ to $\Lambda^{n+1}(X)_x$, the bundle of $(n+1)$-forms
on $X$. A smooth section of $J^1(Y)^\star$ is therefore an
affine bundle map of
$J^1(Y)$ to $\Lambda^{n+1}(X)$ covering $\pi_{XY}$.
\end{defn}
Fiber coordinates on $J^1(Y)^\star$ are $(p, p_A{}^\mu)$, which
correspond to the affine map given in coordinates by
\begin{equation} v^A{}_\mu \mapsto (p + p_A{}^\mu v^A{}_\mu)
d^{n+1} x, \label{B1}
\end{equation} where $d^{n+1}x = dx^1 \wedge \dots \wedge dx^n
\wedge dx^0$.

Analogous to the canonical one- and two-forms on a cotangent
bundle, there exist canonical ($n+1$)- and ($n+2$)-forms on the
dual jet bundle $J^1(Y)^\star$.  In coordinates, with $d^nx_\mu:=
\partial_\mu \intprod d^{n+1}x$, these forms are given by
\begin{equation}
\Theta = {p_A}^\mu dy^A \wedge d^nx_\mu + p d^{n+1}x \text{ 
and  }
\Omega = dy^A\wedge d{p_A}^\mu \wedge d^nx_\mu - dp \wedge
d^{n+1}x.
\label{B2}
\end{equation}

A Lagrangian density $\mathcal{L} : J^1(Y) \rightarrow
\Lambda^{n+1}(X)$ is a smooth bundle map over $X$.  In
coordinates, we write
\begin{equation}
\mathcal{L} (\gamma)  =  L(x^\mu, y^A, v^A{}_\mu) d^{n+1}x. \label{C1}
\end{equation}

The corresponding covariant Legendre transformation for
$\mathcal{L}$ is a fiber preserving map over $Y$, $\mathbb{F}\mathcal{L}: J^1(Y) \rightarrow J^1(Y)^\star$, expressed
intrinsically as the first order vertical Taylor approximation
to $\mathcal{L}$:
\begin{equation}
\mathbb{F}\mathcal{L}(\gamma) \cdot  \gamma'  
= \mathcal{L}(\gamma)
+\left.
\frac{d}{d\varepsilon} \right|_{\varepsilon = 0}
\mathcal{L}(\gamma +
\varepsilon(\gamma' - \gamma)) \label{C3}
\end{equation}
where $\gamma, \gamma' \in  J^1(Y)_y$.
A straightforward calculation shows that the covariant Legendre
transformation is given in coordinates by
\begin{equation}
p_A{}^\mu =
\frac{\partial L}{\partial v^A{}_\mu}, \text{ and }\quad p =
L -\frac{\partial L}{\partial v^A{}_\mu} v^A{}_\mu .\label{C2}
\end{equation}

We can then define the {\bfi Cartan form} as the
$(n+1)$-form $\Theta_{\mathcal{L}}$ on $J^1(Y)$  given by
\begin{equation}
\Theta_{\mathcal{L}} =  (\mathbb F\mathcal{L})^*\Theta,
\label{D1}
\end{equation}
and the $(n+2)$-form $\Omega _\mathcal{L}$ by
\begin{equation}
\Omega_\mathcal{L} = -d \Theta_\mathcal{L}
= (\mathbb F \mathcal{L})^* \Omega,
\label{D2}
\end{equation}
with local coordinate expressions
\begin{equation}
\begin{array}{c}
{\displaystyle
\Theta_\mathcal{L} = \frac{\partial L}{\partial {v^A}_\mu} dy^A
\wedge d^nx_\mu +\left( L -  \frac{\partial L}{\partial
{v^A}_\mu} {v^A}_\mu \right) d^{n+1}x,}  \\[12pt]
{\displaystyle \Omega_\mathcal{L} = dy^A \wedge d\left(
\frac{\partial L } {\partial {v^A}_\mu }\right) \wedge d^nx_\mu
- d\left[  L - \frac{\partial L}{\partial {v^A}_\mu } {v^A}_\mu
\right] \wedge d^{n+1}x. }
\end{array}
\label{D1a}
\end{equation}

This is the differential-geometric formulation of the
multisymplectic structure. Subsequently, we shall
show how we may obtain this structure directly from the
variational principle, staying entirely on the Lagrangian side 
$J^1(Y)$.

\paragraph{The multisymplectic form formula.} In this
subsection we prove a formula that is the multisymplectic
counterpart to the fact that in finite-dimensional mechanics,
the flow of a mechanical system consists of symplectic maps.
Again, we do this by studying the action function.
\begin{defn}
Let $U$ be a smooth manifold with (piecewise) smooth closed
boundary. We define the set of smooth maps
$$
  \mathcal{C}^\infty = \{ \phi: U \rightarrow Y \mid \pi_{XY}
  \circ \phi: U \rightarrow X \text{ is an embedding}\}.
$$ 
For each $\phi\in \mathcal{C}^\infty $, we set $\phi_X :=
\pi_{XY}\circ \phi$ and $U_X : = \pi_{XY}\circ \phi (U)$ so that
$\phi_X:U \rightarrow U_X$ is a diffeomorphism. $\qquad $.
\end{defn}
We may then define the infinite-dimensional manifold $\mathcal{C}$ to be the closure of $\mathcal{C}^\infty$ in either a
Hilbert space or Banach space topology. For example, the
manifold $\mathcal{C}$ may be given the topology of a Hilbert
manifold of bundle mappings,
$H^s(U,Y)$, ($U$  considered a bundle with fiber a point) for
any integer $s \geq  (n + 1)/2$, so that the Hilbert sections
${\phi \circ \phi_X^{-1}}$ in $Y$ are those whose distributional
derivatives up to order $s$ are square-integrable in any
chart.  With our condition on $s$, the Sobolev embedding
theorem makes such mappings well defined. Alternately, one may
wish to consider the Banach manifold
$\mathcal{C}$ as the closure of $\mathcal{C}^\infty$ in the
usual
$C^k$-norm, or more generally, in a Holder space
$C^{k+\alpha}$-norm. See Palais [1968] and Ebin and Marsden [1970] for a detailed
account of manifolds of mappings. The choice of topology for
$\mathcal{C}$ will not play a crucial role in this paper.

\begin{defn}
Let $\mathcal{G}$ be the Lie group of $\pi_{XY}$-bundle
automorphisms
$\eta_{Y}$ covering diffeomorphisms $\eta_{X}$, with Lie algebra
$\mathfrak{g}$.   We define the {\bfi action} 
$\Phi:\mathcal{G} \times \mathcal{C} \rightarrow\mathcal{C}$ by
$$\Phi({\eta_{Y}},\phi) = \eta_{Y}\circ ({\phi \circ
\phi_X^{-1}}) \circ
\eta_{X}^{-1}.$$ Furthermore, if $\phi \circ \phi_X^{-1} \in
\Gamma(\pi_{U_X,Y})$, then
$\Phi(\eta_{Y},\phi) \in \Gamma(\pi_{\eta_{X}(U_X),Y})$.
The {\bfi tangent space} to the manifold $\mathcal{C}$ at a point
$\phi$ is the set $ T_\phi \mathcal{C}$  defined by
\begin{equation}
\left\{  V \in C^\infty(X,TY) \mid \pi_{Y,TY}
\circ V = \phi, \& T \pi_{XY} \circ V = V_X, \text{ a
vector field on } X \right\} .
\end{equation}
\end{defn}
Of course, when these objects are topologized as we have
described, the definition of the tangent space becomes a
theorem, but as we have mentioned, this functional analytic
aspect plays a minor role in what follows.

Given vectors $V,W \in T_\phi\mathcal{C}$ we may extend them
to vector fields ${\mathcal V}, {\mathcal W}$ on
$\mathcal{C}$ by fixing vector fields $v,w \in TY$ such that
$V= v\circ ({\phi\circ\phi_X^{-1}})$ and
$W= w \circ ({\phi \circ \phi_X^{-1}})$, and letting ${\mathcal
V}_\rho = v \circ ({\rho \circ \rho_X^{-1}})$ and
${\mathcal W}_\rho= w \circ ({\rho \circ \rho_X^{-1}})$.  Thus,
the flow of ${\mathcal V}$ on $\mathcal{C}$ is given by
$\Phi(\eta_{Y}^\lambda,
\rho)$ where $\eta_{Y}^\lambda$ covering $\eta_{X}^\lambda$ is
the flow of $v$.  The definition  of the bracket of vector
fields using their flows, then shows that
$$ [{\mathcal V},{\mathcal W}](\rho) = [v,w] \circ (\rho
\circ \rho_X^{-1}) .$$ Whenever it is contextually clear, we
shall, for convenience, write $V$ for
$v \circ ({\phi \circ \phi_X^{-1}})$.

\begin{defn}
The {\bfi action function} $\mathcal{S}$ on $\mathcal{C}$ is defined as follows:
\begin{equation}
\mathcal{S}(\phi) = \int_{U_X}
\mathcal{L}({j^1(\phi \circ \phi_X^{-1}})) \; \mbox{for all} \;
 \phi \in \mathcal{C}. \qquad
\label{s3}
\end{equation}
\end{defn}
Let $\lambda \mapsto \eta_{Y}^\lambda$ 
be an arbitrary smooth path in
$\mathcal{G}$ such that $\eta_{Y}^0=e$, 
and let $V\in T_\phi \mathcal{C}$
be given by
\begin{equation}
V = \left.\frac{d}{d \lambda}\right|_{\lambda=0}
\Phi(\eta_{Y}^\lambda,
\phi), \text{ and } V_X =
\left.\frac{d}{d \lambda}\right|_{\lambda=0} 
\eta_{X}^\lambda \circ \phi.
\label{ig}
\end{equation}

\begin{defn}
We say that $\phi$ is a {\bfi stationary point, critical point, or extremum} 
of $\mathcal{S}$
if
\begin{equation}
\left.\frac{d}{d\lambda}\right|_{\lambda=0} \mathcal{S}
(\Phi( \eta_{Y}^\lambda, \phi)) =0. \qquad
\label{s4}
\end{equation}
\end{defn}
Then,
\begin{eqnarray}
d\mathcal{S}_\phi \cdot V &=&
\left.\frac{d}{d\lambda}\right|_{\lambda=0}
\int_{\eta_{X}^\lambda\circ\phi_X (U)}
\mathcal{L}(\Phi(\eta_{Y}^\lambda,\phi)) \label{s5} \\
&=&\left.\frac{d}{d\lambda}\right|_{\lambda=0}
\int_{\phi_X (U)} j^1({\phi \circ \phi_X^{-1}})^* 
j^1(\eta_{Y}^\lambda)^*
\Theta_\mathcal{L}, \nonumber
\end{eqnarray}
where we have used the fact that 
$\mathcal{L}(z)=z^*\Theta_\mathcal{L}$ 
for all holonomic sections $z$  of $J^1(Y)$
(see Corollary \ref{cor4.2} below), and that
$$j^1(\eta_{Y}\circ{\phi \circ \phi_X^{-1}}\circ\eta_{X}^{-1})=
j^1(\eta_{Y})\circ j^1({\phi \circ 
\phi_X^{-1}})\circ\eta_{X}^{-1}.$$

Using the Cartan formula, we obtain that
\begin{eqnarray}
d\mathcal{S}_\phi \cdot V &=&
\int_{U_X} j^1({\phi \circ \phi_X^{-1}})^* 
{\mathfrak L}_{j^1(V)}
\Theta_\mathcal{L} \nonumber \\
&=& \int_{U_X} j^1({\phi \circ \phi_X^{-1}})^* 
[{j^1(V)} \intprod
\Omega_\mathcal{L}] \nonumber \\
&& \qquad +
\int_{\partial U_X} j^1({\phi \circ 
\phi_X^{-1}})^* [{j^1(V)} \intprod
\Theta_\mathcal{L}] .
\label{s6}
\end{eqnarray}

Hence, a necessary condition for $\phi \in \mathcal{C}$ to be
an extremum of $\mathcal{S}$ is that the first term in
(\ref{s6}) vanish.  One may readily verify that the integrand
of the first term in (\ref{s6}) is equal to zero whenever
${j^1(V)}$ is replaced by $W\in TJ^1(Y)$ which is either
$\pi_{Y,J^1(Y)}$-vertical or tangent to $j^1({\phi \circ
\phi_X^{-1}})$ (see Marsden and Shkoller [1997]), so that using
a standard argument from the calculus of variations, $j^1({\phi
\circ \phi_X^{-1}})^* [W \intprod \Omega_\mathcal{L}]$ must
vanish for all vectors $W$ on
$J^1(Y)$ in order for $\phi$ to be an extremum of the action.
We shall call such elements $\phi \in \mathcal{C}$ covering
$\phi_X$, solutions of the Euler-Lagrange equations.
\begin{defn}
We let
\begin{equation}
{\mathcal P}=\{ \phi \in \mathcal{C} \ | \
j^1(\phi \circ \phi_X^{-1})^* [W \intprod \Omega_\mathcal{L}]=0
\ \; \mbox{for all} \; \ W\in TJ^1(Y)\}.
\label{s8}
\end{equation}
In coordinates, $({\phi \circ \phi_X^{-1}})^A$ 
is an element of
${\mathcal P}$ if
$$ \frac{\partial L}{\partial y^A}
(j^1({\phi \circ \phi_X^{-1}}))
-\frac{\partial }{\partial x^\mu}
\left( \frac{\partial L}{\partial
v^A_\mu} (j^1({\phi \circ \phi_X^{-1}})\right) = 0 
\text{ in } U_X.
\qquad
$$
\end{defn}

We are now ready to prove the multisymplectic form formula, a
generalization of the symplectic flow theorem, but we first
make the following remark. If
${\mathcal P}$ is a submanifold of $\mathcal{C}$, then for any
$\phi\in {\mathcal P}$, we may identify $T_\phi{\mathcal P}$
with the set $\{V \in T_\phi\mathcal{C} \ | \ j^1({\phi \circ
\phi_X^{-1}})^* {\mathfrak L}_{j^1(V)}[W\intprod
\Omega_\mathcal{L}]=0 \ \; \mbox{for all} \; \ W \in TJ^1(Y)
\}$  since such vectors arise by differentiating
$\frac{d}{d \epsilon}|_{\epsilon = 0} j^1(\phi^\epsilon \circ
{\phi^\epsilon_X}^{-1})^* [W\intprod
\Omega_\mathcal{L}] =0$, where $\phi^{\epsilon}$ is a smooth
curve of solutions of the Euler-Lagrange equations in
${\mathcal P}$ (when such solutions exist).  More generally, we
do not require ${\mathcal P}$ to be a submanifold in order to
define the first variation solution of the Euler-Lagrange
equations.

\begin{defn}
For any $\phi \in {\mathcal P}$ ,we  define the set
\begin{equation} 
{\mathcal F}=\{ V \in T_\phi\mathcal{C} \ | \
j^1({\phi \circ \phi_X^{-1}})^* {\mathfrak L}_{j^1(V)}[W
\intprod \Omega_ \mathcal{L}]=0 \ \; \mbox{for all} \; \ W \in
TJ^1(Y)\}.
\label{s10}
\end{equation} 
Elements of ${\mathcal F}$ solve the first
variation equations of the Euler-Lagrange equations.
$\qquad $
\end{defn}

\begin{thm}[Multisymplectic form formula]
\label{thm1}
If $\phi \in {\mathcal P}$, then for all $V$ and $W$ in ${\mathcal F}$,
\begin{equation}
\int_{\partial U_X} j^1({\phi \circ \phi_X^{-1}})^* [{j^1(V)}\intprod
{j^1(W)} \intprod \Omega_\mathcal{L}]=0.
\label{s11}
\end{equation}
\end{thm}

\paragraph{Proof.} We define the $1$-forms $\alpha_1$ and
$\alpha_2$ on $\mathcal{C}$ by
$$ 
  \alpha_1(\phi) \cdot V :=
  \int_{U_X} j^1({\phi \circ \phi_X^{-1}})^* [{j^1(V)}\intprod
  \Omega_\mathcal{L}]
$$
and
$$
  \alpha_2(\phi) \cdot V :=
  \int_{\partial U_X} j^1({\phi \circ \phi_X^{-1}})^*
  [{j^1(V)}\intprod \Theta_\mathcal{L}],
$$
so that by (\ref{s6}),
\begin{equation}
dS_\phi \cdot V = \alpha_1 (\phi)\cdot V + \alpha_2 (\phi)\cdot V  \ \;
\mbox{for all} \; \ V \in T_\phi\mathcal{C}.
\label{s13}
\end{equation}

Recall that for any $1$-form $\alpha$ on $\mathcal{C}$ and vector fields
$V,W$ on $\mathcal{C}$,
\begin{equation}
d \alpha (V,W) = V[\alpha(W)] - W[\alpha(V)] - \alpha([V,W]).
\label{dalpha}
\end{equation}
We let $\phi_\epsilon = \eta_{Y}^\epsilon\circ\phi$ be a curve
in $\mathcal{C}$ through $\phi$, where $\eta_{Y}^\epsilon$ is a curve
in $\mathcal{G}$ through the identity such that
$$W=\frac{d}{d \epsilon}|_{\epsilon =0} \eta_{Y}^\epsilon
\text{  and  } W \in {\mathcal F},$$
and consider equation (\ref{s13}) restricted to all $V \in {\mathcal F}$.

Thus, 
\begin{eqnarray*}
d(\alpha_2(V))(\phi) \cdot W &=& \left.\frac{d}{d \epsilon}\right|
_{\epsilon =0} (\alpha_2(V)(\phi_\epsilon)) \\
 &=& \left.\frac{d}{d \epsilon}\right|_{\epsilon =0}
\int_{\partial U_X} {j^1}({\phi \circ \phi_X^{-1}})^* {j^1}
(\eta_{Y}^\epsilon)[{j^1(V)} \intprod \Theta_\mathcal{L}]\\
&=& \int_{\partial U_X} {j^1}({\phi \circ \phi_X^{-1}})^*
{\mathfrak L}_{j^1(W)}[{j^1(V)} \intprod \Theta_\mathcal{L}] \\
&=& -\int_{\partial U_X} {j^1}({\phi \circ \phi_X^{-1}})^*
[{j^1(W)} \intprod d({j^1(V)} \intprod \Theta_\mathcal{L})] \\
&&\qquad\mbox+
\int_{\partial U_X} {j^1}({\phi \circ \phi_X^{-1}})^*
d [{j^1(W)} \intprod {j^1(V)} \intprod \Theta_\mathcal{L}],
\end{eqnarray*}
where the last equality was obtained using Cartan's formula.  Using
Stoke's theorem, noting that $\partial \partial U$ is empty, and
applying Cartan's formula once again, we obtain that
\begin{eqnarray*}
d(\alpha_2(\phi)(V)) \cdot W &=&
\int_{\partial U_X} {j^1}({\phi \circ \phi_X^{-1}})^* [{j^1(W)}
\intprod {j^1(V)} \intprod \Omega_\mathcal{L}] \\
&& \qquad - \int_{\partial U_X}
{j^1}({\phi \circ \phi_X^{-1}})^* [{j^1(W)} \intprod
{\mathfrak L}_{j^1(V)}\Theta_\mathcal{L}],
\end{eqnarray*}
and
\begin{eqnarray*}
d(\alpha_2(\phi)(W)) \cdot V &=&
\int_{\partial U_X} {j^1}({\phi \circ \phi_X^{-1}})^* [{j^1(V)}
\intprod {j^1(W)} \intprod \Omega_\mathcal{L}] \\
&& \qquad - \int_{\partial U_X}
{j^1}({\phi \circ \phi_X^{-1}})^* [{j^1(V)} \intprod
{\mathfrak L}_{j^1(W)}\Theta_\mathcal{L}].
\end{eqnarray*}
Also, since $[{j^1(V)},{j^1(W)}] = {j^1}([V,W])$, we have
$$
\alpha_2(\phi)([V,W]) = \int_{\partial U_X} {j^1}({\phi \circ
\phi_X^{-1}})^*[{j^1(V)},{j^1(W)}] \intprod \Theta_\mathcal{L}.$$
Now
$$[{j^1(V)},{j^1(W)}]\intprod \Theta_\mathcal{L} =
{\mathfrak L}_{j^1(V)} ({j^1(W)}\intprod \Theta_\mathcal{L})
- {j^1(W)} \intprod {\mathfrak L}_{j^1(V)} \Theta_\mathcal{L},$$
 so that
\begin{eqnarray*}
d\alpha_2 (\phi)(V,W)& =& 2 \int_{\partial U_X}
{j^1}({\phi \circ \phi_X^{-1}})^*[{j^1(V)}\intprod {j^1(W)}
 \intprod \Omega_\mathcal{L}] \\
&& + \int_{\partial U_X} {j^1}({\phi \circ \phi_X^{-1}})^*
[{j^1(V)}\intprod {\mathfrak L}_{j^1(W)} \Theta_\mathcal{L}
- {\mathfrak L}_{j^1(V)}({j^1(W)}\intprod \Theta_\mathcal{L})].
\end{eqnarray*}
But
$$
{\mathfrak L}_{j^1(V)}({j^1(W)}\intprod \Theta_\mathcal{L})=
d({j^1(V)}\intprod {j^1(W)} \intprod \Theta_\mathcal{L})+
{j^1(V)}\intprod d({j^1(W)}\intprod \Theta_\mathcal{L})$$
and
$${j^1(V)}\intprod {\mathfrak L}_{j^1(W)}\Theta_\mathcal{L} =
{j^1(V)}\intprod d({j^1(W)}\intprod \Theta_\mathcal{L})+
{j^1(V)}\intprod {j^1(W)}\intprod \Omega_\mathcal{L}.$$
Hence,
\begin{eqnarray*}
\lefteqn{\int_{\partial U_X} {j^1}({\phi \circ \phi_X^{-1}})^*
[{j^1(V)}\intprod {\mathfrak L}_{j^1(W)} \Theta_\mathcal{L}
- {\mathfrak L}_{j^1(V)}({j^1(W)}\intprod \Theta_\mathcal{L})]} \\
&=&
\int_{\partial U_X} {j^1}({\phi \circ \phi_X^{-1}})^*
( {j^1(V)}\intprod {j^1(W)}\intprod \Omega_\mathcal{L}) \\
&& - \int_{\partial U_X} {j^1}({\phi \circ \phi_X^{-1}})^*
d( {j^1(V)}\intprod {j^1(W)}\intprod \Theta_\mathcal{L}).
\end{eqnarray*}
The last term once again vanishes by Stokes theorem together with the fact
that $\partial \partial U$ is empty, and we obtain that
\begin{equation}
d\alpha_2(\phi)(V,W) = 3\int_{\partial U_X} {j^1}({\phi \circ \phi_X^{-1}})^*
( {j^1(V)}\intprod {j^1(W)}\intprod \Omega_\mathcal{L}) .
\label{s14}
\end{equation}

We now use (\ref{dalpha}) on $\alpha_1$.  A similar computation as above yields
\begin{equation}
d(\alpha_1(\phi)\cdot V) \cdot W = \int_{U_X} j^1({\phi \circ \phi_X^{-1}})^*
{\mathfrak L}_{j^1(W)}[{j^1(V)}\intprod \Omega_\mathcal{L}]
\nonumber
\end{equation}
which vanishes for all $\phi \in {\mathcal P}$ and $W \in {\mathcal F}$.
Similarly, $d(\alpha_1(\phi)\cdot W) \cdot v = 0$ for all
$\phi \in {\mathcal P}$ and $V \in {\mathcal F}$.  Finally,
$\alpha_1(\phi)=0$ for all $\phi \in {\mathcal P}$.

Hence, since
$$0=ddS(\phi)(V,W)= d\alpha_1(\phi)(V,W) + d\alpha_2(\phi)(V,W),$$
we obtain the formula (\ref{s11}).
\quad $\blacksquare$

\paragraph{Symplecticity revisited.}
Let $\Sigma$ be a compact oriented connected boundaryless $n$-manifold
which we think of as our reference Cauchy surface, and consider the
space of embeddings of $\Sigma$ into $X$, Emb$(\Sigma,X)$; again, although
it is unnecessary for this paper, we may topologize Emb$(\Sigma, X)$
by completing the space in the appropriate $C^k$ or $H^s$-norm.

Let $B$ be an $m$-dimensional manifold.
For any fiber bundle $\pi_{BK}:K \rightarrow B$, we shall, in addition
to $\Gamma(\pi_{BK})$, use the
corresponding script letter ${\mathcal K}$ to denote the  space of sections
of $\pi_{BK}$.
The space of sections of a fiber bundle is an
infinite-dimensional manifold; in fact, it can be precisely defined and
topologized as the
manifold $\mathcal{C}$ of the previous section, where the diffeomorphisms
on the base manifold are taken to be the identity map, so
that the tangent space to ${\mathcal K}$ at $\sigma$ is given simply by
$$ T_\sigma{\mathcal K} = \{ W:B \rightarrow VK \ | \pi_{K,TK} \circ W =
\sigma\},$$
where $VK$ denotes the vertical tangent bundle of $K$.  We let
$\pi_{K, L(VK, \Lambda^{m}(B))}:L(VK, \Lambda^{m}(B)) \rightarrow K$ be the
vector bundle over $K$ whose fiber at $k\in K_x$, $x=\pi_{BK}(k)$,
is the set of linear mappings
from $V_kK$ to $\Lambda^{m}(B)_x$.  Then the cotangent space to ${\mathcal
K}$ at
$\sigma$ is defined as
$$ T_\sigma^* {\mathcal K} = \{ \pi : B \rightarrow
L(VK, \Lambda^{m}(B)) \ | \
\pi_{K, L(VK, \Lambda^{n+1}(B))} \circ \pi = \sigma \}.$$
Integration provides
the natural pairing of $T^*_\sigma {\mathcal K}$ with $T_\sigma {\mathcal K}$:
$$ \langle \pi , V \rangle = \int_B \pi \cdot V.$$
In practice, the manifold $B$ will either be $X$ or some
($n+1$)-dimensional
subset of $X$, or the $n$-dimensional manifold $\Sigma_\tau$, where for each
$\tau \in \text{\rm Emb}(\Sigma, X)$, $\Sigma_\tau := \tau (\Sigma)$.  We
shall use the notation $Y_\tau$ for the bundle $\pi_{\Sigma_\tau,Y}$, and
${\mathcal Y}
_\tau$ for sections of this bundle.  For the remainder of this section, we
shall set the manifold $\mathcal{C}$ introduced earlier to
${\mathcal Y}$.

The infinite-dimensional manifold ${\mathcal Y}_\tau$ is called
the {\bfi $\tau$-configuration space}, its tangent bundle is called the
{\bfi $\tau$-tangent space}, and its cotangent bundle
$T^*{\mathcal Y}_\tau$ is called the {\bfi $\tau$-phase space}.  Just as we
described in Section \ref{mechanics}, the cotangent bundle has a canonical
$1$-form $\theta_\tau$ and a canonical $2$-form $\omega_\tau$.  These
differential
forms are given by
\begin{equation}\label{t1}
\theta_\tau (\varphi, \pi) \cdot V = \int_{\Sigma_\tau}
\pi( T\pi_{{\mathcal Y}_\tau,T^*{\mathcal Y}_\tau} \cdot V) \text{  and  }
\omega_\tau = -d \theta_\tau,
\end{equation}
where $(\varphi, \pi) \in {\mathcal Y}_\tau$, $V \in T_{(\varphi, \pi)} T^*
{\mathcal Y}_\tau$, and $\pi_{{\mathcal Y}_\tau,T^*{\mathcal Y}_\tau}:
T^*{\mathcal Y}_\tau \rightarrow {\mathcal Y}_\tau $ is the cotangent
bundle projection map.

An infinitesimal slicing of the bundle $\pi_{XY}$ consists of $Y_\tau$ together
with a vector field $\zeta$ which is everywhere transverse to $Y_\tau$, and
covers $\zeta_X$ which is everywhere transverse to ${\Sigma_\tau}$.  The
existence of an infinitesimal slicing allows us to invariantly decompose the
temporal from the spatial derivatives of the fields.  Let $\phi \in
{\mathcal Y}$,
$\varphi := \phi|_{\Sigma_\tau}$, and let $i_\tau: {\Sigma_\tau} \rightarrow X$
be the inclusion map.  Then we may define the map $\beta_\zeta$
taking $j^1({\mathcal Y})_\tau$  to $j^1({\mathcal Y}_\tau) \times \Gamma(
\pi_{{\Sigma_\tau},VY_\tau})$ over ${\mathcal Y}_\tau$ by
\begin{equation}\label{t2}
\beta_\zeta ({j^1}(\phi) \circ i_\tau) = (j^1(\varphi), \dot{\varphi})
\text{  where  } \dot{\varphi}:= {\mathfrak L}_\zeta \phi  .
\end{equation}
In our notation, $j^1({\mathcal Y})_\tau$ is the collection of restrictions
of holonomic sections of $J^1(Y)$ to ${\Sigma_\tau}$, while
$j^1({\mathcal Y}_\tau)$ are the holonomic sections of $\pi_{{\Sigma_\tau},
J^1(Y)}$.  It is easy to see that $\beta_\zeta$ is an isomorphism;
it then follows that $\beta_\zeta$ is an isomorphism of $j^1({\mathcal
Y})_\tau$
with $T{\mathcal Y}_\tau$, since $j^1(\varphi)$ is completely determined by
$\varphi$.
This bundle map is called the jet decomposition map, and its inverse is called
the jet reconstruction map.  Using this map, we can define the instantaneous
Lagrangian.
\begin{defn}
The instantaneous Lagrangian $L_{\tau,\zeta}:T{\mathcal Y}_\tau \rightarrow
\mathbb{R}$ is given by
\begin{equation}
L_{\tau,\zeta}(\varphi, \dot{\varphi}) =
\int_{\Sigma_\tau} i_\tau^* [\zeta_X \intprod \mathcal{L}
(\beta_\zeta^{-1}(j^1(\varphi),\dot{\varphi})]
\label{t3}
\end{equation}
for all $(\varphi,\dot{\varphi}) \in T{\mathcal Y}_\tau$. $\qquad $
\end{defn}
The instantaneous Lagrangian $L_{\tau,\zeta}$ has an instantaneous Legendre
transform
$$ \mathbb{F}L_{\tau,\zeta} : T{\mathcal Y}_\tau \rightarrow
T^* {\mathcal Y}_\tau; \ \ (\varphi, \dot{\varphi}) \mapsto (\varphi, \pi)$$
which is defined in the usual way by vertical fiber differentiation of
$L_{\tau,\zeta}$ (see, for example,  Abraham and Marsden [1978]).
Using the instantaneous Legendre transformation, we can
pull-back the canonical $1$- and $2$-forms on $T^*{\mathcal Y}_\tau$.
\begin{defn}
Denote, respectively,  the instantaneous Lagrange $1$- and $2$-forms on
 $T{\mathcal Y}_\tau$ by
\begin{equation}\label{t4}
\theta_\tau^L = \mathbb{F} L_{\tau, \zeta}^* \theta_\tau
\text{  and  }
\omega_\tau^L = -d\theta_\tau^L. \qquad
\end{equation}
\end{defn}
Alternatively, we may define $\theta_\tau^L$ using Theorem \ref{thm2.1},
in which case no reference to the cotangent bundle is necessary.

We will show that our covariant multisymplectic form formula can be used to
recover the fact that the flow of the Euler-Lagrange equations in
the bundle
$\pi_{\text{\rm Emb}(\Sigma,X),\cup_{\tau\in\text{\rm Emb}(\Sigma,X)}
T{\mathcal Y}_\tau}$
 is symplectic with respect to $\omega_\tau^L$.  To do so, we must
relate the multisymplectic Cartan ($n+2$)-form
$\Omega_\mathcal{L}$ on $J^1(Y)$ with the symplectic $2$-form $\omega_\tau^L$
on $T{\mathcal Y}_\tau$.

\begin{thm}
\label{thm3.2}
Let $\Theta_\tau^\mathcal{L}$ be the canonical $1$-form on $j^1({\mathcal
Y})_\tau$
given by
\begin{equation}\label{t5}
\Theta_\tau^\mathcal{L}({j^1}(\phi) \circ i_\tau) \cdot V = \int_{\Sigma_\tau}
i_\tau^* {j^1}(\phi)^* [V \intprod \Theta_\mathcal{L}],
\end{equation}
where ${j^1}(\phi) \circ i_\tau \in j^1({\mathcal Y})_\tau$, $V \in
T_{{j^1}(\phi)\circ i_\tau} j^1({\mathcal Y})_\tau$. \newline
(a) If the $2$-form $\Omega_\tau^\mathcal{L}$ on $j^1({\mathcal Y})_\tau$ is
defined by $\Omega_\tau^\mathcal{L} := -d \Theta_\tau^\mathcal{L}$, then
for $V, W \in T_{{j^1}(\phi)\circ i_\tau} j^1({\mathcal Y})_\tau$,
\begin{equation}\label{t6}
\Omega_\tau^\mathcal{L}({j^1}(\phi) \circ i_\tau) (V,W) = \int_{\Sigma_\tau}
i_\tau^* {j^1}(\phi)^* [W \intprod V \intprod \Omega_\mathcal{L}].
\end{equation}
\noindent
(b)Let the diffeomorphism $s_X: \Sigma \times \mathbb{R} \rightarrow X$ be a
slicing of $X$ such that for $\lambda \in \mathbb{R}$,
$$ \Sigma_\lambda := s_X (\Sigma \times \{ \lambda \} ) \text{  and  }
\Sigma_\lambda := \tau_\lambda(\Sigma),$$
where $\tau_\lambda \in \text{\rm Emb}(\Sigma, X)$ is given by
$\tau_\lambda(x) = s_X(x,\lambda)$.
For any $\phi \in {\mathcal P}$, let
$V, W \in T_\phi {\mathcal Y}
\cap {\mathcal F}$ so that for each $\tau \in \text{\rm Emb}(\Sigma,X)$,
$j^1V_\tau, j^1W_\tau  \in
T_{{j^1}(\phi) \circ i_\tau} j^1({\mathcal Y})_\tau$, and let
$\tau_{\lambda_1},
\tau_{\lambda_2} \in \text{\rm Emb}(\Sigma,X)$.  Then
\begin{equation}\label{t7}
\Omega_{\tau_{\lambda_1}}^\mathcal{L}(j^1V_{\tau_{\lambda_1}},j^1W_{\tau_{\lambda_1}}) =
\Omega_{\tau_{\lambda_2}}^\mathcal{L}(j^1V_{\tau_{\lambda_2}},j^1W_{\tau_{\lambda_2}}).
\end{equation}
\end{thm}
\paragraph{Proof.} Part (a) follows from the Cartan formula
together with Stokes theorem using an argument like that in the
proof of Theorem \ref{thm1}.

For part (b), we recall that the multisymplectic form formula on ${\mathcal
Y}$
states that for any subset $U_X \subset X$ with smooth closed boundary and
vectors
$V,W \in T_\phi {\mathcal Y} \cap {\mathcal F}$, $\phi \in {\mathcal Y}$,
\begin{equation}\label{t8}
\int_{\partial U_X} {j^1}(\phi)^* [{j^1(V)}\intprod {j^1(W)} \intprod
\Omega_\mathcal{L}]=0.
\end{equation}
Let
$$U_X = \cup_{\lambda\in [\lambda_1, \lambda_2]} \Sigma_\lambda.$$
Then $\partial U_X = \Sigma_{\lambda_1} - \Sigma_{\lambda_2}$, so that
(\ref{t8}) can
be written as
\begin{eqnarray*}
0 &=& \int_{\Sigma_{\lambda_2}} {j^1}(\phi \circ i_{\tau_{\lambda_2}})^*
[ {j^1V_{\tau_{\lambda_2}}} \intprod {j^1W_{\tau_{\lambda_2}}} \intprod
\Omega_\mathcal{L}] \\
&& \qquad - \int_{\Sigma_{\lambda_1}} {j^1}(\phi \circ i_{\tau_{\lambda_1}})^*
[ {j^1V_{\tau_{\lambda_1}}} \intprod {j^1W_{\tau_{\lambda_1}}} \intprod
\Omega_\mathcal{L}] \\
&=&
\Omega_{\tau_{\lambda_1}}^\mathcal{L}(j^1V_{\tau_{\lambda_1}},j^1W_{\tau_{\lambda_1}})-
\Omega_{\tau_{\lambda_2}}^\mathcal{L}(j^1V_{\tau_{\lambda_2}},j^1W_{\tau_{\lambda_2}}),
\end{eqnarray*}
which proves (\ref{t7}).
\quad $\blacksquare$

\begin{thm}
\label{thm3.3}
The identity $ \Theta_\tau^\mathcal{L} =  \beta_\zeta^*
\theta_\tau^L
$ holds.
\end{thm}
\paragraph{Proof.} Let $W \in T_{{j^1}(\phi)\circ i_\tau}
j^1({\mathcal Y})_\tau$, which we identify
with $w \circ \phi \circ i_\tau$, where $w$ is a $\pi_{X,J^1(Y)}$-vertical
vector.
Choose a coordinate chart which is adapted to the slicing so that
$\partial_0 |_{Y_\tau} = \zeta$.  With $w = (0, W^A, W^A_\mu)$, we see that
$$ \Theta_\tau^\mathcal{L} \cdot W = \int_{{\Sigma_\tau}}
\frac{\partial L}{\partial {v^A}_0} (\phi^B, {\phi^B}_{,\mu})W^A d^nx_0.$$
Now, from (\ref{t3}) we get
\begin{eqnarray*}\label{t9}
\theta_\tau^L (\varphi, \dot{\varphi}) &=&
\frac{\partial L_{\tau,\zeta}}{\partial \dot{y}^A} dy^A \\
&=& \int_{\Sigma_\tau}
\frac{\partial }{\partial \dot{y}^A} i_\tau^*
[\partial _0 \intprod L(x^\mu, \phi^A, {\phi^A}_{,\mu}) d^{n+1}x \otimes
dy^A]\\
&=& \int_{\Sigma_\tau}
\frac{\partial L}{\partial {v^A}_0} (\phi^B, {\phi^B}_{,\mu})dy^A\otimes
d^nx_0,
\end{eqnarray*}
where we arrived at the last equality using the fact that $\dot{y}^A = {v^A}_0$
in this adapted  chart.  Since $(T\beta_\zeta \cdot W)^A = W^A$, we see that
$\Theta_\tau^\mathcal{L} \cdot W = \theta_\tau^L \cdot (T\beta_\zeta \cdot
W)$,
and this completes the proof.
\quad $\blacksquare$
\medskip 

Let the instantaneous energy $E_{\tau,\zeta}$ associated with $L_{\tau,\zeta}$
be given by
\begin{equation}\label{t10}
E_{\tau,\zeta}(\varphi, \dot{\varphi}) = \mathbb{F}
L_{\tau,\zeta} (\dot{\varphi}) \cdot \dot{\varphi} -
L_{\tau,\zeta} (\varphi, \dot{\varphi}),
\end{equation}
and define the ``time''-dependent Lagrangian vector field
$X_{E_{\tau,\zeta}}$ by
$$X_{E_{\tau,\zeta}} \intprod \omega_\tau^L = dE_{\tau,\zeta}.$$
Since $\cup_{\tau \in \text{\rm Emb}(\Sigma,X)}T{\mathcal Y}_\tau$
over Emb$(\Sigma,X)$ is infinite-dimensional
and $w_\tau^L$ is only weakly nondegenerate, the second-order vector field
$X_{E_{\tau,\zeta}}$ does not, in general, exist.  In the case that it
does, we obtain
the following result.
\begin{cor}
Assume $X_{E_{\tau,\zeta}}$ exists and let $F_\tau$ be its semiflow,
defined on some subset ${\mathcal D}$ of  the bundle
$\cup_{\tau \in \text{\rm Emb}(\Sigma,X)} T{\mathcal Y}_\tau$ over
Emb$(\Sigma,X)$. Fix $\bar{\tau}$ so that $F_{\bar{\tau}}(\varphi_1,
\dot{\varphi_1}) = (\varphi_2, \dot{\varphi_2})$ where
$(\varphi_1, \dot{\varphi_1}) \in T {\mathcal Y}_{\tau_1}$ and
$(\varphi_2, \dot{\varphi_2}) \in T {\mathcal Y}_{\tau_2}$.  Then
$F_{\bar{\tau}}^* \omega_{\tau_2}^L = \omega_{\tau_1}^L$.
\end{cor}
\paragraph{Proof.}
This follows immediately from Theorem \ref{thm3.2}(b) and Theorem \ref{thm3.3}
and the fact that $\beta_\zeta$ induces an isomorphism between
$j^1({\mathcal Y})_\tau$
and $T{\mathcal Y}_\tau$.
\quad $\blacksquare$

\paragraph{Example: nonlinear wave equation.} To illustrate
the geometry that we have developed, let us consider the scalar nonlinear wave
equation given by
\begin{equation}
\frac{\partial ^2 \phi}{\partial {x^0}^2} - \triangle \phi - N'(\phi)=0,
\ \ \phi \in \Gamma(\pi_{XY}),
\label{t11}
\end{equation}
where $\triangle$ is the Laplace-Beltrami operator and $N$ is a real-valued
$C^\infty$ function of one variable.  For concreteness, fix $n$$=$$1$ so that
the spacetime manifold $X:= \mathbb{R}^2$, the configuration bundle
$Y:= \pi_{\mathbb{R}^2,\mathbb{R}}$, and the first jet bundle
$J^1(Y):= \pi_{\mathbb{R}^2, \mathbb{R}^3}$.

Equation (\ref{t11}) is governed by the Lagrangian density
\begin{equation}\label{lag}
\mathcal{L} = \left\{
{\frac{1}{2}}\left[ {\frac{\partial \phi}{\partial x^0}}^2 -
{\frac{\partial \phi}{\partial x^1}}^2 \right] + N(\phi)
\right\} dx^1 \wedge dx^0 .
\end{equation}
Using coordinates $(x^0,x^1,\phi, \phi_{,0}, \phi_{,1})$ for $J^1(Y)$, we write
the multisymplectic $3$-form for this nonlinear wave equation on $\mathbb{R}^2$
in coordinates as
\begin{eqnarray}
\Omega_\mathcal{L} &=& -d\phi \wedge d\phi_{,0} \wedge dx^1 -
d\phi \wedge d\phi_{,1} \wedge dx^0 - N'(\phi) d \phi \wedge dx^1 \wedge dx^0
\nonumber \\
&& \qquad + \phi_{,0}d \phi_{,0} \wedge dx^1 \wedge dx^0
- \phi_{,1}d \phi_{,1} \wedge dx^1 \wedge dx^0;
\label{t12}
\end{eqnarray}
a short computation verifies that solutions of (\ref{t11}) are elements of
${\mathcal P}$, or that ${j^1}(\phi \circ \phi_X^{-1})^* [W \intprod
\Omega_\mathcal{L}]
=0$ for all $W \in T J^1(Y)$ (see Marsden and Shkoller [1997]).

We will use this example to demonstrate that our multisymplectic form formula
generalizes the notion of symplecticity given by Bridges [1997].
Since the Lagrangian (\ref{lag}) does not explicitly depend on time,
it is convenient to identify sections of $Y$ as mappings from
$\mathbb{R}^2$ into $\mathbb{R}$, and similarly,
sections of $J^1(Y)$  as mappings from $\mathbb{R}^2$ into $\mathbb{R}^3$.
Thus, for $\phi \in \Gamma(\pi_{XY})$,
${j^1}(\phi)(x^\mu) := (\phi(x^\mu), \phi_{,0}(x^\mu), \phi_{,1}(x^\mu)) \in
\mathbb{R}^3$, and if we set $p^\mu := \phi_{,\mu}$,
then (\ref{t11}) can be reformulated to
\begin{eqnarray}
&&{\bf J}_0 j^1\phi_{, 0}+ {\bf J}_1 j^1\phi_{, 1} : =
\nonumber
\\
&&
\left[\begin{array}{cccccc}
0   & 1  & 0\\
-1   & 0  & 0\\
0   & 0  & 0\\
\end{array} \right]
\left[\begin{array}{cccccc}
 \phi  \\
p^0    \\
p^1   \\
\end{array} \right] _{, 0}
+
\left[\begin{array}{cccccc}
0   & 0  & -1\\
0   & 0  & 0\\
1  & 0  & 0\\
\end{array} \right]
\left[\begin{array}{cccccc}
 \phi  \\
p^0    \\
p^1   \\
\end{array} \right] _{, 1}
=
\left[\begin{array}{cccccc}
N'(\phi)  \\
-p^0    \\
p^1   \\
\end{array} \right]
.
\qquad \qquad
\label{t12a}
\end{eqnarray}
To each degenerate matrix ${\bf J_\mu}$,
we associate the contact form $\omega^\mu$ on $\mathbb{R}^3$
given by $\omega^\mu (u_1, u_2) = \langle {\bf J_\mu} u_1, u_2 \rangle$,
where $u_1,u_2 \in \mathbb{R}^3$ and $\langle \cdot, \cdot \rangle$ is the
standard inner product on $\mathbb{R}^3$.
Bridges obtains the following conservation of symplecticity:
\begin{equation}\label{t14}
\frac{\partial }{\partial x^0} \left[ \omega^0
(j^1(\phi_{,0}),j^1(\phi_{,1})) \right] +
\frac{\partial }{\partial x^1} \left[ \omega^1
(j^1(\phi_{,0}),j^1(\phi_{,1})) \right] =0.
\end{equation}

This result is interesting, but has somewhat limited scope in that the
vector fields in (\ref{t14}) upon which  the contact forms act
are not general solutions to the first variation equations; rather,
they are the specific first variation solutions $\phi_{,\mu}$.  Bridges
obtains this result by crucially relying on the multi-Hamiltonian
structure of (\ref{t11}); in particular, the vector $(N'(\phi), -p^0, p^1)$
on the right-hand-side of (\ref{t12a})
is the gradient of a smooth multi-Hamiltonian function  $H(\phi,p^0,p^1)$
(although the multi-Hamiltonian formalism is not important for this
article, we refer the reader to Marsden and Shkoller [1997] for the
Hamiltonian version of our covariant framework, and to Bridges [1997]).
Using equation (\ref{t12a}), it is clear that
$$H_{,0}=\omega^0 (j^1(\phi_{,0}),j^1(\phi_{,1})) \text{  and  }
H_{,1}=-\omega^1 (j^1(\phi_{,0}),j^1(\phi_{,1}))$$
so that (\ref{t14}) follows from the relation $H_{,0,1} = H_{,1,0}$.

\begin{prop} \label{prop1}
The multisymplectic form formula is an intrinsic generalization of
the conservation law (\ref{t14}); namely, for any $V,W \in {\mathcal F}$
that are $\pi_{X,J^1(Y)}$-vertical,
\begin{equation}\label{t15}
\frac{\partial }{\partial x^0} \left[ \omega^0
(j^1(V),j^1(W)) \right] + \frac{\partial }{\partial x^1} \left[ \omega^1
(j^1(V),j^1(W)) \right] =0.
\end{equation}
\end{prop}
\paragraph{Proof.}
Let ${j^1(V)}$ and ${j^1(W)}$ have the  coordinate expressions $(V, V^0, V^1)$
and $(W, W^0, W^1)$, respectively.  Using (\ref{t12}), we compute
$$ j^1(W)\intprod {j^1(V)} \intprod \Omega_\mathcal{L} =
\left( V W^0 - V^0 W \right) dx + \left( V W^1 - V^1 W \right) dt , $$
so that with Theorem \ref{thm1} and the definition of $\omega^\mu$, we
have, for $U_X \subset X$,
$$ \int_{\partial U_X} \omega^0 ({j^1(V)},{j^1(W)}) dx -
\omega^1 ({j^1(V)},{j^1(W)}) dt = 0,$$
and hence by Green's theorem,
$$ \int_{U_X}
\left\{ \frac{\partial }{\partial x^0} \left[ \omega^0
(j^1(V),j^1(W)) \right] + \frac{\partial }{\partial x^1} \left[ \omega^1
(j^1(V),j^1(W)) \right] \right\}dx^1 \wedge dx^0 =0.$$
Since $U_X$ is arbitrary,  we obtain the desired result.
\quad $\blacksquare$
\medskip

In general, when $V$ is $\pi_{XY}$-vertical, $j^1(V)$ has the
coordinate expression $(V, V_{,\mu} + \partial V / \partial \phi \cdot
\phi_{,\mu}$, but for the special case that $V = \phi_{,\mu}$,
$j^1(\phi_{\mu}) = (j^1\phi)_{,\mu}$, and Proposition \ref{prop1}
gives
$$
\frac{\partial }{\partial x^0} \left[ \phi_{0} \phi_{,0,1} -
\phi_{1} \phi_{,0,0} \right] -
\frac{\partial }{\partial x^1} \left[ \phi_{0} \phi_{,1,1} -
\phi_{1} \phi_{,0,1} \right] = 0,$$
which simplifies to the trivial statement that
$$ \phi_{,0} N(\phi)_{,1} - \phi_{,1} N(\phi)_{,0} =0.$$

\paragraph{The variational route to the Cartan form.}
We may alternatively define the Cartan form by beginning with equation
(\ref{s5}).  Using the infinitesimal generators defined in (\ref{ig}),
we obtain that
\begin{eqnarray}
d\mathcal{S}_\phi\cdot V &=& \left.\frac{d}{d \lambda}\right|_{\lambda=0}
 \mathcal{S}(\Phi(\eta_Y^\lambda,\phi)) \nonumber\\
&=& \int_{\eta_{X}^\lambda(U_X)} \left.\frac{d}{d \lambda}\right|_{\lambda=0}
\mathcal{L}(j^1(\Phi(\eta_{Y}^\lambda, \phi)))\nonumber\\
&=& \int_{U_X} \left.\frac{d}{d \lambda}\right|_{\lambda=0}
\mathcal{L}(j^1(\Phi(\eta_{Y}^\lambda, \phi))) \nonumber \\
&& \qquad
 + \int_{U_X} {\mathfrak L}_{V_X}\left[
\mathcal{L}(j^1({\phi \circ \phi_X^{-1}})) \right]. \label{s18}
\end{eqnarray}
Using the natural splitting of $TY$, any vector $V\in T_\phi\mathcal{C}$
may decomposed as
\begin{equation}
V=V^h+V^v, \text{  where  } V^h= T({\phi \circ \phi_X^{-1}})
\cdot V_X \text{ and } V^v = V - V^h,
\label{s18a}
\end{equation}
where we recall that $V_X = T\pi_{XY} \cdot V$.
\begin{lemma}
\label{lemma1}
For any $V \in T_\phi\mathcal{C}$,
\begin{equation}
d\mathcal{S}_\phi\cdot V^h = \int_{\partial U_X} V_X \intprod [\mathcal{L}
(j^1({\phi \circ \phi_X^{-1}}))],
\label{s19}
\end{equation}
and
\begin{equation}
d\mathcal{S}_\phi  \cdot V^v =
\int_{U_X} \left.\frac{d}{d \lambda}\right|_{\lambda=0}
\mathcal{L}(j^1(\Phi(\eta_{Y}^\lambda, \phi)))
\label{s20}
\end{equation}
\end{lemma}
\paragraph{Proof.}
The equality (\ref{s20}) is obvious, since the second term in
(\ref{s18}) clearly vanishes for all vertical vectors. For
vectors $V^h$, the first term in (\ref{s18}) vanishes; indeed,
using the chain rule, we need only compute that
\begin{equation}
\left.\frac{d}{d \lambda}\right|_{\lambda=0}
\eta_{Y}^\lambda \circ \phi \circ {\eta_{X}^\lambda}^{-1}
= V^h - T({\phi \circ \phi_X^{-1}}) \cdot V_X,
\nonumber
\end{equation}
which is zero by (\ref{s18a}).  We then apply the Cartan formula
to the second term in (\ref{s18}) and note that $d\mathcal{L}$
is an ($n+2$)-form on the ($n+1$)-dimensional manifold
$U_X$ so that we obtain (\ref{s19}).
\quad $\blacksquare$

\begin{thm}
\label{thm_cartan}
Given a smooth Lagrangian density $\mathcal{L}:J^1(Y) \rightarrow
\Lambda^{n+1}(X)$, there exist a unique smooth section $D_{EL}
\mathcal{L} \in C^\infty (Y'', \Lambda^{n+1}(X)
\otimes T^*Y))$ and a unique differential form $\Theta_\mathcal{L} \in
\Lambda^{n+1}(J^1(Y))$ such that for any $V \in T_\phi\mathcal{C}$,
and any open subset $U_X$ such that $\overline{U}_X \cap \partial X
= \emptyset$,
\begin{equation}
d \mathcal{S}_\phi \cdot V =  \int_{U_X}
D_{EL} \mathcal{L}(j^2(\phi \circ \phi_X^{-1})) \cdot V
+ \int_{\partial U_X}
j^1({\phi \circ \phi_X^{-1}})^* [{j^1(V)} \intprod \Theta_\mathcal{L}].
\label{s21}
\end{equation}
Furthermore,
\begin{equation}
D_{EL} \mathcal{L}(j^2(\phi \circ \phi_X^{-1})) \cdot V
= j^1(\phi \circ \phi_X^{-1})^*[j^1(V) \intprod \Omega_\mathcal{L}]
\text{ in } U_X.
\label{s21b}
\end{equation}
In coordinates, the action of the Euler-Lagrange derivative
$D_{EL}\mathcal{L}$ on $Y''$ is given by
\begin{eqnarray}
&& D_{EL} \mathcal{L} (j^2({\phi \circ \phi_X^{-1}})) =
\left[
\frac{\partial L }{\partial y^A} (j^1({\phi \circ \phi_X^{-1}}))
 -
\frac{\partial^2 L }{\partial x^\mu \partial {v^A}_\mu} (j^1({\phi \circ
\phi_X^{-1}})) \right.  \nonumber\\
&&  \qquad -
\frac{\partial^2 L }{\partial y^B \partial {v^A}_\mu} (j^1({\phi \circ
\phi_X^{-1}}))
\cdot (\phi \circ \phi_X^{-1})^B_{,\mu} \nonumber \\
&& \qquad - \left.
\frac{\partial^2 L }{\partial {v^B}_\nu \partial {v^A}_\mu}
 (j^1({\phi \circ \phi_X^{-1}}))
\cdot (\phi\circ\phi_X^{-1})^B_{,\mu \nu} \right] dy^A \wedge d^{n+1}x ,
\label{s21c}
\end{eqnarray}
while the form $\Theta_\mathcal{L}$ matches the definition of the Cartan
form given
in (\ref{D1a}) and has the coordinate expression
\begin{equation}
\Theta_\mathcal{L} = \frac{\partial L}{\partial {v^A}_\mu} dy^A \wedge
d^nx_\mu
+\left( L - \frac{\partial L}{\partial {v^A}_\mu} {v^A}_\mu \right) d^{n+1}x.
\label{s21d}
\end{equation}
\end{thm}
\paragraph{Proof.}
Choose $U_X:= \phi_X(U)$ small enough so that it is contained in a
coordinate chart, say $O_1$.  In these coordinates, let $V=(V^\mu, V^A)$
so that along ${\phi \circ \phi_X^{-1}}$, our decomposition
(\ref{s18a}) may be written as
$$V_X= V^\mu \frac{\partial }{\partial x^\mu}  \text{  and  }
V^v= (V^v)^A \frac{\partial }{\partial y^A}:=
\left(V^A - V^\mu \frac{\partial ({\phi \circ \phi_X^{-1}})^A}
{\partial x^\mu}\right) \frac{\partial }{\partial y^A},$$
and equation (\ref{s20}) gives
\begin{equation}
d\mathcal{S}_\phi \cdot V^v =\int_{U_X} \left[
\frac{\partial L}{\partial y^A} (j^1({\phi \circ \phi_X^{-1}})) \cdot
(V^v)^A +
\frac{\partial L}{\partial {v^A}_\mu} (j^1({\phi \circ \phi_X^{-1}})) \cdot
\frac{\partial (V^v)^A}{\partial x^\mu} \right] d^{n+1}x,
\label{s70}
\end{equation}
where we have used the fact that in coordinates along $j^1({\phi \circ
\phi_X^{-1}})$,
$$ {\{j^1(V)\}^A}_\mu = \partial_\mu [ (V^v)^A (j^1({\phi \circ
\phi_X^{-1}}))].$$
Integrating (\ref{s70}) by parts, we obtain
\begin{eqnarray}
d \mathcal{S}_\phi \cdot V^v &=& \int_{U_X} \left\{ \left[
\frac{\partial L}{\partial y^A}(j^1({\phi \circ \phi_X^{-1}})) -
\frac{\partial}{\partial x^\mu}
\frac{\partial L}{\partial {v^A}_\mu}(j^1({\phi \circ \phi_X^{-1}}))\right]
\cdot V^A \right\} d^{n+1}x \nonumber \\
&& \qquad + \int_{\partial U_X}\left\{
\frac{\partial L}{\partial {v^A}_\mu}(j^1({\phi \circ \phi_X^{-1}}))
\cdot V^A d^nx_\mu \right. \nonumber \\
&&\qquad \qquad + \left.
\frac{\partial L}{\partial {v^A}_\mu}(j^1({\phi \circ \phi_X^{-1}}))
\frac{\partial ({\phi \circ \phi_X^{-1}})^A}{\partial x^\nu} \cdot
V^\nu d^{n}x_\mu \right\} .
\label{s71}
\end{eqnarray}
Let $\alpha$ be the $n$-form integrand of the boundary integral in (\ref{s71});
then $\int_{\partial U_X} \alpha = \int_{\partial j^1(\phi
 \circ \phi_X^{-1} )(U_X)}$ since $\alpha$ is invariant under this lift.
Additionally, from equation (\ref{s19}), we obtain the horizontal contribution
\begin{equation}
d \mathcal{S}_\phi \cdot V^h = \int_{\partial U_X} (V^\mu \partial _\mu)
\intprod
(Ld^{n+1}x) = \int_{\partial {j^1}(\phi \circ \phi)X^{-1})(U_X) }V^\mu L d^nx_\mu ,
\label{s72}
\end{equation}
so combining equations (\ref{s71}) and (\ref{s72}), a simply computation verifies that
\begin{eqnarray}
d \mathcal{S}_\phi \cdot V &=& \int_{U_X} \left\{ \left[
\frac{\partial L}{\partial y^A}(j^1({\phi \circ \phi_X^{-1}})) -
\frac{\partial}{\partial x^\mu}
\frac{\partial L}{\partial {v^A}_\mu}(j^1({\phi \circ \phi_X^{-1}}))
 \right] d^{n+1}x \otimes dy^A \right\} \cdot V \nonumber \\
&& \qquad + \int_{\partial {j^1}(\phi\circ \phi_X^{-1})(U_X)}V \intprod \left\{
\frac{\partial L}{\partial {v^A}_\mu}(j^1({\phi \circ \phi_X^{-1}}))
dy^A \wedge d^nx_\mu  \right. \nonumber \\
&& \qquad \qquad \qquad \left.
+ \left[ L -
\frac{\partial L}{\partial {v^A}_\mu}(j^1({\phi \circ \phi_X^{-1}}))
\frac{\partial ({\phi \circ \phi_X^{-1}})^A}{\partial x^\mu} \right]
 d^{n+1}x \right\}.
\label{s73}
\end{eqnarray}
The vector $V$ in the second term of (\ref{s73}) may be replaced by
${j^1(V)}$ since $\pi_{Y,J^1(Y)}$-vertical vectors are clearly in the kernel of
the form that $V$ is acting on.
This shows that (\ref{s21c}) and (\ref{s21d}) hold, and hence that the
boundary integral in (\ref{s73}) may be written as
\begin{equation}
\int_{\partial U_X} j^1({\phi \circ \phi_X^{-1}})^* [j^1(V) \intprod
\Theta_\mathcal{L}]. \nonumber
\end{equation}

Now, if we choose another coordinate chart $O_2$, the coordinate expressions
of $D_{EL} \mathcal{L}$ and $\Theta_\mathcal{L}$ must agree on the overlap
$O_1 \cap O_2$ since the left-hand-side of (\ref{s21}) is intrinsically
defined.  Thus, we have uniquely defined $D_{EL}\mathcal{L}$ and
$\Theta_\mathcal{L}$ for any $U_X$ such that $\overline{U}_X \cap
\partial X = \emptyset$.

Finally, (\ref{s21b}) holds, since $\Omega_\mathcal{L} = d
\Theta_\mathcal{L}$ is also intrinsically defined and both
sides of the equation yield the same coordinate representation,
the Euler-Lagrange equations in $U_X$.
\quad $\blacksquare$
\paragraph{Remark} To prove Theorem \ref{thm_cartan} for the
case $U_X =X$, we must modify the proof to take into account
the boundary conditions which are prescribed on $\partial X$.

\begin{cor}
\label{cor4.2}
The ($n+1$)-form $\Theta_\mathcal{L}$ defined by the variational
principle satisfies the relationship
$$\mathcal{L}(j^1(z)) = z^* \Theta_\mathcal{L}$$
for all holonomic sections $z\in\Gamma(\pi_{X,J^1(Y)})$.
\end{cor}
\paragraph{Proof.}
This follows immediately by substituting (\ref{s21b}) into (\ref{s21})
and integrating by parts using Cartan's formula.
\quad $\blacksquare$

\paragraph{Remark} We have thus far focused on holonomic
sections of $J^1(Y)$, those that are the first jets of sections
of $Y$, and correspondingly, we have restricted the general
splitting of $TY$ given by
$$ TY = \text{image } \gamma \oplus VY \text{ for any } \gamma
\in \Gamma(J^1(Y)),$$ to $TY = T\phi \oplus VY$, $\phi \in
\Gamma(Y)$ as we specified in (\ref{s18a}).  For general
sections $\gamma \in \Gamma(J^1(Y))$, the horizontal bundle is
given by  image $\gamma$, and the Frobenius theorem guarantees
that $\gamma$ is locally holonomic if the connection is flat,
or equivalently if the curvature of the connection $R_\gamma$
vanishes.  Since this is a local statement, we may assume that
$Y=U \times \mathbb{R}^N$, where $U \subset \mathbb{R}^{n+1}$
is open, and that $\pi_{XY}$ is simply the projection onto the
first factor. For  $\phi \in \Gamma(Y)$, and $\gamma \in
\Gamma(J^1(Y))$, $\gamma(x, \phi(x)): \mathbb{R}^{n+1}
\rightarrow \mathbb{R}^N$ is a linear operator which is
holonomic if $\phi'(x) = \gamma(x,\phi(x))$, where $\phi'(x)$
is the differential of $\phi$, and this is the case whenever
the operator $\phi''(x)$ is symmetric.  Equivalently, the
operator
$$ S_\gamma(x,y) \cdot (v,w) : = D_1 \gamma(x,y) \cdot (v,w) +
D_2\gamma(x,y) \cdot (\gamma(x,y) \cdot v, w)$$ is symmetric
for all $v,w \in \mathbb{R}^{n+1}$.  One may easily verify that
the local curvature is given by
$$ R_\gamma(x,y) \cdot (v,w) := S_\gamma(x,y) \cdot (v,w) -
S_\gamma(x,y) \cdot (w,v)$$ and that $\gamma = {j^1}(\phi)$
locally for some $\phi \in \Gamma(Y)$, if and only if $R_\gamma
=0$.

\paragraph{The variational route to Noether's Theorem.}
Suppose the Lie group $\mathcal{G}$ acts on $\mathcal{C}$ and
leaves the action $\mathcal{S}$ invariant so that
\begin{equation} \label{s22a}
\mathcal{S}(\Phi(\eta_{Y}, \phi)) = \mathcal{S}(\phi) \ \; 
\mbox{for
all} \; \
\eta_{Y} \in \mathcal{G}.
\end{equation}
This implies that for each $\eta_{Y}\in \mathcal{G}$,
$\Phi(\eta_{Y},
\phi) \in {\mathcal P}$ whenever $\phi \in {\mathcal P}$.  We
restrict the action of $\mathcal{G}$ to ${\mathcal P}$, and let
$\xi_\mathcal{C}$ be the corresponding infinitesimal generator
on $\mathcal{C}$ restricted to points in ${\mathcal P}$; then
\begin{eqnarray*}
0 = (\xi_\mathcal{C} \intprod d\mathcal{S})_\phi &=&
\int_{\partial U_X} j^1({\phi \circ \phi_X^{-1}})^*
[j^1(\xi) \intprod \Theta_\mathcal{L}] \\
&=& \int_{U_X} j^1({\phi \circ \phi_X^{-1}})^*
[j^1(\xi) \intprod \Omega_\mathcal{L}],
\end{eqnarray*}
since ${\mathfrak L}_{j^1(\xi)} \Theta_\mathcal{L} = 0$ 
by (\ref{s22a}) and Corollary \ref{cor4.2}.

We denote the covariant momentum map on $J^1(Y)$ by
$\mathcal{J}^\mathcal{L}\in L(\mathfrak{g}, \Lambda^n(J^1(Y))$
which we define as
\begin{equation}
j^1(\xi) \intprod \Omega_\mathcal{L} =
d \mathcal{J}^\mathcal{L}(\xi).
\label{s22}
\end{equation}

Using (\ref{s22}), we find that
$ \int_{U_X} d [ j^1({\phi \circ \phi_X^{-1}})^*
\mathcal{J}^\mathcal{L}(\xi)] =0$, and since this must hold
for all infinitesimal generators $\xi_\mathcal{C}$ at $\phi \in
\mathcal{C}$, the integrand must also vanish so that
\begin{equation}
d[ j^1({\phi \circ \phi_X^{-1}})^* \mathcal{J}^\mathcal{L}(\xi)]
=0,
\label{s23}
\end{equation}
which is precisely a restatement of the covariant Noether
Theorem.

\section{Veselov-type Discretizations 
of Multisymplectic Field Theory}
\subsection{General Theory}
We  now  generalize the Veselov discretization given in Section (3) to
multisymplectic field theory, by discretizing the
spacetime~$X$. For simplicity we restrict to the discrete
analogue of $\dim X=2$; i.e. $n=1$. Thus, we take
$X=\mathbb Z\times\mathbb Z=\{(i,j)\}$ and the fiber bundle $Y$
to be $X\times \mathcal F$ for some smooth manifold
$\mathcal F$.

\paragraph{Notation.}
The development in this section is aided by a small amount of
notation and terminology. Elements of $Y$ over the base point
$(i,j)$ are  written as $y_{ij}$ and the projection $\pi_{XY}$
acts on $Y$ by $\pi_{XY}(y_{ij})=(i,j)$. The fiber over
$(i,j)\in X$ is denoted $Y_{ij}$. A {\bfi triangle\/} $\Delta$
of $X$ is an  ordered triple of the form
\begin{displaymath}
    \Delta=\bigl((i,j),(i,j+1),(i+1,j+1)\bigr).
\end{displaymath}
The first component $(i,j)$ of $\Delta$ is the {\bfi first
vertex\/} of the triangle, denoted $\Delta^1$, and similarly
for the {\bfi second\/} and  {\bfi third vertices.} The set of
all triangles in $X$ is denoted $X^\Delta$. By abuse of
notation the same symbol is used for a triangle and the
(unordered) set of its vertices. A point $(i,j)\in X$ is {\bfi
touched\/} by a triangle if it is a vertex of that  triangle.
If $U\subseteq X$, then $(i,j)\in U$ is an {\bfi interior
point\/} of $U$ if  $U$~contains all three triangles of $X$
that touch $(i,j)$. The {\bfi interior\/} $\operatorname{int}U$
of $U$ is the collection of the interior points of $U$. The
{\bfi closure\/} $\operatorname{cl} U$ of $U$ is the union of
all triangles touching interior points of $U$. A {\bfi boundary
point\/} of $U$ is a point in $U$ and $\operatorname{cl} U$
which  is not an interior point. The {\bfi boundary\/} of $U$
is the set of boundary points of $U$, so that
\begin{displaymath}\partial U\equiv
(U\cap\operatorname{cl}U)\setminus
\operatorname{int}U
\end{displaymath} Generally, $U$  properly contains the
union of its interior and boundary, and we call $U$  {\bfi
regular\/} if it is exactly that union. A {\bfi section\/} of
$Y$ is a map $\phi:U\subseteq X\rightarrow Y$ such that
$\pi_{XY}\circ\phi=\operatorname{id}_U$.

\begin{figure}[ht]
\epsfxsize=0.9\textwidth
\centerline{\epsfbox{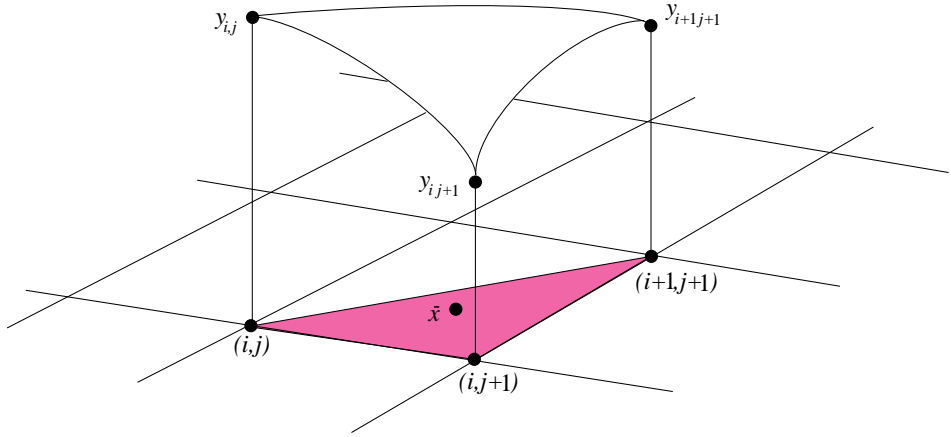}}
\caption{\label{25}
{\it\small Depiction of the heuristic interpretation of
an element of $J^1Y$ when $X$ is discrete.} \normalsize}
\end{figure}

\paragraph{Multisymplectic phase space.}
We define the  {\bfi first jet bundle}\footnote{
Using three vertices is the simplest choice for approximating
the two partial derivatives of the field $\phi$, but may
not lead to a good numerical scheme.  Later, we shall also use
four vertices together with averaging to define the partial
derivatives  of the fields.} 
of $Y$ to be
\begin{eqnarray*}
J^1Y&\equiv&\bigl\{(y_{ij},y_{i\,j+1},
  y_{i+1\,j+1})
\bigm|(i,j)\in X,\;
y_{ij},y_{i\,j+1},y_{i+1\,j+1}\in \mathcal F\bigr\}\\
&\equiv&X^\Delta\times \mathcal F^3.
\end{eqnarray*}
Heuristically (see Figure~(\ref{25})), $X$ corresponds to some
grid of elements $x_{ij}$ in continuous spacetime, say $\tilde
X$, and $\bigl(y_{ij},y_{i\,j+1},y_{i+1\,j+1}\bigr)\in J^1Y$
corresponds to $j^1\phi(\bar x)$, where $\bar x$ is ``inside''
the triangle bounded by $x_{ij},x_{i\,j+1},x_{i+1\,j+1}$, and
$\phi$ is some smooth section of $\tilde X\times \mathcal F$ 
interpolating the field values
$y_{ij},y_{i\,j+1},y_{i+1\,j+1}$. The {\bfi first jet
extension\/} of a section $\phi$ of $Y$ is the map
$j^1\phi:X^\Delta\rightarrow J^1Y$  defined by
\begin{displaymath}
      j^1\phi(\Delta)\equiv
   \bigl(\Delta,\phi(\Delta^1),
     \phi(\Delta^2),\phi(\Delta^3)\bigr).
\end{displaymath}
Given a vector field $Z$ on $Y$, we denote its restriction to
the fiber $Y_{ij}$ by $Z_{ij}$, and similarly for vector fields
on $J^1Y$. The {\bfi first jet extension\/} of a vector field
$Z$ on $Y$ is the vector field $j^1Z$ on $J^1Y$ defined by
\begin{displaymath}j^1Z(
  y_{\Delta^1},y_{\Delta^2},y_{\Delta^3})\equiv\big(
  Z_{\Delta^1}(y_{\Delta^1}),
  Z_{\Delta^2}(y_{\Delta^2}),
  Z_{\Delta^3}(y_{\Delta^3})\bigr),
\end{displaymath}
for any triangle $\Delta$.

\begin{figure}[ht]
\epsfxsize=0.5\textwidth
\centerline{\epsfbox{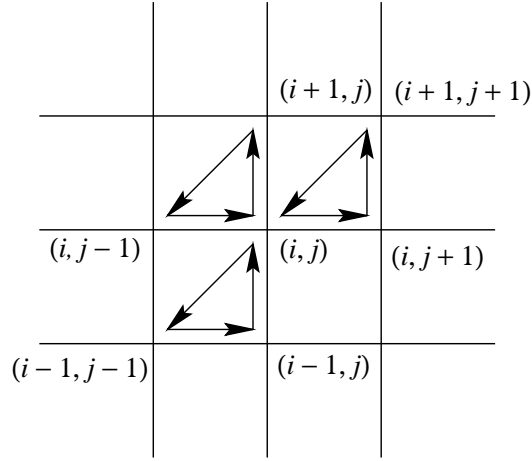}}
\caption{\label{26}
{\it\small The triangles which touch $(i,j)$.}\normalsize}
\end{figure}

\paragraph{The variational principle.}
Let us posit a {\bfi discrete Lagrangian}
$L:J^1Y\rightarrow\mathbb R$.
Given a triangle $\Delta$, define the function
$L_\Delta:\mathcal F^3\rightarrow\mathbb R$ by
\begin{displaymath}
L_\Delta(y_1,y_2,y_3)\equiv L(\Delta,y_1,y_2,y_3),
\end{displaymath}
so that we may view the Lagrangian $L$ as being a choice of
a function $L_\Delta$ for each triangle $\Delta$ of $X$. The
variables on the domain of $L_\Delta$ will be labeled
$y^1,y^2,y^3$, irrespective of the particular $\Delta$. Let $U$
be regular and let $\mathcal{C}_U$ be the set of sections of
$Y$ on  $U$, so $\mathcal{C}_U$ is the manifold $\mathcal
F^{|U|}$. The {\bfi action\/} will assign real numbers to
sections in $\mathcal{C}_U$ by the rule
\begin{equation}\label{30}
    S(\phi)\equiv\sum_{\Delta;\Delta\subseteq U}L\circ
       j^1\phi(\Delta).
\end{equation}
Given $\phi\in\mathcal{C}_U$ and a  vector field $V$, there is
the 1-parameter family of sections
\begin{displaymath}
    (F_\epsilon^{V}\phi)(i,j)\equiv F_\epsilon^{V_{ij}}
        (\phi(i,j)),
\end{displaymath}
where $F^{V_{ij}}$ denotes the flow of $V_{ij}$ on $\mathcal F$.
The {\bfi variational principle\/} is to seek those $\phi$ for
which
\begin{displaymath}
    {\displaystyle \left. \frac{d}{d\epsilon }\right|_{\epsilon
       = 0}} S(F_\epsilon^V\phi)=0
\end{displaymath}
for all  vector fields $V$.

\paragraph{The discrete Euler-Lagrange equations.}
The  variational principle gives certain field equations, the
{\bfi discrete Euler-Lagrange field equations} (DELF
equations), as follows. Focus upon some
$(i,j)\in\operatorname{int}U$, and abuse notation by writing
$\phi(i,j)\equiv y_{ij}$. The action, written with its summands
containing $y_{ij}$ explicitly, is (see Figure~(\ref{26}))
\begin{displaymath}
S=\cdots+L(y_{ij},y_{i\,j+1},y_{i+1\,j+1})
+L(y_{i\,j-1},y_{ij},y_{i+1\,j})
+L(y_{i-1\,j-1},y_{i-1\,j},y_{ij})+\cdots
\end{displaymath}
so by differentiating in $y_{ij}$, the DELF equations are
\begin{displaymath}
 \frac{\partial L}{\partial y^1}(y_{ij},y_{i\,j+1},y_{i+1\,j+1})
+
 \frac{\partial L}{\partial y^2}(y_{i\,j-1},y_{ij},y_{i+1\,j})
+
  \frac{\partial L}{\partial
       y^3}(y_{i-1\,j-1},y_{i-1\,j},y_{ij})=0,
\end{displaymath}
for all $(i,j)\in\operatorname{int}U$. Equivalently, these
equations may be written
\begin{equation}\label{41}
    \sum_{l;\Delta;(i,j)=\Delta^l}
       \frac{\partial L_\Delta}{\partial y^l}(
          y_{\Delta^1},y_{\Delta^2},y_{\Delta^3})=0,
\end{equation}
for all $(i,j)\in\operatorname{int}U$.

\paragraph{The discrete Cartan form.}
Now suppose we allow nonzero variations on the boundary
$\partial U$, so we consider the effect on $S$ of a vector
field $V$ which does not necessarily vanish on $\partial U$.
For each $(i,j)\in\partial U$ find the triangles in $U$
touching $(i,j)$. There is at least one such triangle since
$(i,j)\in\operatorname{cl}U$; there are not three such
triangles since $(i,j)\not\in\operatorname{int}U$. For each
such triangle $\Delta$, $(i,j)$ occurs as the
$l^{\mbox{\scriptsize th}}$ vertex, for one or two  of 
$l=1,2,3$, and  those
$l^{\mbox{\scriptsize th}}$ expressions from the list
\begin{eqnarray*}
&&\displaystyle
    \frac{\partial L}{\partial
        y^1}(y_{ij},y_{i\,j+1},y_{i+1\,j+1})V_{ij}
           (y_{ij}),\\
&&\displaystyle
 \frac{\partial L}{\partial y^2}
       (y_{i\,j-1},y_{ij},y_{i+1\,j})V_{ij}(y_{ij}),\\
&&\displaystyle
\frac{\partial L}{\partial y^3}
       (y_{i-1\,j-1},y_{i-1\,j},y_{ij})V_{ij}(y_{ij}),
\end{eqnarray*}
yielding one or two  numbers. The contribution to $dS$ from the
boundary is the sum of  all such numbers. To bring this into a
recognizable format, we take our cue from discrete Lagrangian
mechanics, which featured {\em two} $1$-forms. Here the above
list suggests the {\em three} $1$-forms on $J^1Y$, the first of
which we define to be
\begin{eqnarray*}
   &&\Theta_L^1(y_{ij},y_{i\,j+1},y_{i+1\,j+1}) \cdot
     (v_{y_{ij}},v_{y_{i\,j+1}},
      v_{y_{i+1\,j+1}})\\ &&\hspace*{1.5in} \equiv\frac{\partial
        L}{\partial y^1}(y_{ij},y_{i\,j+1},y_{i+1\,j+1})
         \cdot(v_{y_{ij}},0,0),
\end{eqnarray*}
$\Theta_L^2$ and $\Theta_L^3$ being defined analogously. With
these notations, the contribution to $dS$ from the boundary can
be written  $\theta_L(\phi)\cdot V$, where
$\theta_L$ is the $1$-form on the space of  sections
$\mathcal{C}_U$ defined by
\begin{equation}\label{31}
\theta_L(\phi)\cdot V
  \equiv\sum_{\Delta;\Delta\cap\partial U\not=\emptyset}\left(
      \sum_{l;\Delta^l\in\partial U}
  \left[(j^1\phi)^*(j^1V\intprod\Theta^l_L)\right](\Delta)\right).
\end{equation}
In comparing~(\ref{31}) with~(\ref{s21}), the analogy with the
multisymplectic formalism of Section (4) is immediate.
\paragraph{The discrete multisymplectic form formula.} Given  a
triangle  $\Delta$ in $X$, we define the projection
$\pi_\Delta:\mathcal{C}_U\rightarrow J^1Y$ by
\begin{displaymath}
\pi_\Delta(\phi)\equiv 
   (\Delta,y_{\Delta^1},y_{\Delta^2},y_{\Delta^3}).
\end{displaymath}
In this notation, it is easily verified that~(\ref{31}) takes
the convenient form
\begin{equation}\label{43}
\theta_L
  =\sum_{\Delta;\Delta\cap\partial U\not=\emptyset}
 \left(\sum_{l;\Delta^l\in\partial U}
  \pi_\Delta^*\Theta^l_L\right).
\end{equation}
A {\bfi first-variation} at a solution  $\phi$ of the DELF
equations is a vertical vector field $V$ such that the
associated flow $F^V$ maps $\phi$ to other solutions of the
DELF equations. Set $\Omega^l_L=-d\Theta_L^l$. Since
\begin{equation}\label{35}
\Theta_L^1+\Theta_L^2+\Theta_L^3=dL,
\end{equation}
one obtains
\begin{displaymath}
\Omega_L^1+\Omega_L^2+\Omega_L^3=0,
\end{displaymath}
so that only two of the three $2$-forms $\Omega_L^l$, $l=1,2,3$
are essentially distinct. Exactly as in Section~(2), the
equation $d^2S=0$, when specialized to two first-variations $V$
and $W$ now gives, by taking one exterior derivative
of~(\ref{43}),
\begin{displaymath}
0 = d\theta_L(\phi)(V,W)
  = \sum_{\Delta;\Delta\cap\partial U\not=\emptyset}
    \left(\sum_{l;\Delta^l\in\partial U}
     V \intprod W \intprod \pi_ \Delta^* \Omega^l_L \right),
\end{displaymath}
which in turn is equivalent to
\begin{equation}\label{32}
\sum_{\Delta;\Delta\cap\partial U\not=\emptyset}
\left(\sum_{l;\Delta^l\in\partial U}
  \left[(j^1\phi)^*(j^1V\intprod
j^1W\intprod\Omega^l_L)\right](\Delta)\right)=0.
\end{equation}
Again, the analogy with the  multisymplectic form formula for
continuous spacetime~(\ref{s11}) is immediate.
\paragraph{The discrete Noether theorem.} Suppose that a Lie
group $G$ with Lie algera ${\mathfrak g}$ acts on $F$ by vertical 
symmetries in such a way that the Lagrangian~$L$ is
$G$-invariant. Then $G$ acts on $Y$ and $J^1Y$ in the
obvious ways. Since there are three Lagrange $1$-forms, there
are three momentum maps $J^l$, $l=1,2,3$,  each one a
$\mathfrak{g}^*$-valued function on triangles in $X$, and
defined by
\begin{displaymath}J^l_\xi\equiv\xi_{J^1Y}\intprod\Theta_L^l,
\end{displaymath}
for any $\xi\in\mathfrak{g}$.
Invariance of $L$ and~(\ref{35}) imply that
\begin{displaymath}
    J^1+J^2+J^3=0,
\end{displaymath}
so, as in the case of the $1$-forms, only two of the three
momenta are essentially distinct. For any $\xi$, the
infinitesimal generator $\xi_Y$ is a first-variation, so
invariance of $S$, namely $\xi_Y\intprod dS=0$ , becomes
$\xi_Y\intprod\theta_L=0$. By left insertion into~(\ref{31}),
this becomes the discrete version of Noether's  theorem:
\begin{equation}\label{38}
    \sum_{\Delta;\Delta\cap\partial
      U\not=\emptyset}\left(\sum_{l;\Delta^l\in\partial U}
        J^l(\Delta)\right)=0.
\end{equation}

\begin{figure}[ht]
\epsfxsize=0.7\textwidth
\centerline{\epsfbox{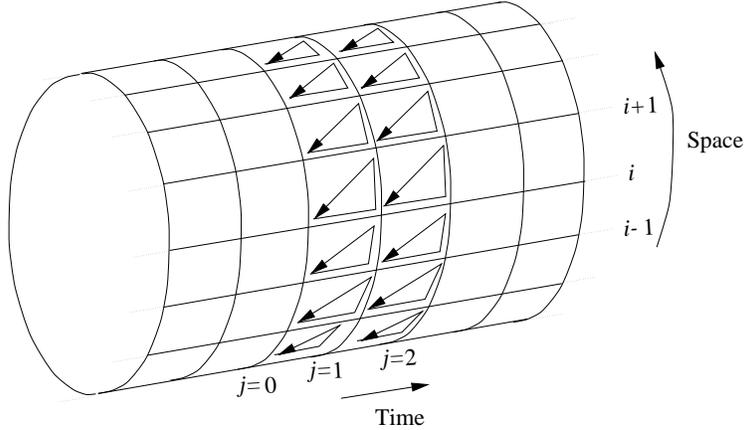}}
\caption{\label{39}
{\it \small Symplectic flow and conservation of momentum  from the
discrete Noether theorem when the spatial boundary is empty and
the temporal boundaries agree.}\normalsize}
\end{figure}

\paragraph{Conservation in a space and time split.}
To understand the significance of~(\ref{32}) and~(\ref{38})
consider a discrete field theory with space a discrete version
of the circle and time the real line, as depicted
in~Figure~(\ref{39}), where space is split into space and
time, with ``constant time'' being constant $j$ and the
``space index'' $1\le i\le N$ being cyclic.
Applying~(\ref{38}) to the region $\{(i,j)\mid j=0,1,2\}$
shown in the Figure, Noether's theorem takes the conservation
form
\begin{eqnarray*}
\sum_{i=1}^NJ^1(y_{i0},y_{i1},y_{i+1\,1})&=&
  -\sum_{i=1}^N\bigl(J^2(y_{i1},y_{i2},y_{i+1\,2})+
  J^3(y_{i1},y_{i2},y_{i+1\,2})\bigr)\\
&=&\sum_{i=1}^NJ^1(y_{i1},y_{i2},y_{i+1\,2}).
\end{eqnarray*}
Similarly, the discrete multisymplectic form formula also
takes a conservation form. When there is spatial boundary, the
discrete Noether theorem and the discrete multisymplectic form
formulas automatically account for it, and thus form
nontrivial generalizations of these conservation results.

Furthermore, as in the continuous case, we can achieve
``evolution type'' symplectic systems (i.e. discrete
Moser-Veselov mechanical systems) if we define $Q$ as the
space of fields at constant $j$, so $Q\equiv \mathcal F^N$, and
take as the discrete Lagrangian
\begin{displaymath}
\tilde L([q_j^0],[q_j^1])\equiv\sum_{i=1}^N
L(q_i^0,q_i^1,q_{i+1}^1).
\end{displaymath}
Then the  Moser-Veselov DEL evolution-type equations~(\ref{40}) are
equivalent to the DELF equations~(\ref{41}), the
multisymplectic form formula implies symplecticity of the
Moser-Veselov evolution map, and conservation of momentum
gives identical results in both the ``field'' and
``evolution'' pictures.

\paragraph{Example: nonlinear wave equation.}
To illustrate the discretization method we have developed, let
us consider the Lagrangian~(\ref{lag}) of Section~(4), which describes
the nonlinear sine-Gordon wave equation. This is a completely integrable
system with an extremely interesting hierarchy of soliton solutions,
which we shall investigate by developing for it a variational 
multisymplectic-momentum integrator; see the recent article by Palais
[1997] for a wonderful discussion on soliton theory. 

To discretize the continuous Lagrangian, we visualize each
triangle $\Delta$ as having base length~$h$ and height~$k$, and
we think of the discrete jet
$(y_{\Delta^1},y_{\Delta^2},y_{\Delta^3})$ as corresponding to
the continuous jet
\begin{displaymath}
\frac{\partial\phi}{\partial x^0}(\bar y_{ij})=\frac{y_{i\,j+1}-y_{ij}}h,\quad
\frac{\partial\phi}{\partial x^1}(\bar y_{ij})=\frac{y_{i+1\,j+1}-y_{i\,j+1}}k,
\end{displaymath}
where $\bar y_{ij}$ is a the center of the triangle
\footnote{
Other discretizations based on triangles are possible; for
example, one could use the  value $y_{ij}$ for insertion into
the nonlinear term instead of $\bar y_{ij}$.}. This leads to
the discrete Lagrangian
\begin{displaymath}
L=\frac12\left(\frac{y_2-y_1}h\right)^2-
\frac12\left(\frac{y_3-y_2}k\right)^2
+N\left(\frac{y_1+y_2+y_3}{3}\right),
\end{displaymath}
with corresponding DELF equations
\begin{eqnarray}
\lefteqn{\frac {y_{{i+1\,j}}-2\,y_{{ij}}+y_{{i-1\,j}}}{{k}^{2}}
-\frac {y_{{i\,j+1}}-2\,y_{{ij}}+y_{{i\,j-1}}}{{h}^{2}}}&& \nonumber\\
&&\qquad\mbox{}+
\frac13N^\prime\left(\frac{y_{{ij}}+y_{{i\,j+1}}
+y_{{i+1\,j+1}}}3\right) \nonumber\\
&&\qquad\qquad\mbox{}+
\frac13N^\prime\left(\frac{y_{{i\,j-1}}+y_{{ij}}+y_{{i+1\,j}}}3\right) \nonumber\\
&&\qquad\qquad\qquad\mbox{}
+\frac13N^\prime\left(\frac{y_{{i-1\,j-1}}+y_{{i-1\,j}}+y_{{ij}}}3\right)=0.
\label{explicit}
\end{eqnarray}
When $N=0$ (wave equation) this gives the explicit method
\begin{displaymath}
y_{i\,j+1}=\frac{h^2}{k^2}
(y_{{i+1\,j}}-2\,y_{{ij}}+y_{{i-1\,j}})+2\,y_{{ij}}-y_{{i\,j-1}},
\end{displaymath}
which is stable whenever the Courant stability condition is
satisfied.

\begin{figure}[ht]
\epsfxsize=0.4\textwidth
\centerline{\epsfbox{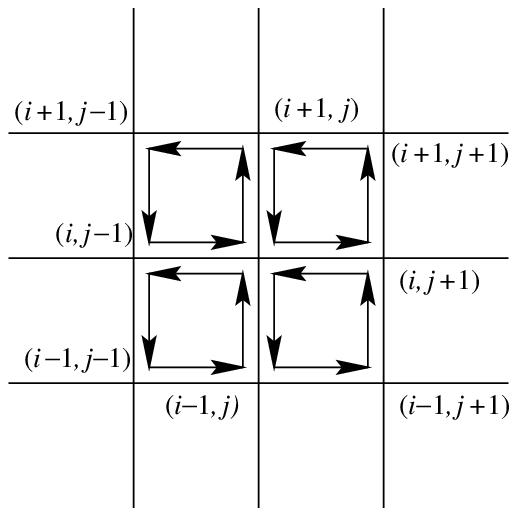}\hspace*{.3in}\epsfbox{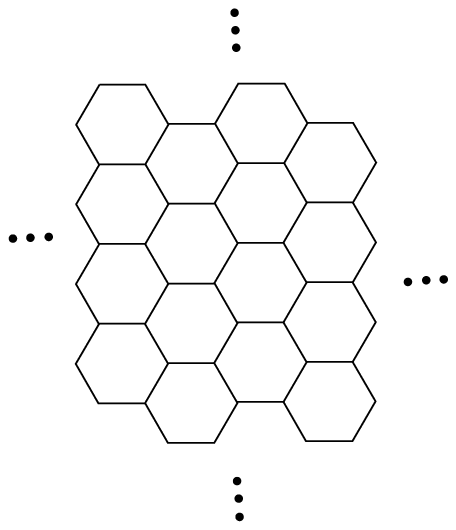}}
\caption{\label{60} 
{\it \small On the left, the method 
based on rectangles; on the right,
a possible method based on hexagons.}\normalsize}
\end{figure}

\paragraph{Extensions: Jets from rectangles and other polygons.}
Our choice of discrete jet bundle is obviously not restricted to
triangles, and can be extended to rectangles or more general
polygons (left of Figure(~\ref{60})).  A {\bfi rectangle} is a
quadruple of the form,
\begin{displaymath}
\Delta=\bigl((i,j),(i,j+1),(i+1,j+1),(i+1,j)\bigr),
\end{displaymath}
a point is an {\bfi interior point} of a subset~$U$ of
rectangles if 
$U$ contains all four rectangles touching that point, the
discrete Lagrangian depends on variables $y_1,\cdots ,y_4$, 
and the DELF equations become
\begin{eqnarray*}
\lefteqn{\frac{\partial L}{\partial y^1}(y_{ij},y_{i\,j+1},y_{i+1\,j+1},
  y_{i+1\,j})
  +\frac{\partial L}{\partial y^2}(y_{i\,j-1},y_{ij},y_{i+1\,j},
  y_{i+1\,j-1})}&&\\
&&\qquad\mbox+
\frac{\partial L}{\partial y^3}(y_{i-1\,j-1},y_{i-1\,j},y_{ij},y_{i\,j-1})
+\frac{\partial L}{\partial y^4}(y_{i-1\,j},y_{i-1\,j+1},y_{i\,j+1},y_{ij})=0.
\end{eqnarray*}
The extension to polygons with even higher numbers of sides is
straightforward; one example is illustrated on the right of
Figure(~\ref{60}).
The motivation for consideration of these extensions is
enhancing the stability of the triangle-based method in
the nonlinear wave example just above.

\paragraph{Example: nonlinear wave equation, rectangles.}
Think of each rectangle $\Delta$ as having length~$h$ and
height~$k$, and each discrete jet
$(y_{\Delta^1},y_{\Delta^2},y_{\Delta^3},y_{\Delta^4})$ being
associated to the continuous jet
\begin{displaymath}
\frac{\partial\phi}{\partial x^0}(p)=
\frac{y_{i\,j+1}-y_{ij}}h
,\quad
\frac{\partial\phi}{\partial x^1}(p)=\frac12\left(
\frac{y_{i+1\,j}-y_{i\,j}}k
   +\frac{y_{i+1\,j+1}-y_{i\,j+1}}k\right),
\end{displaymath}
where $p$ is a the center of the rectangle. This leads to the discrete
Lagrangian
\begin{eqnarray}
L&=&\frac12\left(\frac{y_2-y_1}h\right)^2
-\frac12\left(\frac{y_4-y_1}{2k}
  +\frac{y_3-y_2}{2k}\right)^2\nonumber\\
&&\qquad\mbox{}+N\left(\frac{y_1
    +y_2+y_3+y_4}{4}\right).\label{100}
\end{eqnarray}
If, for brevity, we set
\begin{displaymath}
\bar y_{ij}\equiv\frac{y_{ij}+y_{i\,j+1}+y_{i+1\,j+1}
    +y_{i+1\,j}}4,
\end{displaymath}
then one verifies  that the DELF equations become
\begin{eqnarray*}
&&\Biggl[\frac12\frac {y_{{i+1\,j}}-2\,y_{{ij}}
   +y_{{i-1\,j}}}{{k}^{2}}
   +\frac14\frac {y_{{i+1\,j+1}}-2\,y_{{i\,j+1}}
   +y_{{i-1\,j+1}}}{{k}^{2}}\\
&&\qquad\mbox{}+\frac14\frac {y_{{i+1\,j-1}}-2\,y_{{i\,j-1}}
  +y_{{i-1\,j-1}}}{{k}^{2}}\Biggr]
    -\Biggl[\frac {y_{{i\,j+1}}
    -2\,y_{{ij}}+y_{{i\,j-1}}}{{h}^{2}}\Biggr]\\
&&\qquad\qquad\mbox{}
+\frac14\Biggl[N^\prime(\bar y_{ij})+
N^\prime(\bar y_{i\,j-1})+
N^\prime(\bar y_{i-1\,j-1})+N^\prime(\bar y_{i-1\,j})
\Biggr]=0,
\end{eqnarray*}
which, if we make the definitions
\begin{eqnarray*}
&\displaystyle \partial^2_hy_{ij}\equiv
  y_{{i\,j+1}}-2\,y_{{ij}}+y_{{i\,j-1}},\quad
 \partial^2_ky_{ij}\equiv
 y_{{i+1\,j}}-2\,y_{{ij}}+y_{{i-1\,j}},&\\
&\displaystyle \bar N^\prime(\bar y_{ij})\equiv
  \frac14\Biggl[N^\prime(\bar y_{ij})+
  N^\prime(\bar y_{i\,j-1})+
  N^\prime(\bar y_{i-1\,j-1})
  +N^\prime(\bar y_{i-1\,j})\Biggr],&
\end{eqnarray*}
is (more compactly)
\begin{equation}
\frac1{k^2}\Biggl[\frac14\partial^2_ky_{i\,j+1}
  +\frac12\partial^2_ky_{ij}
  +\frac14\partial^2_ky_{i\,j-1}\Biggr]
-\frac1{h^2}\partial^2_hy_{ij} 
+ \bar N^\prime(\bar y_{ij})=0.\label{101}
\end{equation}
These are implicit equations which must  be solved for 
$y_{i\,j+1}$,
$1\le i\le N$, given  $y_{i\,j}$, $y_{i\,j-1}$, $1\le i\le N$;
rearranging, an iterative form equivalent to~(\ref{101}) is
\begin{eqnarray*}
\lefteqn{-\left(\frac{h^2}{2(h^2+2k^2)}\right)
   y_{i+1\,j+1}+y_{i\,j+1}-
  \left(\frac{h^2}{2(h^2+2k^2)}\right)y_{i-1\,j+1}}&&\\
&&\qquad\mbox{}=\,
\frac{h^2}{h^2+2k^2}
\Bigl((y_{{i+1\,j}}-2\,y_{{ij}}+y_{i-1\,j})\\
&&\qquad\qquad\qquad\qquad\mbox{}
+\frac12(y_{{i+1\,j-1}}-2\,y_{{i\,j-1}}+y_{{i-1\,j-1}})\Bigr)\\
&&\qquad\quad\mbox{}+\frac{2k^2}{h^2+2k^2}
\left(2\,y_{{ij}}-y_{{i\,j-1}}\right)\\
&&\qquad\quad\mbox{}+\frac{h^2k^2}{2(h^2+2k^2)}
\bigl(N^\prime(\bar y_{ij})+
N^\prime(\bar y_{i\,j-1})+
N^\prime(\bar y_{i-1\,j-1})+N^\prime(\bar y_{i-1\,j})
\bigr).
\end{eqnarray*}
In the case of the sine-Gordon equation the values of the field
ought to be considered as lying in ${\mathbb S}^1$, by virtue of the
vertical symmetry  $y\mapsto y+2\pi$.  Soliton solutions for
example will have a jump of $2\pi$ and  the method will fail
unless field values at close-together spacetime points  are
differenced modulo $2\pi$. As a result it becomes important to
calculate using integral multiples of small field-dependent
quantities, so that it is clear when to discard  multiples of
$2\pi$, and for this the above iterative form is inconvenient.
But if we define
\begin{displaymath}
\partial^1_hy_{ij}\equiv y_{{i\,j+1}}-\,y_{{ij}},\qquad
\partial^1_ky_{ij}\equiv y_{{i+1\,j}}-\,y_{{ij}},
\end{displaymath}
then there is the following iterative form, 
again equivalent to~(\ref{101})
\begin{eqnarray}
&&
  y_{i\,j+1}=y_{ij} 
  + \partial^1_hy_{i\,j},\quad\mbox{and}\nonumber\\\nonumber\\
&&
  -\left(\frac{h^2}{2(h^2+2k^2)}\right)\partial^1_hy_{i+1\,j}
  +\partial^1_hy_{i\,j}-
  \left(\frac{h^2}{2(h^2+2k^2)}\right)\partial^1_hy_{i-1\,j}
  \nonumber\\
&&\qquad\mbox{}=\,\frac{h^2}{h^2+2k^2}(
    3\partial^2_ky_{ij}+\partial^2_ky_{i\,j-1})
    +\frac{2k^2}{h^2+2k^2}\partial^1_hy_{ij}\nonumber\\
&&\qquad\qquad\mbox{}
    +\frac{h^2k^2}{2(h^2+2k^2)}
    \bar N^\prime(\bar y_{ij}).\label{110}
\end{eqnarray}

One can also modify~(\ref{100}) so as to treat
space and time symmetrically, 
which leads to the discrete Lagrangian
\begin{eqnarray*}
L&=&\frac12\left(\frac{y_2-y_1}{2h}+\frac{y_3-y_4}{2h}\right)^2
-\frac12\left(\frac{y_4-y_1}{2k}+\frac{y_3-y_2}{2k}\right)^2\\
&&\qquad+N\left(\frac{y_1+y_2+y_3+y_4}{4}\right),
\end{eqnarray*}
and one verifies that the DELF equations become
\begin{eqnarray}
&&\frac1{k^2}\biggl[\frac14\partial^2_ky_{i\,j+1}
+\frac12\partial^2_ky_{ij}
  +\frac14\partial^2_ky_{i\,j-1}\biggr]\nonumber\\
&&\qquad\mbox{}-\frac1{h^2}\biggl[
\frac14\partial^2_hy_{i+1\,j}+\frac12\partial^2_hy_{ij}+
\frac14\partial^2_hy_{i-1\,j}\biggr] 
+ \bar N^\prime(\bar y_{ij})=0,\label{102}
\end{eqnarray}
an equivalent iterative form of which is
\begin{eqnarray}
&&
  y_{i\,j+1}=y_{ij} 
  + \partial ^1_hy_{i\,j},\quad\mbox{and}\nonumber\\
&&
  -\left(\frac{h^2-k^2}{2(h^2+k^2)}\right)\partial^1_hy_{i+1\,j}
  +\partial^1_hy_{i\,j}-
  \left(\frac{h^2-k^2}{2(h^2+k^2)}\right)\partial^1_hy_{i-1\,j}\nonumber\\
&&\qquad\mbox{}=\,\frac{h^2}{2(h^2+k^2)}(
    3\partial^2_ky_{ij}+\partial^2_ky_{i\,j-1})\nonumber\\
&&\qquad\qquad\mbox{}
 +\frac{h^2}{2(h^2+k^2)}(
 2\partial^1_hy_{ij}
 +\partial^1_hy_{i+1\,j}+\partial^1_hy_{i-1\,j})\nonumber\\
&& \qquad\qquad\qquad\mbox{}   
+\frac{h^2k^2}{2(h^2+k^2)}\bar N^\prime(\bar y_{ij}).\label{111}
\end{eqnarray}

\subsection{Numerical checks.}
\begin{figure}[ht]
\vspace*{-.25in}\centerline{
\epsfxsize=2.5in\epsfbox{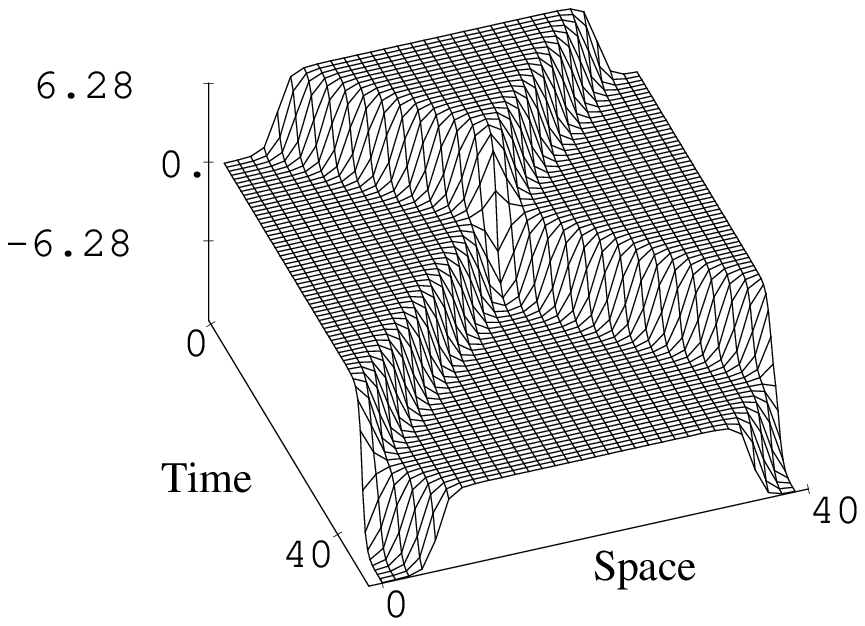}\hspace*{-.4in}
\epsfxsize=2.5in\epsfbox{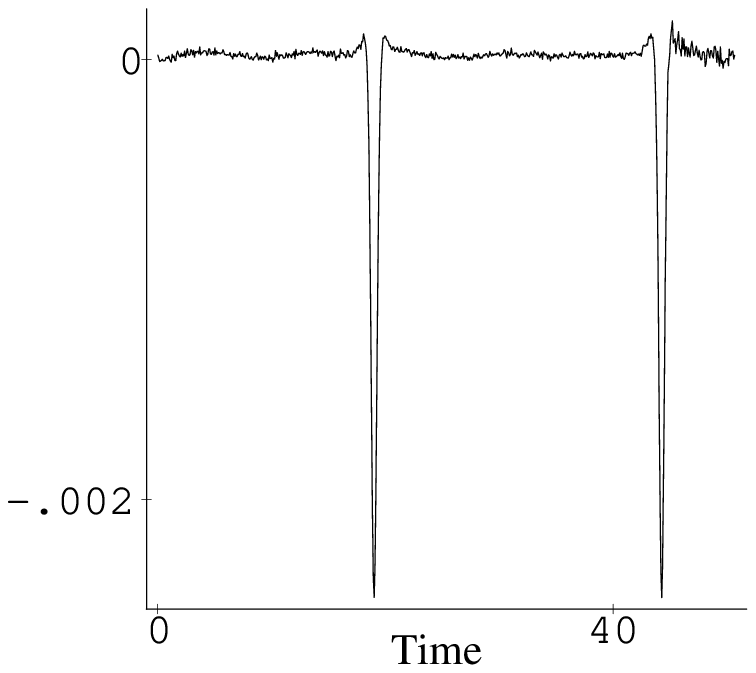}}\vspace*{-.6in}
\centerline{
\epsfxsize=2.5in\epsfbox{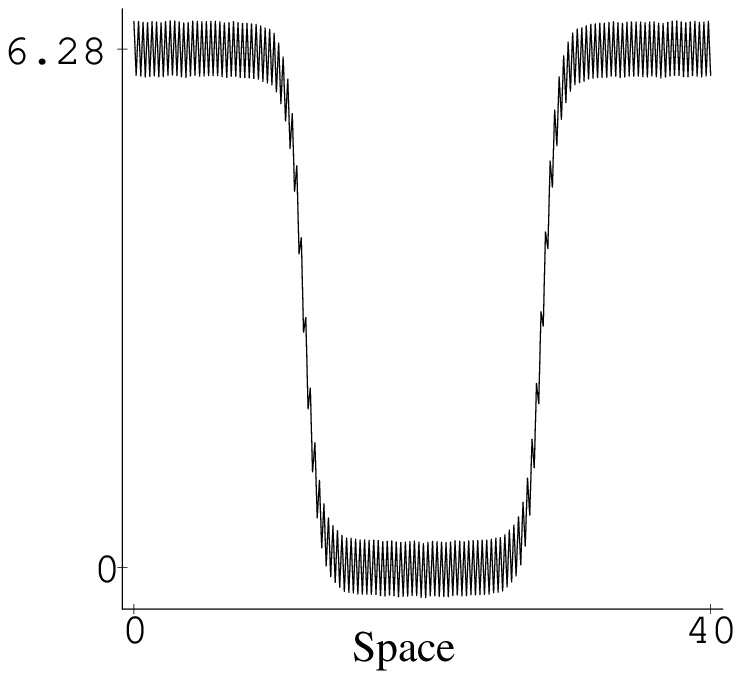}\hspace*{-.2in}
\epsfxsize=2.5in\epsfbox{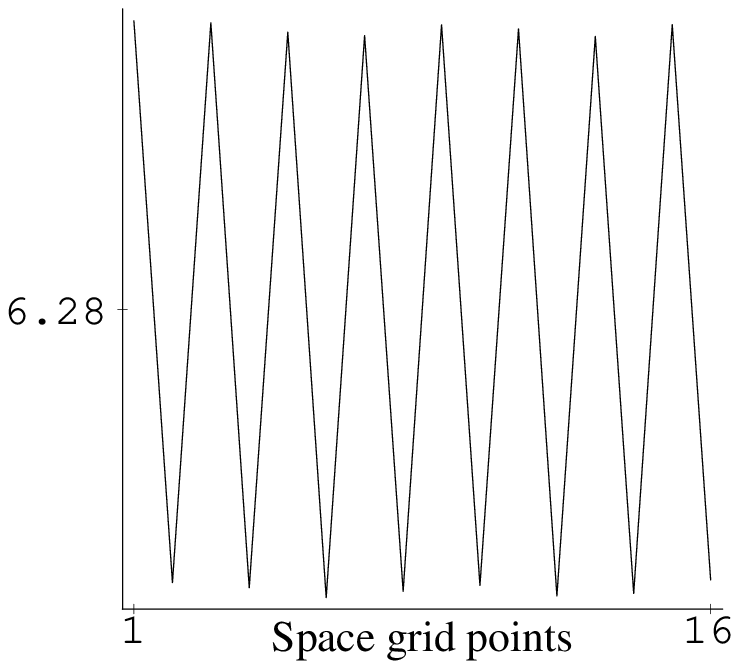}}\vspace*{-.25in}
\caption{\label{200} 
{\it \small Top left: the wave forms for a two soliton kink and 
antikink collision using~(\ref{102}).
Top right: the energy error.
Bottom left: the  wave form  at time
$t\approx11855$. Bottom right: the portion of the bottom left graph 
for spatial grid points $1\ldots 16$.}\normalsize}
\end{figure}
While the focus of this article is not the numerical
implementation of the  integrators which we have derived, we
have, nevertheless, undertaken some preliminary numerical 
investigations of our multisymplectic methods in the context of the
sine-Gordon equation with periodic boundary conditions.

\begin{figure}[p]
\vspace*{-.25in}\centerline{
\epsfxsize=2.5in\epsfbox{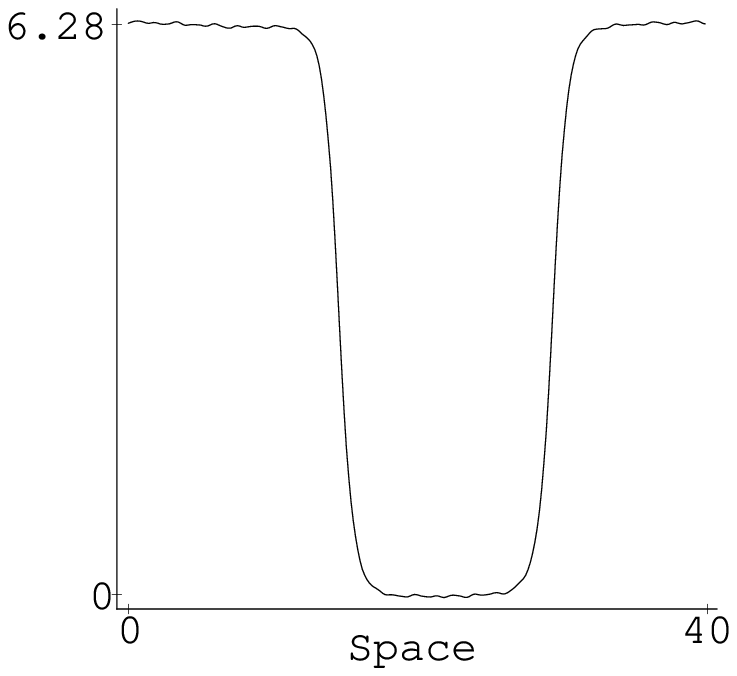}\hspace*{-.3in}
\epsfxsize=2.5in\epsfbox{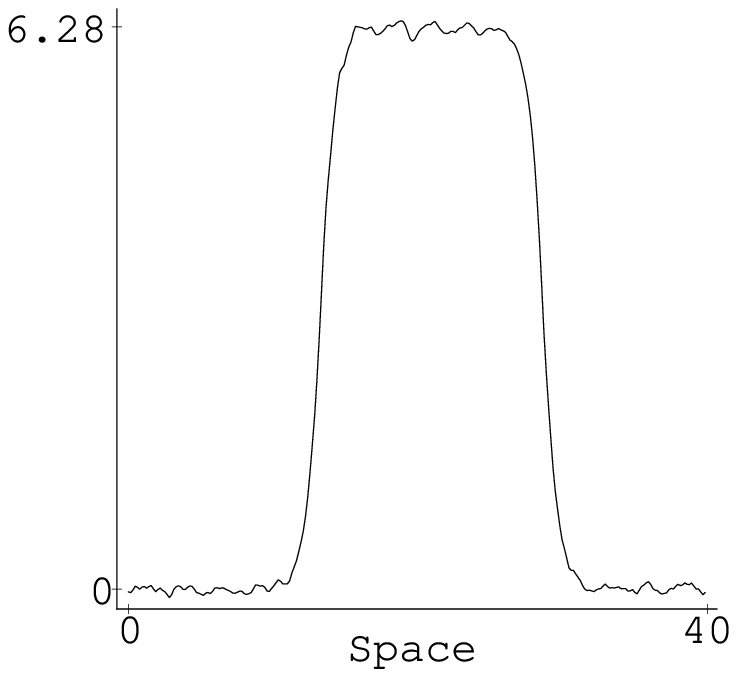}}\vspace*{-.25in}
\caption{\label{202} {\it \small On the left, the  final wave form (after about 
5000~soliton collisions at $t\approx129133$) obtained using the
rectangle-based multisymplectic
method~(\ref{102}). 
On the right, the final waveform (at $t\approx129145$)
from the  energy-conserving method~(\ref{155}) of Vu-Quoq and Li.
In both simulations, temporal drift is occurring.
For this reason the waveforms are inverted with respect to one 
another; moreover,
the separate solitons are drifting  at slightly
different rates, as indicated by the off-center waveforms.}
\normalsize }
\end{figure}

\begin{figure}[p]
\vspace*{-.25in}\centerline{
\epsfxsize=2.5in\epsfbox{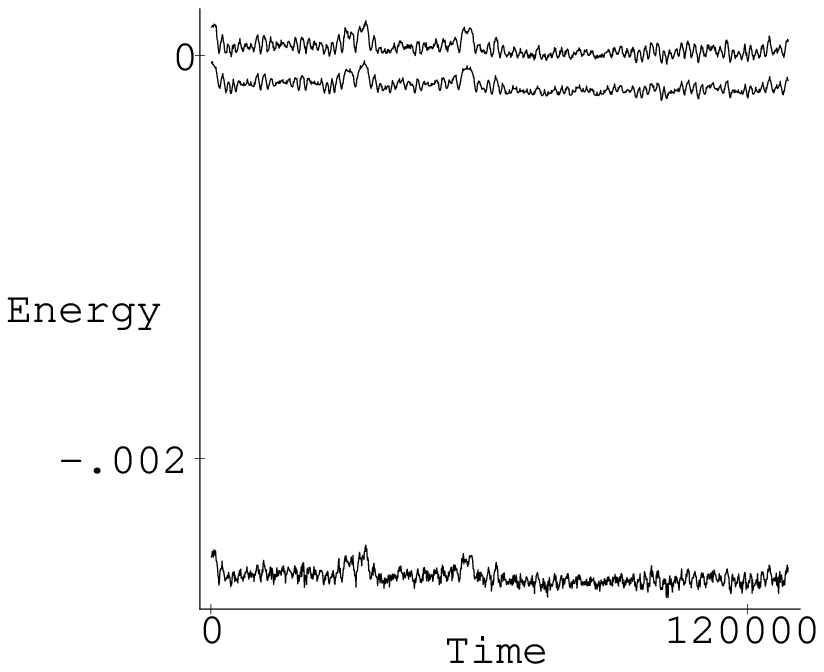}\hspace*{-.3in}
\epsfxsize=2.5in\epsfbox{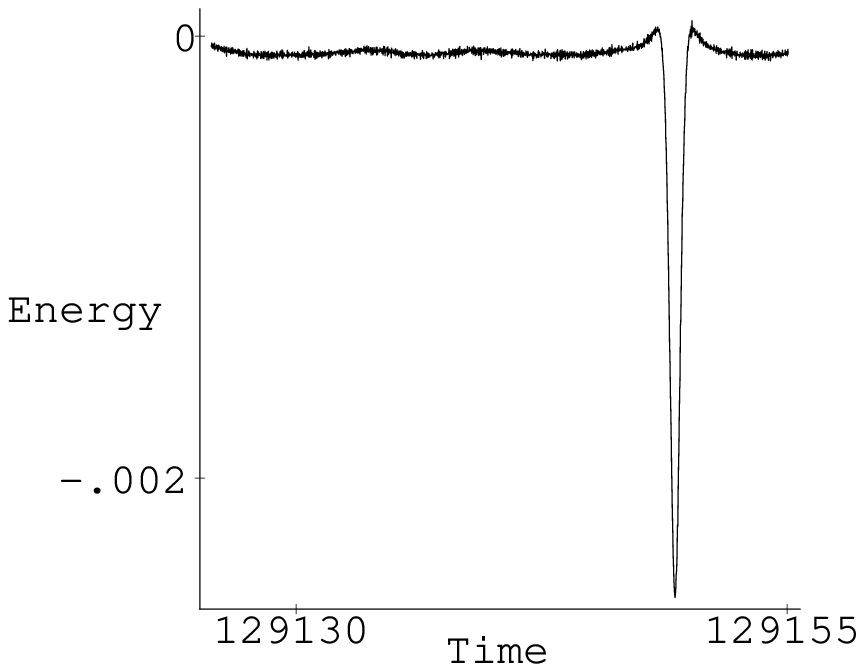}}\vspace*{-.25in}
\caption{\label{203} {\it \small On the left,
the energy error corresponding to our multisymplectic method (\ref{111})
for 5000 solition collisions; the three graphs correspond to the 
minimum, average, and maximum energy error over consecutive 5000 time
step regions. On the right, the final energy error (i.e. the energy
error after about 5000 soliton collisions), which can be compared with the 
initial energy error plot in the top left of~Figure~(\ref{200}).}\normalsize}
\end{figure}

\paragraph{The rectangle-based multisymplectic method.}
The top half of~Figure~(\ref{200})
shows a simulation of the collision of  ``kink''
and ``antikink'' solitons for the sine-Gordon 
equation, using the rectangle-based multisymplectic method~(\ref{102}).
In the bottom half of that figure we show the result of 
running that simulation until
the solitons have undergone about $460$ collisions; shortly
after this the simulation stops because the
iteration~(\ref{111}) diverges. The anomalous spatial
variations in the waveform of the bottom left of Figure~(\ref{200}) have period
$2$ spatial grid divisions and are shown in finer scale on the
bottom right of that figure. These variations are reminiscent of those
found in Ablowitz, Herbst and Schober~[1996] for the completely
integrable discretization of Hirota, where the
variations are attributed to independent evolution of waveforms
supported on even vs. odd grid points. Observation
of~(\ref{102}) indicates what is wrong: the nonlinear term $N$ 
contributes to~(\ref{102}) in a way that will average out these
variations, and consequently,  once they have begun,~(\ref{102}) 
tends to continue such variations  via the linear wave
equation. In  Ablowitz et. al.,  the  situation is rectified
when the number of  spatial grid points is not even,  and this
is the case for~(\ref{102}) as well. This is indicated on the
left of  Figure~(\ref{202}), which shows the waveform  after 
about~$5000$ soliton collisions when $N=255$ rather than
$N=256$.  Figure~(\ref{203})  summarizes the evolution of
energy 
error\footnote{The discrete energy that we calculated was
\begin{eqnarray*}
&&\sum_{i=1}^N\Biggl(
\frac12\left(\frac{y_{i\,j+1}-y_{ij}}{2h}
+\frac{y_{i+1\,j+1}-y_{i+1\,j}}
  {2h}\right)^2\\
&&\qquad\qquad\mbox{}
+\frac12\left(\frac{y_{i+1\,j}-y_{ij}}{2k}
+\frac{y_{i+1\,j+1}-y_{i\,j+1}}{2k}\right)^2
-N\left(\bar y_{ij}\right)\Biggr).
\end{eqnarray*}}
for that simulation.

\paragraph{Initial data.}

For the two-soliton-collision simulations,
we used the following initial data:
$h=k/8$ (except $h=k/16$ where noted), where $k=40/N$ and $N=255$
spatial grid points  (except Figure~(\ref{200}) where $N=256$).
The circle that is space should be visualized as having
circumference~$L=40$. Let $\kappa=1-\epsilon$ where
$\epsilon=10^{-6}$, $\tilde L=L/4=10$,
\begin{displaymath}
P=2\int_0^{1/\kappa}\!\!\!\!\frac{1}{\sqrt{1-y^2}\sqrt{1-\kappa^2y^2}}\,dy
\approx15.90,
\qquad c=\sqrt{1-\frac{\tilde L^2}{\kappa^2P^2}}\approx.7773,
\end{displaymath}
and 
\begin{displaymath}
\tilde\phi(x)\equiv2\arcsin
\left(\operatorname{sn}\left(\frac{x}{\kappa
\sqrt{1-c^2}};\kappa\right)\right).
\end{displaymath}
Then $\tilde\phi(x-ct)$ is a kink solution if space has a
circumference of $\tilde L$. This kink and an oppositely moving
antikink (but placed on the last quarter of space) made up the
initial field, so that $y_{i0}=\phi(40(i-1)/N)$,
$i=1,\ldots,N$, where
\begin{displaymath}
\phi(x)\equiv\left\{\begin{array}{ll}
\phi(x)&0\le x< L/4\\
2\pi&L/4\le x< 3L/4\\
2\pi-\phi(x-3L/4)&3L/4\le x < L
\end{array}\right.,\end{displaymath}
while $y_{i1}=y_{i0}+\dot\phi(40(i-1)/N)h$ where
\begin{displaymath}
\dot{\phi}(x)\equiv\left\{\begin{array}{ll}
(\phi(x-hc)-\phi(x))/h&0\le x< L/4\\
0&L/4\le x< 3L/4\\
-(\phi(x-hc)-\phi(x))/h&3L/4\le x < L
\end{array}\right..\end{displaymath}


\paragraph{Comparison with energy-conserving methods.}
As an example of how our method compares with an existing method, 
we considered the energy-conserving method  of
Vu-Quoc and Li [1993], page~354:
\begin{eqnarray}
&&\frac1{k^2}\biggl[\frac14\partial^2_ky_{i\,j+1}
+\frac12\partial^2_ky_{ij}
  +\frac14\partial^2_ky_{i\,j-1}\biggr]
  -\frac1{h^2}\partial^2_hy_{ij}\nonumber\\
&&\qquad\mbox{}+\frac12\left(
  \frac{N(y_{i\,j+1})-N(y_{ij})}{y_{i\,j+1}-y_{ij}}
  +\frac{N(y_{ij})-N(y_{i\,j-1})}{y_{ij}-y_{i\,j-1}}
  \right)=0.\label{155}
\end{eqnarray}
This 
has an iterative form similar to~(\ref{111}) and 
is quite comparable with~(\ref{101}) and~(\ref{102}) in
terms of the computation required. 
Our method  seems to preserve the soliton
waveform better than~(\ref{155}), as is indicated by comparison
of the left and right  Figure~(\ref{202}).

\begin{figure}[p]
\vspace*{-.25in}\centerline{
\epsfxsize=5.5in\epsfbox{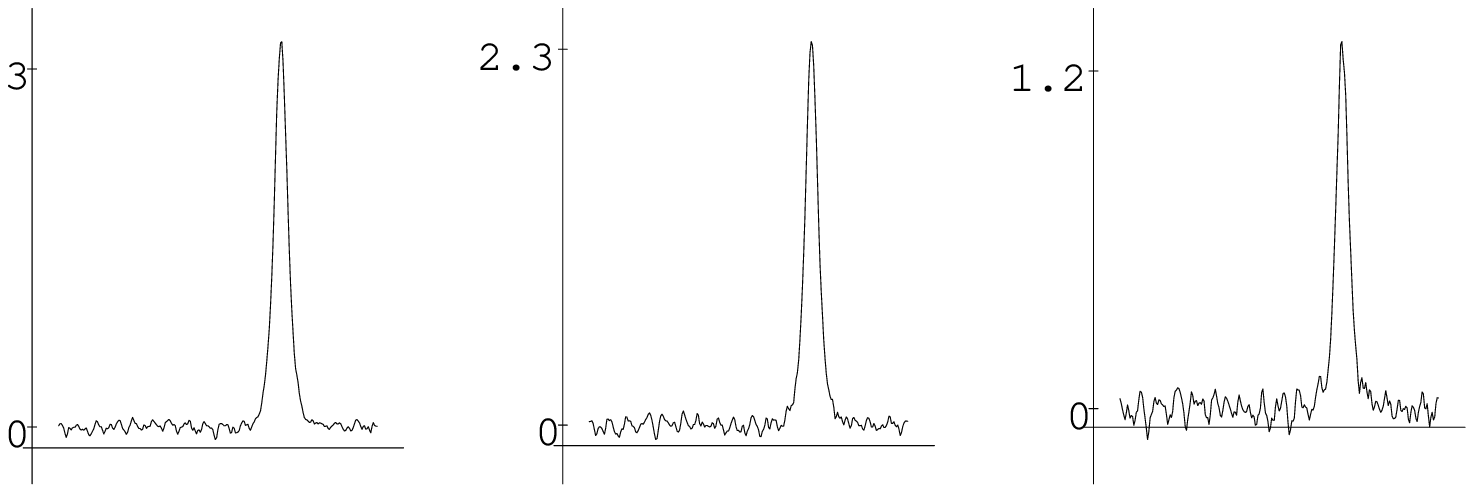}\vspace*{-.8in}}
\centerline{
\epsfxsize=5.5in\epsfbox{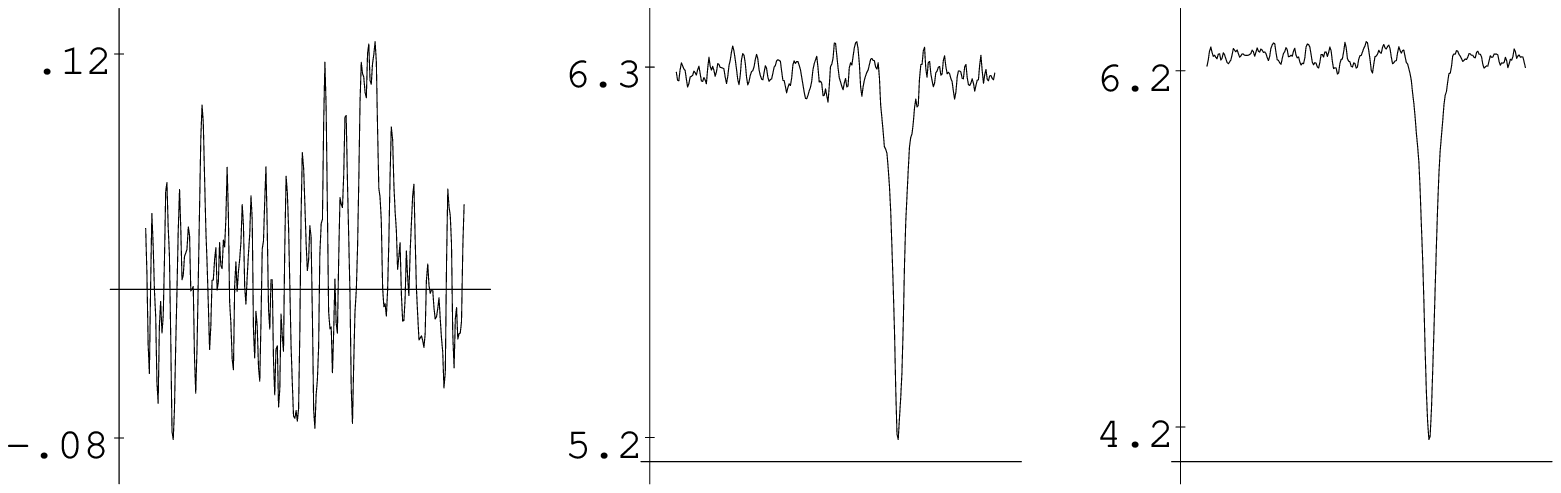}\vspace*{-.5in}}
\caption{\label{706} {\it \small The soliton collision
at time $t\approx129130$, after the energy-conserving method~(\ref{500})
of Guo et al.  has simulated about 5000 soliton collisions.
The solitons collide beginning at the 
top left and proceed to the top right, then to the bottom left, and finally
to the bottom right. The vertical scales are not constant and visually
exaggerate the high frequency oscillations, which are small on the scale 
$0$ to $2\pi$.
}\normalsize}
\end{figure}

In regards to the closely related papers  Vu-Quoc and~Li~[1993] 
and  Li~and~Vu-Quoc~[1995], we could not verify 
in our simulations that their method conserves energy, nor could we 
verify their {\em proof} that their
method conserves energy. So, as a further check, we implemented
the following   energy-conserving method
of  Guo, Pascual, Rodriguez, and Vazquez [1986]:
\begin{equation}
\partial^2_ky_{ij}-\partial^2_hy_{ij}
+\frac{N(y_{i\,j+1})-N(y_{i\,j-1})}{y_{i\,j+1}-y_{i\,j-1}},\label{500}
\end{equation}
which conserves the discrete energy
\begin{displaymath}
\sum_{j=1}^N\left(
\frac12\frac{(y_{i\,j+1}-y_{ij})(y_{ij}-y_{i\,j-1})}{h^2}
+\frac12\left(\frac{y_{i+1\,j}-y_{ij}}{k}\right)^2-N(y_{ij})\right).
\end{displaymath}
This method diverged after just~345~soliton collisions.
\begin{figure}[t]
\vspace*{-.25in}\centerline{
\epsfxsize=2.5in\epsfbox{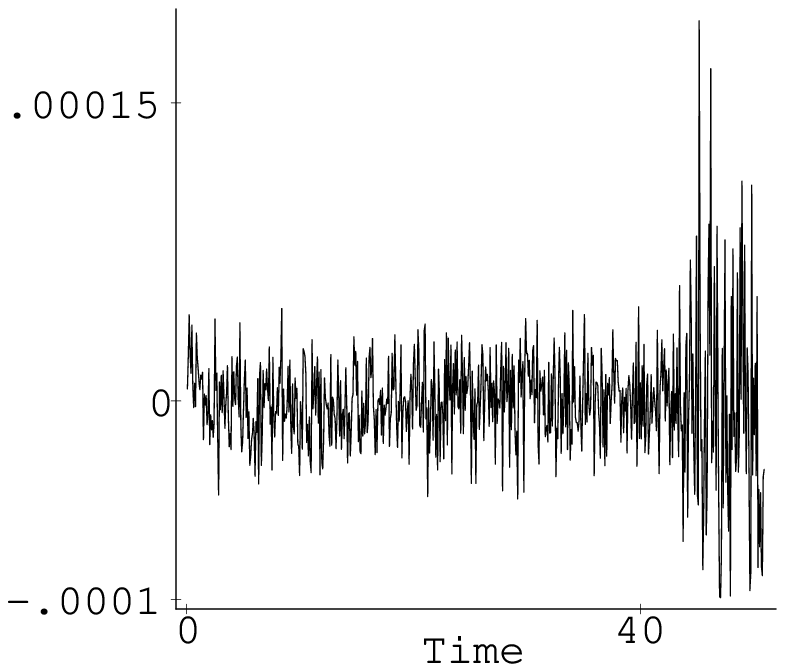}\hspace*{-.3in}
\epsfxsize=2.5in\epsfbox{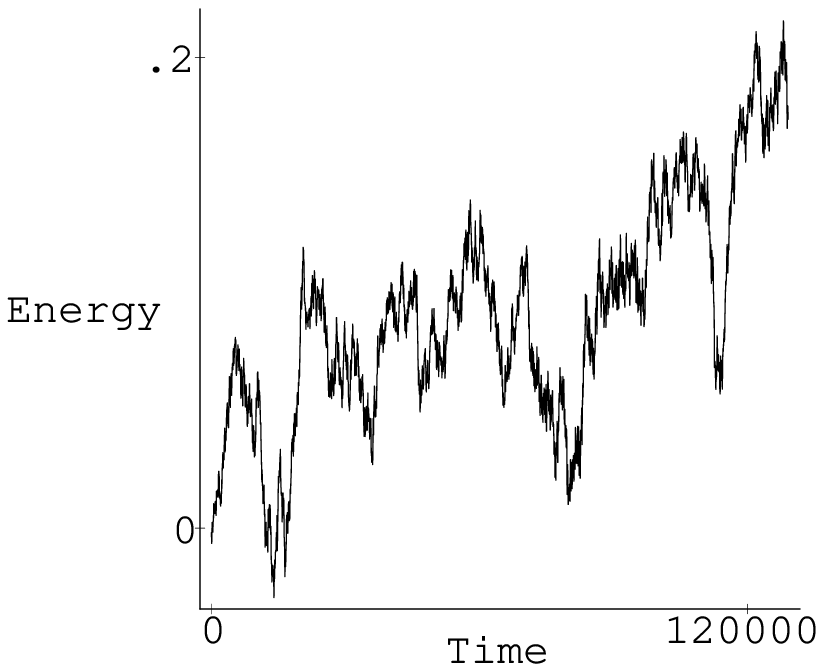}}
\vspace*{-.6in}\centerline{
\epsfxsize=2.5in\epsfbox{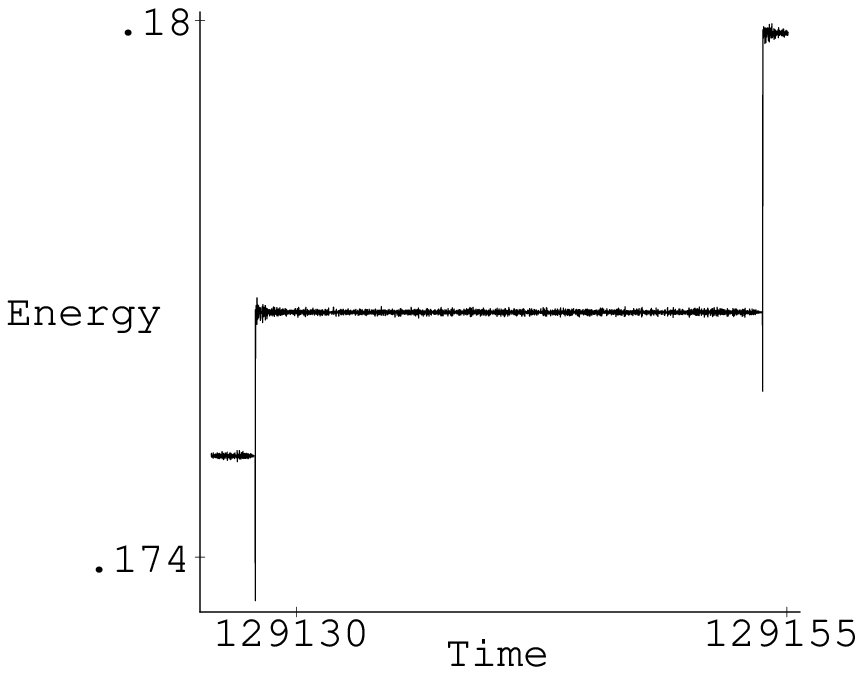}\hspace*{-.3in}
\epsfxsize=2.5in\epsfbox{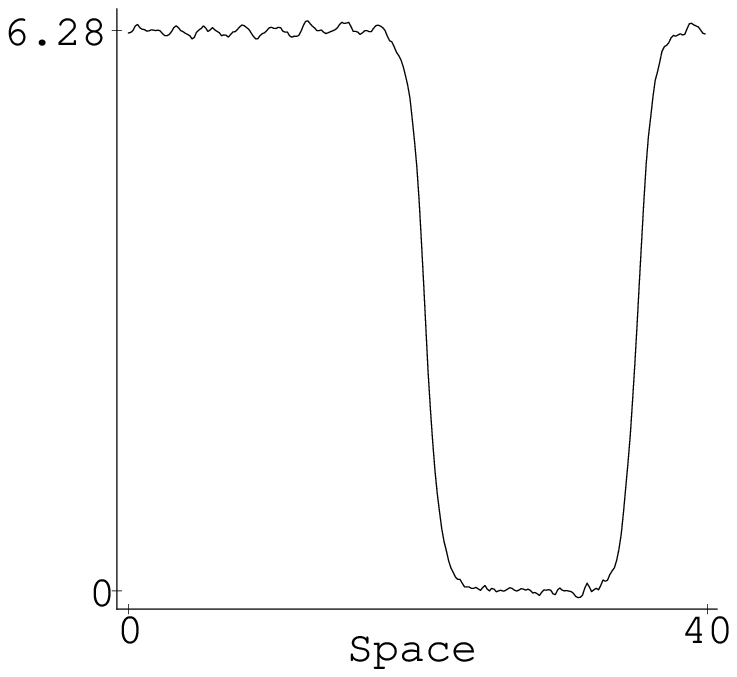}}\vspace*{-.2in}
\caption{\label{705} {\it \small A long-time simulation using the energy-conserving 
method~(\ref{500}) of Guo et al.
Above left: the initial energy error. Above right: the average energy error
over consecutive 5000 time step regions (the
maximum and minimum closely parallel the average).
Below left: the final energy error.
Below right: the final waveform at $t\approx129149$. 
}\normalsize}
\end{figure}
As can be seen from~(\ref{500}), the nonlinear potential $N$ enters
as a difference over two grid spacings, which suggests that 
halving   the time step might result in a more 
fair comparison with the methods~(\ref{102}) or~(\ref{155}).
With this advantage, method~(\ref{500}) was able to
simulate 5000~soliton collisions, with a waveform degradation
similar to the  energy-conserving method~(\ref{155}), as shown at
the bottom right of~Figure~(\ref{705}).   
The same figure also shows that, although the energy behavior 
of~(\ref{500}) is 
excellent for 
short time simulations,
it drifts significantly over long times,
and the final energy error has a  peculiar appearance.
Figure~(\ref{706}) shows the time
evolution of the waveform through the soliton collision that occurs
just before
the simulation stops. Apparently,  at the soliton collisions, significant 
high frequency
oscillations are present, and these are causing the jumps in the energy
error in the bottom left plot of Figure~(\ref{705}). This error then 
accumulates
due to the energy-conserving property of the method.
In these simulations, so as to guard against
the possibility that this behavior of the energy
was due to inadequately solving the implicit equation~(\ref{500}), 
we imposed a minimum limit of 3~iterations in the corresponding
iterative loop, whereas this loop would
otherwise have converged after just 1~iteration.

\begin{figure}[p]
\vspace*{-.25in}\centerline{
\epsfxsize=5.5in\epsfbox{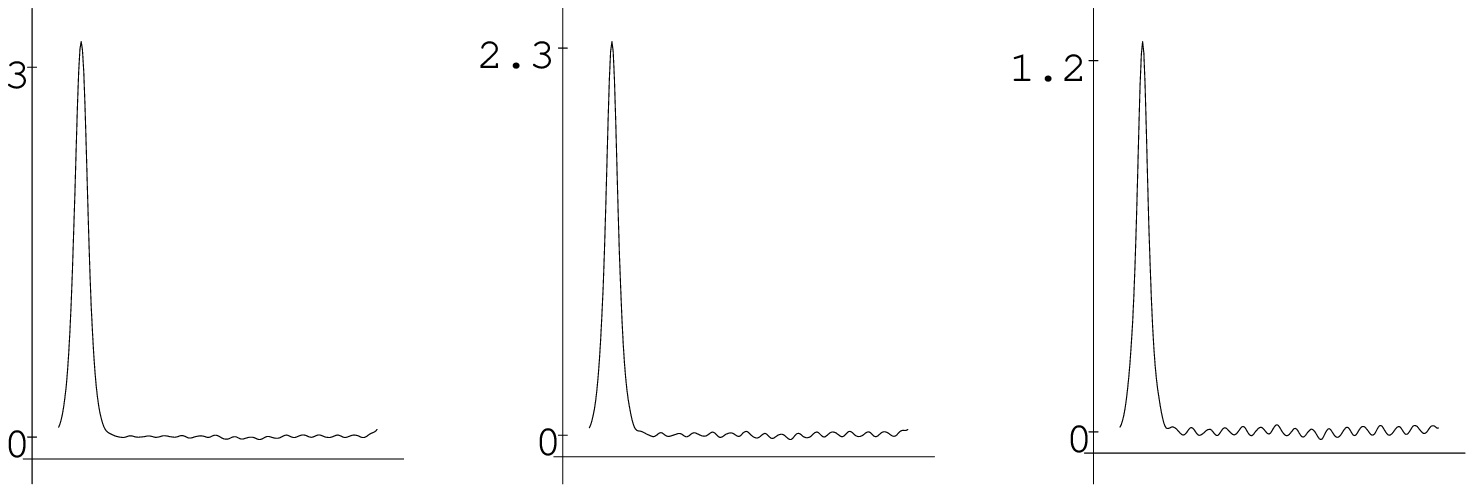}\vspace*{-.8in}}
\centerline{
\epsfxsize=5.5in\epsfbox{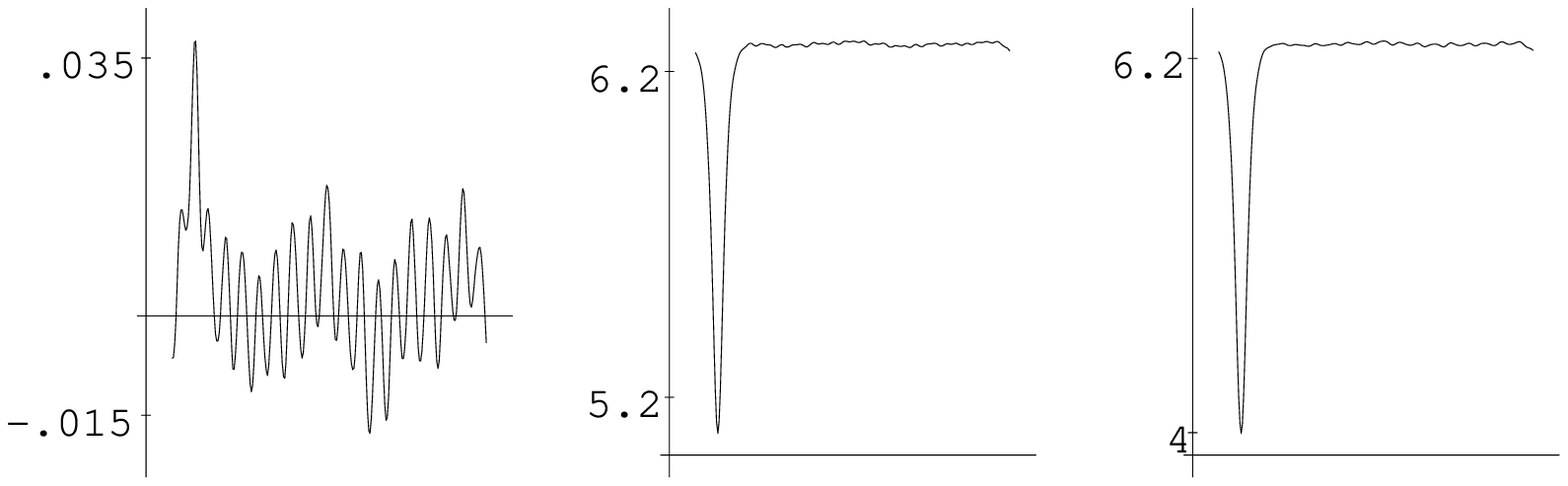}\vspace*{-.5in}}
\caption{\label{790} {\it \small Similar to  the above plot but for our
triangle-based multisymplectic method~(\ref{explicit}).
}\normalsize}
\end{figure}
\begin{figure}[ht]
\vspace*{-.25in}\centerline{
\epsfxsize=2.5in\epsfbox{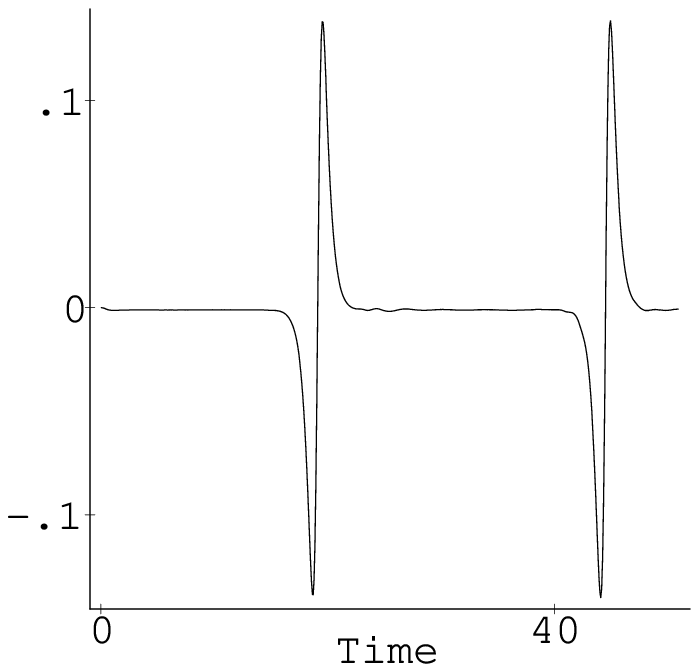}\hspace*{-.3in}
\epsfxsize=2.5in\epsfbox{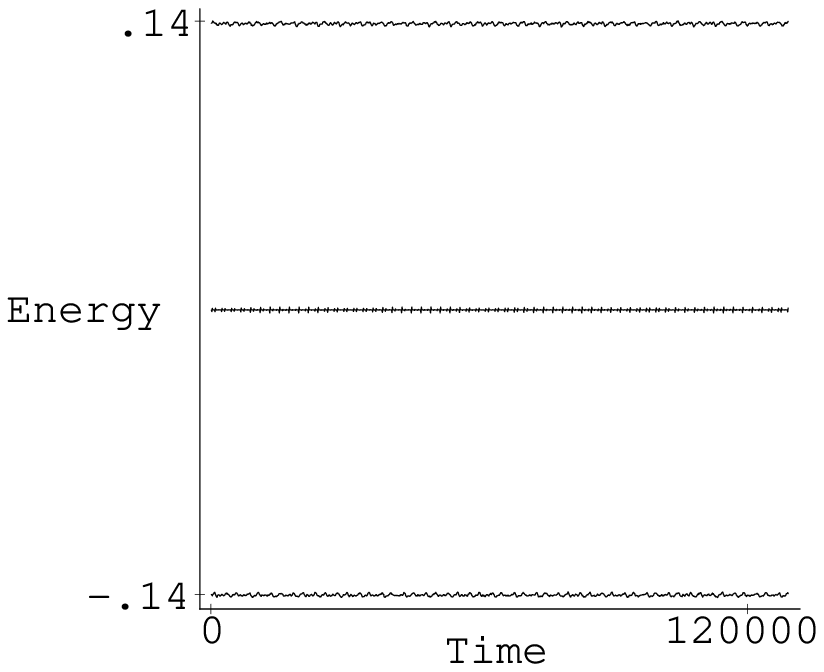}}
\vspace*{-.6in}\centerline{
\epsfxsize=2.5in\epsfbox{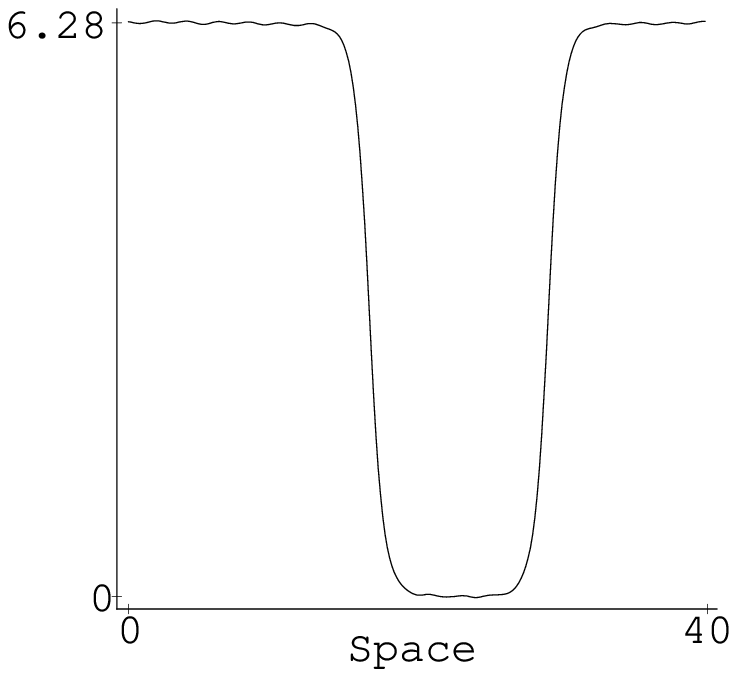}\hspace*{-.3in}
\epsfxsize=2.5in\epsfbox{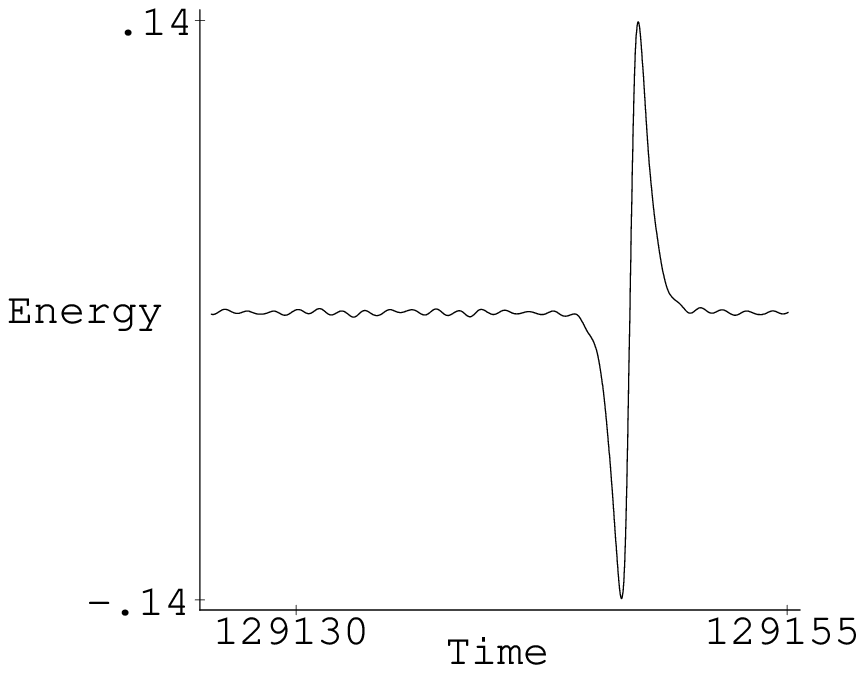}}\vspace*{-.2in}
\caption{\label{703} 
{\it \small A simulation  of 5000~soliton collisions using the 
triangle-based multisymplectic method~(\ref{explicit}).
Above left: The initial energy error. Above right: The minimum, average 
and maximum energy as in the left of
Figure~(\ref{203}).  Below left: the final waveform (at $t\approx 129130$). 
Below right:  the final energy error.}\normalsize}
\end{figure}

\paragraph{Comparison with the triangle-based multisymplectic method.}
The discrete second derivatives in the method~(\ref{500})
are the same as in the triangle-based multisymplectic method~(\ref{explicit});
these derivatives are simpler than either our  
rectangle-based multisymplectic method~(\ref{102})
or the  energy-conserving method of Vu-Quoc and Li~(\ref{155}).
To explore  this we  implemented the 
triangle-based multisymplectic method~(\ref{explicit}). 
Even with the less complicated discrete second derivatives
our triangle-based multisymplectic method
simulated 5000~soliton~collisions with comparable energy
\footnote{
The discrete energy that we calculated was
\begin{displaymath}
\sum_{i=1}^N\left(
\frac12\left(\frac{y_{i\,j+1}-y_{ij}}{h}\right)^2
+\frac12\left(\frac{y_{i+1\,j}-y_{ij}}{2k}\right)^2
-N(\bar y_{ij})\right).
\end{displaymath}}
and waveform
preservation properties as the rectangle-based multisymplectic
method~(\ref{102}), as shown in Figure~(\ref{703}).
Figure~(\ref{790}) shows the time
evolution of the waveform through the soliton collision just before
the simulation stops, and may be compared to Figure~(\ref{706}).
As can be seen, the high frequency oscillations that are
present during the soliton collisions are smaller and more smooth for
the triangle-based multisymplectic method than for the energy-conserving
method~(\ref{500}). A similar statement is true irrespective which of
the two multisymplectic or two energy conserving methods we tested, and
is true all along the waveform, irrespective of whether or not a
soliton collision is occurring.

\paragraph{Summary.}
Our multisymplectic methods
are finite difference methods that are  computationally 
competitive with existing   finite difference methods.
Our methods show promise for long-time simulations of conservative partial
differential equations, in that, for  long-time simulations of the
sine-Gordon equation, our method 1) had superior  energy-conserving 
behavior, {\em even when compared with energy-conserving methods};
2) better preserved the waveform than energy-conserving methods;
and 3) exhibited superior  stability, in that our methods
excited smaller and more smooth high frequency oscillations than
energy-conserving methods.  However, further
numerical investigation is certainly necessary to  make any
lasting conclusions about the long-time behavior of our
integrator.

\paragraph{The programs.}
The programs that were used in the preceding simulations are ``C''~language
implementations of the various methods. A simple tridiagonal LUD
method was used to solve  the linear equations (e.g. the left
side of~(\ref{111})), as in Vu-Quoc and~Li~[1993], page~379. 
An $8^{\mbox{\scriptsize th}}$~order  extrapolator was used to
provide a  seed for the implicit step. All calculations were
performed in double precision while the  implicit step was
terminated when the fields ceased to change to single precision;
the program's output was in single precision. 
The extrapolation usually provided a seed accurate enough so that the
methods became practically {\em explicit}, in that for many of the
time-steps  the first or second run through the iterative loops
solving the implicit equations
solved those equations to  single precision.
However,  in
the absence of a regular spacetime grid the expenses of the
extrapolation and solving the linear equation would grow. Our
programs  are freely available at URL
\verb|http://www.cds.caltech.edu/cds|.

\section{Concluding Remarks}

Here we make a few miscellaneous comments and remark on some
work planned for the future.

\paragraph{Lagrangian reduction.} As mentioned in the text, it
is useful to have a discrete counterpart to the Lagrangian
reduction of Marsden and Scheurle [1993a,b], Holm, Marsden and
Ratiu [1998a] and Cendra, Marsden and Ratiu [1998]. We sketch
briefly how this theory might proceed. This reduction can be
done for both the case of ``particle mechanics'' and for field
theory. 

For particle mechanics, the simplest case to start with is 
an invariant (say left)
Lagrangian on the tangent bundle of a Lie group: $L : TG
\rightarrow \mathbb{R}$. The reduced Lagrangian is
$l : \mathfrak{g} \rightarrow \mathbb{R}$ and the
corresponding Euler--Poincar\'e equations have a 
variational principle of La\-gran\-ge d'Al\-em\-bert type in
that there are constraints on the allowed variations. This
situation is described in Marsden and Ratiu [1994].

The discrete analogue of this would be to replace a discrete
La\-gran\-gian $ \mathbb{L} : G \times G \rightarrow \mathbb{R}
$ by a reduced discrete La\-gran\-gian $ \ell : G
\rightarrow \mathbb{R} $ related to $\mathbb{L}$ by
\[ 
\ell (g _1 g _2 ^{-1} ) = \mathbb{L} (g_1, g_2) 
\] 
In this situation, the algorithm from $G \times G$ to 
$G \times G$ reduces to one from $G$ to $G$ and it is generated
by $\ell$ in a way that is similar to that for $ \mathbb{L} $.
In addition, the discrete variational principle for $
\mathbb{L} $ which states that one should find
critical points of
\[ \mathbb{L} ( g_1 , g_2 ) + \mathbb{L} ( g_2 , g_3 )
\] 
with respect to $g_2$ to implicitly define the map
$( g_1 , g_2 ) \mapsto (g_2, g_3)$, reduces naturally to the
following principle:
Find critical points of
$\ell (g) + \ell (h)$ with respect to variations of 
$g$ and $h$ of the form $ g \xi : = L _g \xi $ and
$ \xi h = R _h \xi $ where $ L _g $ and $ R _h $ denote left
and right translation and where $ \xi \in \mathfrak{g} $. In
other words, one sets to zero, the derivative of 
the sum $\ell (g g _\epsilon^{-1} ) + \ell ( g _\epsilon h)$ with
respect to $\epsilon$ at $\epsilon = 0$ for a curve $
g_\epsilon $ in $G$ that passes through the identity
at $\epsilon = 0 $. This defines (with caveats of regularity as
before) a map of $G$ to itself, which is the reduced algorithm.
This algorithm can then be used to advance points in $G \times
G$ itself, by advancing each component equally, reproducing the
algorithm on $G \times G$. In addition, this can be used with
the adjoint or coadjoint action to advance points in
$\mathfrak{g}$ or $\mathfrak{g}^\ast$ to approximate the
Euler--Poincar\'e or Lie--Poisson dynamics. 

These equations for a discrete map, say $ \phi_{\ell} : G
\rightarrow G $  generated by $\ell$ on $G$ are called the {\bfi
discrete Euler--Poincar\'e equations} as they are the discrete
analogue of the Euler--Poincar\'e equations on $\mathfrak{g}$.
Notice that, at least in theory, computation can be done for this
map first and then the dynamics on $G \times G$ is easily
reconstructed by simply advancing each pair as follows: $(g_1,
g_2) \mapsto (hg_1, hg_2)$,  where $h = \phi_{\ell}(g_1^{-1} g_2)$.

If one identifies the discrete Lagrangians with generating
functions (as explained in Wendlandt and Marsden [1997]) then
the reduced La\-gran\-gian generates the reduced algorithm in
the sense of Ge and Marsden [1988], and this in turn is closely
related to the Lie--Poisson--Hamilton--Jacobi theory.

Next, consider the more general case of $TQ$ with its
discretization $Q\times Q$ with a group action (assumed to be free
and proper) by a Lie group $G$. The reduction of $TQ$ by the action
of $G$ is $TQ/G$, which is a bundle over $ T(Q/G)$ with fiber
isomorphic to $\mathfrak{g}$. The discrete analogue of this is
$(Q \times Q) / G$ which is a bundle over
$(Q/G) \times (Q/G)$ with fiber isomorphic to $G$ itself. The
projection map $\pi : (Q \times Q) / G \rightarrow (Q/G) \times
(Q/G)$ is given by $ [(q_1, q_2)] \mapsto ([q_1], [q_2])$ where $
[ \; ]$ denotes the relevant equivalence class. Notice that in the
case in which $ Q = G $ this bundle is ``all fiber''. The reduced
discrete Euler-Lagrange equations are similar to those in the
continuous case, in which one has shape equations couples with a
version of the discrete Euler--Poincar\'e equations. 

Of course all of the machinery in the continuous case can be
contemplated here too, such as stability theory, geometric phases,
etc. In addition, it would be useful to generalize this Lagrangian
reduction theory to the multisymplectic case. All of these topics
are planned for other papers.

\paragraph{Role of uniformity of the grid.} 
Consider an autonomous, continuous Lagrangian 
$\mathcal L:TQ\rightarrow\mathbb R$
where, for simplicity, $Q$ is an open submanifold of Euclidean space. 
Imagine
some {\em not necessarily uniform} temporal grid
($t_0,t_1,\cdots$) of $\mathbb R$,
so that $t_0<t_1<t_2<\cdots$. In this situation, it is natural to consider
the discrete action
\begin{equation}\label{800}
S=\sum_{k=1}^{n}L_k(q_{k},q_{k-1})\equiv
\sum_{k=1}^{n}\mathcal L\left(\frac{q_{k}+q_{k-1}}2,
  \frac{q_{k}-q_{k-1}}{t_{k}-t_{k-1}}\right)(t_{k}-t_{k-1}).
\end{equation}
This action principle  deviates from the action principle~(\ref{850})
of Section~3 in that the discrete Lagrangian density depends explicitly 
on $k$. Of course nonautonomous continuous
Lagrangians also yield $k$-dependent discrete Lagrangian densities,
irrespective of uniformity of the grid.
Thus, nonuniform temporal grids or nonautonomous Lagrangians
give rise to discrete Lagrangian densities which are
more general those those we have considered in Section~(3).
For field theories,
the Lagrangian in the action~(\ref{30})
depends on the spacetime variables already, through its explicit 
dependence on the triangle~$\Delta$. However, it is 
only in the context of a uniform grid that we
have  experimented numerically and only in that context that we have
discussed the 
significance of 
the discrete  multisymplectic form formula and the discrete Noether
theorem.  

Using~(\ref{800}) as an 
example, will now indicate
why the issue of grid uniformity may not be serious.
The DEL equations corresponding to the action~(\ref{800}) are
\begin{equation}\label{802}
\frac{\partial L_{k}}{\partial q_1}(q_{k},q_{k-1})+
\frac{\partial L_{k+1}}{\partial q_2}(q_{k+1},q_{k})=0,\quad k=1,2,,\cdots,
\end{equation}
and this gives evolution maps $F_{k+1,k}:Q\times Q\rightarrow Q\times Q$
defined so that
\begin{displaymath}
F_{k+1,k}(q_{k},q_{k-1})=(q_{k+1},q_{k}),\quad k=1,2,\cdots
\end{displaymath}
when~(\ref{802})
holds. For the canonical 1-forms corresponding to~(\ref{13}) and~(\ref{14})
we have the $k$-dependent one forms
\begin{equation}
 \theta_{L,k}^-(q_1,q_0)\cdot(\delta q_1,\delta q_0)\equiv
 \frac{\partial L_k}{\partial q_0}(q_1,q_0)\delta q_0,
\label{813}\end{equation}
and
\begin{equation}
  \theta_{L,k}^+(q_1,q_0)\cdot(\delta q_1,\delta q_0)\equiv
  \frac{\partial L_k}{\partial q_1}(q_1,q_0)\delta q_1,
\label{814}\end{equation}
and Equations~(\ref{16}) and~(\ref{2010}) become
\begin{equation}\label{806}
F_{k+1,k}^*(d\theta_{L,k}^+)=-d\theta_{L,k+1}^-,\quad 
d\theta_{L,k}^-+d\theta_{L,k}^+=0
\end{equation}
respectively. Together, these two equations give
\begin{equation}\label{807}
F_{k+1,k}^*(d\theta_{L,k}^+)=d\theta_{L,k+1}^+,
\end{equation}
and if we set 
\begin{displaymath}
F_k\equiv F_{k,k-1}\circ F_{k-1,k-2}\circ\cdots\circ F_{2,1}
\end{displaymath}
then ~(\ref{807}) chain together to imply
$F_k^*(d\theta_{L,1}^+)=d\theta_{L,k}^+$. 
This appears less than adequate since it merely says 
that the pull back by the evolution of a certain 2-form is, in general, a
different 2-form. The significant point to note, however, is that  
{\em this situation
may be repaired at any $k$  simply by choosing $L_k=L_1$.} It is easily
verified that the analogous statement is true with respect to momentum 
preservation via the discrete Noether theorem.

Specifically, imagine integrating a symmetric autonomous
mechanical system in a timestep adaptive way with
Equations~(\ref{802}). As the integration proceeds, various timesteps
are chosen, and if momentum is monitored it will show a dependence
on those choices. {\em A momentum-preserving symplectic simulation
may be obtained by simply choosing the last timestep to be of
equal duration to the first.} This is the highly desirable situation 
which gives us some confidence that grid uniformity is a nonissue.
There is one caveat: symplectic integration algorithms are evolutions
which are high frequency perturbations of the actual system, the frequency
being the inverse of the timestep, which is generally far smaller than
the time scale of any process in the simulation. However,  timestep
adaptation schemes will make choices on a much larger time scale than
the timestep itself, and then  drift in the energy will
appear on this larger time scale. 
A meaningful long-time simulation cannot be expected in the unfortunate 
case that the 
timestep adaptation makes
repeated choices in a way that resonates with some process of the
system being simulated.

\paragraph{The sphere.} The sphere cannot be generally
uniformly  subdivided into spherical triangles; however, a good
approximately uniform grid is obtained as follows: start  from
an inscribed icosahedron which produces a uniform subdivision 
into twenty spherical isosceles triangles; these are further
subdivided  by halving their sides and joining the resulting
points by short geodesics.

\paragraph{Elliptic PDEs.}
The variational approach we have developed allows us to examine
the multisymplectic structure of elliptic boundary value
problems as well. For a given Lagrangian, we form the
associated action function, and by computing its first
variation, we obtain the unique multisymplectic form of the
elliptic operator.  The multisymplectic form formula contains
information on how symplecticity interacts with spatial
boundaries.  In the case of two spatial dimensions, $X=\mathbb{R}^2$, $Y=\mathbb{R}^3$, we see that equation (\ref{t15})
gives us the conservation law 
$$ \text{div} {\mathcal X}=0,$$ where the vector ${\mathcal
X}=(\omega^0(j^1V,j^1W), \omega^1(j^1V, j^1W))$.

Furthermore, using our generalized Noether theory, we may define
momentum-mappings of the elliptic operator associated with its
symmetries.  It turns out that for important problems of spatial
complexity arising in, for example, pattern formation systems,
the covariant Noether current intrinsically contains the
constrained toral variational principles whose solutions are
the complex patterns (see Marsden and Shkoller [1997]).

There is an interesting connection between our variational
construction of multi\-symplectic\--momentum integrators and the
finite element method (FEM) for elliptic boundary value
problems.  FEM is also a variationally derived numerical
scheme, fundamentally differing from our approach in the
following way:  whereas we form a discrete action sum and
compute its first variation to obtain the discrete
Euler-Lagrange equations, in FEM, it is the original continuum
action function which is used together with a projection of the
fields and their variations onto appropriately chosen
finite-dimensional spaces.  One varies the projected fields and
integrates such variations over the spatial domain to recover
the discrete equations.  In general, the two discretization
schemes do not agree, but for certain classes of finite element
bases with  particular integral approximations, the resulting
discrete equations match the discrete Euler-Lagrange equations
obtained by our method, and are hence naturally multisymplectic.

To illustrate this concept, we consider the Gregory and Lin
method of solving two-point boundary value problems in optimal
control.  In this scheme, the discrete equations are obtained
using a finite element method with a basis of linear
interpolants.  Over each one-dimensional element, let $N_1$ and
$N_2$ be the two linear interpolating functions.  As usual, we
define the action function by $S(q)=\int_0^T
L(q(t),\dot{q}(t))dt$. Discretizing the interval $[0,T]$ into
$N$$+$$1$ uniform elements, we may write the action with fields
projected onto the linear basis as 
$$S(q) = \sum_{k=0}^{N-1} \int_k^{k+1} L( \{ N_1 \phi_k + N_2
\phi_{k+1}\},
\{ \dot{N}_1 \phi_k + \dot{N}_2 \phi_{k+1}\} )dt.$$ Since the
Euler-Lagrange equations are obtained by linearizing the action
and hence the Lagrangian, and as the functions $N_i$ are
linear, one may easily check that by evaluating the integrals
in the linearized equations  using a trapezoidal rule, the
discrete Euler-Lagrange equations given in (\ref{12}) are
obtained.  Thus, the Gregory and Lin method is actually a  
multisymplectic-momentum algorithm.

\paragraph{Applicability to fluid problems.} Fluid problems are
not literally covered by the theory presented here because
their symmetry groups (particle relabeling symmetries) are not
vertical. A generalization is needed to cover this case and we
propose to work out such a generalization in a future paper,
along with numerical implementation, especially for geophysical
fluid problems in which conservation laws such as conservation
of enstrophy and Kelvin theorems more generally are quite
important.
\paragraph{Other types of integrators.} It remains to link the
approaches here with other types of integrators, such as volume
preserving integrators (see, eg, Kang and Shang [1995], Quispel
[1995]) and reversible integrators (see, eg, Stoffer [1995]). In
particular  since volume manifolds may be regarded as
multisymplectic manifolds, it seems reasonable that there is an
interesting link.

\paragraph{Constraints.} One of the very nice things about the
Veselov construction is the way it handles constraints, both
theoretically and numerically (see Wendlandt and Marsden
[1997]). For field theories one would like to have a similar
theory. For example, it is interesting that for fluids, the
incompressibility constraint can be expressed as a
pointwise constraint on the first jet of the particle placement
field, namely that its Jacobian be unity. When viewed this way,
it appears as a holonomic constraint and it should be amenable
to the present approach. Under reduction by the particle
relabeling group, such a constraint of course becomes the
divergence free constraint and one would like to understand how
these constraints behave under both reduction and
discretization. 

\subsection*{Acknowledgments}
\addcontentsline{toc}{section}{Acknowledgments}
We would like to extend our gratitude to Darryl Holm, Tudor
Ratiu and Jeff Wendlandt for their time, encouragement and
invaluable input. Work of J. Marsden was supported by the
California Institute of Technology and NSF grant DMS
96--33161. Work by G. Patrick was partially supported by NSERC
grant OGP0105716 and that of S. Shkoller was partially
supported by the Cecil and Ida M. Green Foundation and DOE.
We also thank the Control and Dynamical Systems Department at
Caltech for providing a valuable setting for part of this
work.

\subsection*{References}
\addcontentsline{toc}{section}{References}
\begin{description}

\item Arms, J.M., J.E. Marsden, and V. Moncrief [1982]  
The structure of the space solutions of Einstein's equations:
II Several Killings fields and the  Einstein-Yang-Mills
equations.  {\it Ann. of Phys.} {\bf 144}, 81--106.

\item Ablowitz, M.J., Herbst, B.M., and C. Schober [1996] On
the numerical solution of the Sine-Gordon equation 1. Integrable
discretizations and homoclinic manifolds.  {\it J. Comp. Phys.}
{\bf 126}, 299--314.

\item Abraham, R. and J.E. Marsden [1978]
{\it Foundations of Mechanics.} Second Edition,
Addison-Wesley.

\item Arms, J.M., J.E. Marsden, and V. Moncrief [1982]  The
structure of the space solutions of Einstein's equations:
II Several Killings fields and the  Einstein-Yang-Mills
equations.  {\it Ann. of Phys.} {\bf 144}, 81--106.

\item Arnold, V.I. [1978]
{\it Mathematical Methods of Classical Mechanics.}
Graduate Texts in Math. {\bf 60}, Springer Verlag.
(Second Edition, 1989).

\item Ben-Yu, G, P. J. Pascual, M. J. Rodriguez, and L. Vazquez [1986]
Numerical solution of the sine-Gordon equation.
{\it Appl. Math. Comput.} {\bf 18}, 1--14

\item  Bridges, T.J. [1997]
Multi-symplectic structures and wave propagation.
Math. Proc. Camb. Phil. Soc., {\bf 121}, 147--190.

\item  Calvo, M.P. and E. Hairer [1995]
Accurate long-time integration of dynamical systems.
Appl. Numer. Math., {\bf 18}, 95.

\item  Cendra, H.,  J. E. Marsden and T.S. Ratiu [1998]
Lagrangian reduction by stages. {\it preprint.}

\item de Vogela\'{e}re, R. [1956] Methods of integration which 
preserve the contact transformation property of the Hamiltonian
equations. Department of Mathematics, University of Notre Dame,
Report No. {\bf 4}.

\item Dragt, A.J. and J.M. Finn [1979] 
Normal form for mirror machine Hamiltonians, {\it J. Math.
Phys.} {\bf 20}, 2649--2660.

\item Ebin, D. and J. Marsden [1970] Groups of diffeomorphisms
and the motion of an incompressible fluid. {\it Ann. of Math.}
{\bf 92}, 102-163.

\item Ge, Z. and J.E. Marsden [1988]
Lie-Poisson integrators and Lie-Poisson Hamilton-Jacobi theory.
{\it Phys. Lett. A} {\bf 133}, 134--139.

\item Gotay, M., J. Isenberg, and J.E. Marsden [1997]
Momentum Maps and the Hamiltonian Structure of
Classical Relativistic Field
Theories, I. {\it Preprint.}

\item Gregory, J and C. Lin [1991], The numerical solution of
variable endpoint problems in the calculus of variations,
Lecture Notes in Pure and Appl. Math., {\bf 127}, 175--183.

\item Guo, B. Y., P. J.  Pascual, M. J. Rodriguez, and L. Vazquez [1986]
Numerical solution of the sine-Gordon equation. 
{\it Appl. Math. Comput.} {\bf 18}, 1--14.

\item Holm, D. D., Marsden, J. E. and Ratiu, T. [1998a]  The
Euler-Poincar\'{e} equations and semidirect products with
applications to continuum theories. {\it Adv. in Math.} (to
appear).

\item Holm, D. D., J.E. Marsden and T. Ratiu [1998b] The
Euler-Poincar\'{e} equations in geophysical fluid dynamics, in 
{\it Proceedings of the Isaac Newton Institute Programme on the
Mathematics of Atmospheric and Ocean Dynamics}, Cambridge
University Press (to appear).

\item Li, S. and L., Vu-Quoc [1995] Finite-difference calculus
invariant structure of a class of algorithms for the nonlinear
Klein-Gordon equation. {\it SIAM J. Num. Anal.} {\bf 32},
1839--1875.

\item Marsden, J.E., G.W. Patrick, and W.F. Shadwick (Eds.)
[1996] {\it Integration Algorithms and Classical Mechanics.}
Fields Institute Communications {\bf 10}, Am. Math. Society.

\item Marsden, J.E. and T.S. Ratiu [1994]  {\it Introduction to
Mechanics and Symmetry.} Texts in Applied Mathematics {\bf  17},
Springer-Verlag.

\item Marsden, J.E. and S. Shkoller [1997]  Multisymplectic
geometry, covariant Hamiltonians and water waves. To appear in
{\it Math. Proc. Camb. Phil. Soc.} {\bf 124}.

\item Marsden, J.E. and J.M. Wendlandt [1997] Mechanical systems
with symmetry, variational principles and integration algorithms.
{\it Current and Future Directions in Applied Mathematics}, 
Edited by M. Alber, B. Hu, and J. Rosenthal, Birkh\"{a}user, 219--261.

\item McLachlan, R.I. and C. Scovel [1996]
A survey of open problems in symplectic integration.
{\it Fields Institute Communications\/}
{\bf 10}, 151--180.

\item Moser, J. and A.P. Veselov [1991]
Discrete versions of some classical integrable systems and
factorization of matrix polynomials.
{\it Comm. Math. Phys.} {\bf 139}, 217--243.

\item Neishtadt, A. [1984] The separation of motions in systems
with rapidly rotating phase. {\it P.M.M. USSR} {\bf 48},
133--139.

\item Palais, R.S. [1968] Foundations of global nonlinear analysis,
Benjamin, New York.

\item Palais, R.S. [1997] The symmetries of solitons.
{\it Bull. Amer. Math. Soc.} {\bf 34}, 339--403.

\item Quispel, G.R.W. [1995] Volume-preserving integrators. {\it
Phys. Lett. A.} {\bf 206}, 26.

\item Sanz-Serna, J. M. and M. Calvo [1994]
  {\it Numerical Hamiltonian Problems}.
  Chapman and Hall, London.

\item Simo, J.C.and O. Gonzalez. [1993],  Assessment of
Energy-Momentum and Symplectic Schemes for Stiff Dynamical
Systems. {\it Proc. ASME Winter Annual
Meeting, New Orleans}, Dec. 1993. 

\item Stoffer, D. [1995] Variable steps for reversible integration
methods, {\it Computing} {\bf 55}, 1.

\item Veselov, A.P. [1988]
Integrable discrete-time systems and difference operators.
{\it  Funkts. Anal. Prilozhen.} {\bf 22}, 1--13.

\item Veselov, A.P. [1991]
Integrable {Lagrangian} correspondences and the factorization of
matrix polynomials.
{\it Funkts. Anal. Prilozhen.} {\bf 25}, 38--49.

\item Vu-Quoc, L. and S. Li [1993] Invariant-conserving finite
difference algorithms for the nonlinear  Klein-Gordon equation.
{\it Comput. Methods Appl. Mech. Engrg.} {\bf 107}, 341--391.

\item Wald, R.M. [1993] 
Variational principles, local symmetries and
black hole entropy. {\it Proc. Lanczos Centennary volume\/} SIAM,
231--237.

\item Wendlandt, J.M. and J.E. Marsden [1997] Mechanical
integrators derived from a discrete variational principle. {\it
Physica D} {\bf 106}, 223--246.
\end{description}
\end{document}